\documentclass[final,1p,times]{elsarticle}
\usepackage{amssymb,latexsym,amsmath}
\usepackage{color,caption,cases,bm}
\usepackage[colorlinks,linkcolor=blue, anchorcolor=blue, citecolor=blue]{hyperref}
\usepackage{graphicx,graphics,booktabs}
\graphicspath{{./Figures/}}
\usepackage{tabularx}
\usepackage{stmaryrd}
\usepackage[T1]{fontenc}
\usepackage{multirow}
\biboptions{numbers,sort&compress} 
\numberwithin{equation}{section}

\allowdisplaybreaks
\newtheorem{thm}{Theorem}[section]
\newtheorem{lem}[thm]{Lemma}

\newtheorem{rem}[thm]{Remark}
\newtheorem{example}[thm]{Example}

\newcommand{\mo}{\mathcal{O}}

\journal{XXX}
\begin{document}

%=================== Text begin here =============================================
\begin{frontmatter}

\title{Maximum-Norm Error Estimates of Fourth-Order Compact and ADI Compact Finite Difference Methods for Nonlinear Coupled Bacterial Systems }

\author[OUC]{Jie Xu} \ead{jxu129@163.com}
\author[OUC,LMS]{Shusen Xie} \ead{shusenxie@ouc.edu.cn}
\author[OUC,LMS]{Hongfei Fu\corref{Fu}}\ead{fhf@ouc.edu.cn}

\address[OUC]{School of Mathematical Sciences, Ocean University of China, Qingdao  266100, China}
\address[LMS]{Laboratory of Marine Mathematics, Ocean University of China, Qingdao, Shandong 266003, China}
\cortext[Fu]{Corresponding author.}

\begin{abstract}
In this paper, by introducing two temporal-derivative-dependent auxiliary variables, a linearized and decoupled fourth-order compact finite difference method is developed and analyzed for the nonlinear coupled bacterial systems. The temporal-spatial error splitting technique and discrete energy method are employed to prove the unconditional stability and convergence of the method in discrete maximum norm. Furthermore, to improve the computational efficiency, an alternating direction implicit (ADI) compact difference algorithm is  proposed, and the unconditional stability and optimal-order maximum-norm error estimate for the ADI scheme are also strictly established. Finally, several numerical experiments are conducted to validate the theoretical convergence and to simulate the phenomena of bacterial extinction as well as the formation of endemic diseases.
\end{abstract}


\end{frontmatter}

%%%%%%%%%%%%%%%%%%%%%%%%%%%%%%%%%%%%%%%%%%%%%%%%%%%%%%%%%%%%%%%%%%%%%%

\section{Introduction}
In 1979, Capasso et al. proposed an ordinary differential equation (ODE) system to model the 1973 cholera epidemic that spread across the European Mediterranean regions \cite{CP}, which is given by
\begin{equation}\label{bact:ode}
	\begin{cases}
	\dfrac{d z^1}{d t}=-a_{11} z^1+a_{12} z^2, \\
	\dfrac{d z^2}{d t}=-a_{22} z^2+g(z^1),
	\end{cases}
\end{equation}
supplemented with appropriate initial conditions. Here, $z^1$ represents the average concentration of bacteria, and $z^2$ represents the infective human population in an urban community. The term $-a_{11} z^1$ describes the natural growth rate of the bacterial population, while $a_{12} z^2$ represents the contribution of the infective humans to the growth rate of bacteria. In the second equation, the term $-a_{22} z^2$ describes the natural damping of the infective population due to the finite mean duration of the infectiousness of humans. The last term $g(z^1)$ is the infection rate of humans under the assumption that the total susceptible human population remains constant during the epidemic. This type of mechanism is also suitable for interpreting other epidemics with oro-faecal transmission, such as typhoid fever, infectious hepatitis, schistosomiasis, and poliomyelitis, with appropriate modifications \cite{C, NH}. 

In order to make the model more realistic, people assume that bacteria diffuse randomly in the habitat due to the particular habits in the regions where these kinds of epidemics typically spread. In this context, the above ODE system is modified as the spatial propagation model of bacteria under human environmental pollution \cite{CVM, CVA}:
\begin{equation}\label{bact:eq}
	\begin{cases}
		u_{t}=d_{1} \Delta u-a_{11} u+a_{12} v, &\qquad  (\bm x, t) \in \Omega \times J, \\
		v_{t}=d_{2} \Delta v-a_{22} v+g(u), &\qquad  (\bm x, t) \in \Omega \times J, \\ 
		 u(\bm x, t)= v(\bm x, t)=0, &\qquad  (\bm x, t) \in \partial \Omega \times J, \\ 
		u(\bm x, 0)=u_{0}(\bm x), \quad v(\bm x, 0)=v_{0}(\bm x), & \qquad \bm x \in \Omega,
	\end{cases}
\end{equation} %$\bm x = (x, y)$, 
where $\Omega = (0,1)^{2}$ with boundary denoted by $\partial \Omega$, $J=(0, T]$ with $T>0$, $a_{11}$, $a_{22}$, and $a_{12}$ are real, positive constants. The variables $u$ and $v$ respectively denote the average concentrations of bacteria and infective people. 
The natural mortality and transmission rate  are expressed as $-a_{11} u$ and $a_{12} v$, respectively. 
The diffusion coefficients $d_{1}$ and $d_{2}$ are assumed to be greater than or equal to zero. 
While random diffusion in the human population may be negligible compared to that in the bacteria, it will generally be allowed for completeness, which means $d_1$, $d_2$ are greater than zero.
The cure rate of infected people is expressed as $-a_{22} v$, while the rate of population infection is indicated as $g(u) \in C^2(\mathbb{R})$, satisfying the Lipschitz continuous condition $\left| g(u_1)-g(u_2)\right|  \le L \left| u_1-u_2\right|$ for $u_1, u_2\in \mathbb{R}$ and $g(0)=0$. 

Several studies have focused on the theoretical analysis and numerical computation of the bacterial system \eqref{bact:eq}. For example, in \cite{SL} and \cite{SL22}, the finite element method and the $H^1$-Galerkin finite element method, combined with the two-grid method were respectively proposed. 
The authors proved the existence and uniqueness of solutions of these fully discrete schemes, and also derived second-order accurate spatial superclose and superconvergence estimates with respect to the $H^1$-norm. 
Additionally, in \cite{CJ}, two finite element schemes were considered for models with time delay, and optimal second-order error estimates in the $L^2$-norm were proved.
However, up to now, there has been no consideration of a high-order compact (HOC)  difference method for model \eqref{bact:eq}, and theoretical analyses of finite difference methods are also scarce.

Over the past few decades, HOC difference methods have attracted extensive research interests \cite{CFL, LSS, BGM, HPW, SXL, ZLL, WXF}, as the compact finite difference method can reach high-order accuracy on a compact stencil with few grid points.
Moreover, in practical computations, it is preferable to measure errors using error estimates in the grid-independent maximum-norm in numerical analysis of the proposed methods. As is well known, the Sobolev embedding inequality \cite{A} implies that $H^2(\Omega) \hookrightarrow C(\bar{\Omega})$ for a bounded convex domain $\Omega\subset \mathbb{R}^d$ $(d=2,3)$, where $\bar{\Omega}$ is the closure of the domain $\Omega$.
Thus, it is necessary to employ a technique of $H^2$-error analysis at the discretization level to obtain a maximum-norm error estimate, and yet many researchers have studied such error analysis to obtain the maximum-norm error estimate for various schemes (see, \cite{ ZL, LS, LSS, RC, CC, CFL}).
In particular, Zhao and Li \cite{ZL} developed a fourth-order compact difference scheme for the one-dimensional Benjamin-Bona-Mahony-Burgers’ equation. They studied the conservative invariant, boundedness and unique solvability of the scheme and also its convergence in the sense of  maximum-norm. 
In \cite{CFL}, Cai et al. proposed a fourth-order compact finite difference scheme for solving the quantum Zakharov system, and  optimal-order error estimate in the maximum-norm was proved without any requirement on the temporal stepsize ratio. 
The main purposes of this paper are to construct a  linearized and decoupled Crank-Nicolson (CN) type HOC difference scheme for the nonlinear coupled bacterial systems \eqref{bact:eq} and to analyze its  unconditional maximum-norm error estimates. To be specifically, inspired by the idea of \cite{LSS},  this paper introduces two auxiliary temporal-derivative-dependent variables, $p = u_t + a_{12}v$ and $q = v_t$, to reformulate the original problem \eqref{bact:eq} into a four-variable coupled system of equations. Then, a fourth-order compact difference method  is easily developed.
The auxiliary variables  lead to a splitting of temporal-spatial errors, allowing for the proof of corresponding unconditional stability and optimal-order error estimates in the maximum-norm.

It is well known that computational costs and CPU time consumption are of great consideration when numerical discretizations are used to model and simulate multi-dimensional PDEs. 
The alternating direction implicit (ADI) method, which reduces the solution of a multi-dimensional problem into a series of one-dimensional subproblems, has shown powerful abilities in solving multi-dimensional parabolic and hyperbolic problems (see Refs. \cite{LS, FZLZ, ZFH, WCLD, XZ, GL, DZ, EJS} and references therein). 
However, to the best of our knowledge, there has been no previous research on the construction of an ADI scheme for the bacterial systems \eqref{bact:eq}, let alone its maximum-norm error analysis.
In this paper, a ADI compact difference method for model \eqref{bact:eq} is developed, and the aforementioned error splitting technique is applied to the theoretical analysis of the ADI compact scheme, where the unconditional stability and convergence of the ADI compact scheme are rigorously proved under the discrete maximum norm. 

The main contribution of this paper is the construction of linearized and decoupled fourth-order compact and ADI compact finite difference schemes, along with their rigorous maximum-norm error analyses with second-order accuracy in time and fourth-order accuracy in space. This achievement is accomplished by utilizing a temporal-spatial error splitting technique and the discrete energy method. To the best of our knowledge, this is the first time that linearized and decoupled compact and ADI compact schemes for the nonlinear coupled bacterial systems have been constructed and rigorously analyzed under maximum-norm error.

The contents of the paper are organized as follows. Firstly, some  notations and auxiliary lemmas are presented in the next section. 
Section 3 is devoted to the construction and theoretical investigation of a fourth-order compact difference scheme. 
In Section 4, to achieve highly efficient computation, an ADI algorithm is further introduced, and its analysis under the discrete maximum-norm is carried out. 
In Section 5, several numerical experiments are presented to support our theoretical results. 
Finally, some comments are given in the concluding section. Throughout this paper, we use $C$ to denote a generic positive constant that may depend on the given data but is independent of the mesh parameters.

%%%%%%%%%%%%%%%%%%%%%%%%%%%%%%%%%%%%%%%%%%%%%%%%%%%%%%%%%%%%%%%%%%%%%%

\section{Some notations and auxiliary lemmas}
 First, for positive integers $N$, $M_x$ and $M_y$, let $t_n=n \tau$ $(0 \leq n \leq N)$, $x_i=i h_x$ $(0 \leq i \leq M_x)$ and $y_j= j h_y$ $(0 \leq j \leq M_y)$, with temporal stepsize $\tau=T / N$ and spatial stepsizes $h_x=1/ M_x$ and $h_y=1/ M_y$.  Let $h=\max  \{h_x, h_y \} $ and $t_{n-1/2}=\frac{1}{2} (t_{n-1}+t_{n})$. Denote the sets of spatial grids $\bar{\omega}=\left\{(i, j) \mid 0 \leq i \leq M_x, 0 \leq j \leq M_y\right\}$, $\omega=\left\{(i, j) \mid 1 \leq i \leq M_x-1, 1 \leq j \leq M_y-1\right\}$ and $\gamma=\bar{\omega} \setminus \omega$.
 Given temporal discrete grid functions $\left\{w^n=w(t_n) \mid 0 \leq n \leq N \right\}$ and $\{w^{n-1/2}=w(t_{n-1/2}) \mid 1 \leq n \leq N \}$,  denote 
\begin{align*}
 \bar{w}^{n-1/2}=\frac{1}{2}\left(w^n+w^{n-1}\right) ,
 \quad \delta_t w^{n-1/2}=\frac{1}{\tau}\left(w^n-w^{n-1}\right),
 \quad \triangle_t w^{n}=\frac{1}{\tau}\left(w^{n+1/2}-w^{n-1/2}\right).
\end{align*}
Besides, given spatial discrete grid function $v=\left\{v_{i,j} \mid (i, j)\in \bar{\omega}\right\}$, denote the difference operators
\begin{align*}
	&\delta_x v_{i-1/2, j}=\frac{1}{h_x}\left(v_{i,j}-v_{i-1, j}\right), \quad \delta_x^2 v_{i,j}=\frac{1}{h_x}\left(\delta_x v_{i+1/2, j}-\delta_x v_{i-1/2, j}\right).
\end{align*}
Similarly, one can define $\delta_y v_{i, j-1/2}$, $\delta_y^2 v_{i,j}$, $\delta_x \delta_y v_{i-1/2, j-1/2}$, $\delta_y \delta_x^2 v_{i, j-1/2}$,  $\delta_x \delta_y^2 v_{i-1/2, j}$, and $\delta_x^2 \delta_y^2 v_{i,j}$. Moreover, we also introduce  the following compact difference operators:
\begin{equation}
\begin{aligned}
	&\mathcal{H}_{x} v_{i,j}=\left\{\begin{array}{cl}
		\left(I+\frac{h_x^2}{12} \delta_x^2\right) v_{i,j}, & 1 \leq i \leq M_x-1, \quad 0 \leq j \leq M_y, \\
		v_{i,j}, & i=0 \text { or } M_x,\quad0 \leq j \leq M_y,
	\end{array}\right. 
\end{aligned}
\end{equation}
\begin{equation}
	\begin{aligned}
	& \mathcal{H}_{y} v_{i,j}=\left\{\begin{array}{cl}
		\left(I+\frac{h_y^2}{12} \delta_y^2\right) v_{i,j}, & 0 \leq i \leq M_x, \quad 1 \leq j \leq M_y-1, \\
		v_{i,j}, & 0 \leq i \leq M_x, \quad j=0 \text { or } M_y.
	\end{array}\right. 
\end{aligned}
\end{equation}
Denote
$
\mathcal{H}_h =\mathcal{H}_{x} \mathcal{H}_{y}$ and $ \Lambda_h = \mathcal{H}_{y} \delta_x^2+\mathcal{H}_{x} \delta_y^2.  
$

Next, define the spaces of  grid functions  
\[
\mathcal{V}_h=\left\{v= \{v_{i,j}  \} \mid  (i, j) \in \bar{\omega}\right\},
~~
\mathcal{V}_h^0=\left\{v\in\mathcal{V}_h \mid v_{i,j}=0\ \text{if}\ (i, j) \in \gamma\right\}.
\]
Then, for grid functions $v, w \in \mathcal{V}_h$, we introduce the discrete $L^2$ inner product and norm
\begin{align*}
\left( v, w\right)=h_x h_y \sum_{i=1}^{M_x-1} \sum_{j=1}^{M_y-1} v_{i,j} w_{i,j},\quad &\|v\|=\sqrt{\left(  v, v\right)}.
\end{align*}
Besides, the following discrete  $H^1$ and $H^2$ semi-norms and  norms are also used:
\begin{align*}
	|v|_1=\sqrt{ \|\delta_x v \|_x^2+ \|\delta_y v \|_y^2},\quad 
	& \|v\|_1=\sqrt{\|v\|^2+|v|_1^2}, \\
	\left|\delta_x v\right|_1=\sqrt{ \|\delta_x^2 v \|^2+ \|\delta_x \delta_y v \|_{xy}^2}, \quad 
	&  |\delta_y v |_1=\sqrt{ \|\delta_y \delta_x v \|_{yx}^2+ \|\delta_y^2 v \|^2}, \\
	|v|_2=\sqrt{ |\delta_x v |_1^2+ |\delta_y v |_1^2},\quad  
	& \|v\|_2=\sqrt{\|v\|_1^2+|v|_2^2},
\end{align*}
where
\begin{align*}
	  \|\delta_x v \|_x=\sqrt{h_x h_y \sum_{i=1}^{M_x} \sum_{j=1}^{M_y-1} \left|  \delta_x v_{i-1/2,j}\right|^2}, \quad
	& \|\delta_y v \|_y=\sqrt{h_x h_y \sum_{i=1}^{M_x-1} \sum_{j=1}^{M_y} \left|  \delta_y v_{i,j-1/2}\right|^2},\\
	 \|\delta_x \delta_y v \|_{xy}=\sqrt{h_x h_y \sum_{i=1}^{M_x} \sum_{j=1}^{M_y} \left|  \delta_x \delta_y v_{i-1/2, j-1/2}\right|^2},\quad
	& \|\delta_x^2 v \|=\sqrt{h_x h_y \sum_{i=1}^{M_x-1} \sum_{j=1}^{M_y-1} \left|  \delta_x^2 v_{i,j}\right|^2},
\end{align*}	
 \[	
\]
and  $ \|\delta_y\delta_x  v \|_{yx}$, $ \|\delta_y^2 v \|$ can  be defined similarly, in particular, $	 \|\delta_x \delta_y v \|_{xy} =  \|\delta_y\delta_x  v \|_{yx}$.
Moreover, we denote  the discrete maximum-norm
\begin{align*}
\|v\|_{\infty}=\max _{\left(i, j\right) \in \omega} |v_{i,j}|.
\end{align*}

 %%%%%%%%%%%%%%%%%%%%%%%
\begin{lem}[\cite{SZG}]\label{lem:trunerr:t1}
Let $f(t) \in C^1([t_{0}, t_{1/2}])$,   then we have
\begin{align*}
	f(t_{1/2}) = f(t_0) +\frac{\tau}{2} \int_0^1 f'\left(t_{1/2}-\frac{s \tau}{2}\right)\mathrm{d} s.
\end{align*}
\end{lem}
 %%%%%%%%%%%%%%%%%%%%%%%
\begin{lem}[\cite{SZG}]\label{lem:trunerr:t2}
Let $f(t) \in C^2\left(\left[t_{n-1}, t_n\right]\right)$,   then we have
\begin{align*}
		f(t_{n-1/2})
		&= \frac{f(t_{n}) + f(t_{n-1}) }{2}
		+\frac{\tau^2}{8} \int_0^1\left[f''\left( t_{n-1/2}-\frac{s \tau}{2}\right) +f''\left( t_{n-1/2}+\frac{s \tau}{2}\right) \right] (1-s) \mathrm{d}s.
\end{align*}
Besides, if $f(t) \in C^2\left(\left[t_{n-2}, t_n\right]\right)$, then we have
\begin{align*}
	f(t_{n-1/2})
	= \frac{3}{2}f(t_{n-1}) - \frac{1}{2} f(t_{n-2}) 
	+\frac{9\tau^2}{8} \int_0^1\left[f''\left( t_{n-1/2}-\frac{3s \tau}{2}\right) -3f''\left( t_{n-1/2}-\frac{s \tau}{2}\right) \right]  (1-s) \mathrm{d}s .
\end{align*}
\end{lem}
%%%%%%%%%%%%%%%%%%%%%%%
\begin{lem}[\cite{SZG}]\label{lem:trunerr:t3}
Let $f(t) \in C^3\left(\left[t_{n-1}, t_n\right]\right)$,  then  we have
\begin{align*}
	&f^{\prime}(t_{n-1/2})
	= \frac{f\left(t_n\right)-f\left(t_{n-1}\right)}{\tau} 
	+\frac{\tau^2}{16} \int_0^1\left[f^{(3)}\left( t_{n-1/2}-\frac{s \tau}{2}\right) +f^{(3)}\left( t_{n-1/2}+\frac{s \tau}{2}\right) \right] (1-s)^2 \mathrm{d}s .
\end{align*}
Moreover, if $f(t) \in C^3\left(\left[t_{n-1/2}, t_{n+1/2}\right]\right)$, then we have
\begin{align*}
	&f^{\prime}(t_{n})
	= \frac{f\left(t_{n+1/2}\right)-f\left(t_{n-1/2}\right)}{\tau} 
	+\frac{\tau^2}{16} \int_0^1\left[f^{(3)}\left( t_{n}-\frac{s \tau}{2}\right) +f^{(3)}\left( t_{n}+\frac{s \tau}{2}\right) \right] (1-s)^2 \mathrm{d}s .
\end{align*}
\end{lem}
%%%%%%%%%%%%%%%%%%%%%%%
\begin{lem}[\cite{BH}]\label{lem:trunerr:t4}
 For  $w(\bm x) \in H^6(\Omega)$, let
%\begin{align*}
	$\Xi_{i,j}=\mathcal{H}_{h} \Delta w(x_i,y_j) -\Lambda_h w(x_i,y_j)$  for  $ (i, j) \in  \omega$,
%\end{align*}
then we have
\begin{align*}
\left\|\Xi \right\| \leq C \left( h_x^4+h_y^4\right) \|w\|_{H^6(\Omega)}.
\end{align*}
\end{lem}
%%%%%%%%%%%%%%%%%%%%%%%
\begin{lem}[\cite{XYK}]\label{lem:normEquivalent}
For any grid function $v \in \mathcal{V}_h$,  
\begin{enumerate}[(i)]
	\item if  $v_{0,j} = v_{M_x,j} = 0$, it holds that
	\begin{align*}
	\frac{2}{3}\|v\| \leq  \|\mathcal{H}_{x} v \|\leq \|v\|, 
	\quad \frac{2}{3}\|v\|^2 \le \left(\mathcal{H}_{x} v, v\right) \leq\|v\|^2,
	\quad -\left(\mathcal{H}_{y} \delta_x^2 v, v\right) \geq \frac{2}{3} \|\delta_x v \|_x^2;
	\end{align*}
	\item  if  $v_{i,0} = v_{i,M_y} = 0$, it holds that
	\begin{align*}
	\frac{2}{3}\|v\| \leq  \|\mathcal{H}_y v \|\leq \|v\|, 
	~~\frac{2}{3}\|v\|^2 \le \left( \mathcal{H}_{y} v, v\right)  \leq\|v\|^2,
	~~-\left(\mathcal{H}_{x} \delta_y^2 v, v\right) \geq \frac{2}{3} \|\delta_y v \|_y^2.
	\end{align*}
\end{enumerate}
Furthermore, for any grid function $v \in \mathcal{V}_h^0$, it holds that
$
\frac{4}{9}\|v\| \leq \left\|\mathcal{H}_{h} v\right\|\leq \|v\|.
$
\end{lem}
%%%%%%%%%%%%%%%%%%%%%%%
\begin{lem}[\cite{LSS}]\label{lem:infHL}
 For any grid function $v \in \mathcal{V}_h^0$, it holds that  
$
   \|v\|_{\infty}\leq C\left( \|\mathcal{H}_h v \|+ \|\Lambda_h v \|\right). 
$
\end{lem}
%%%%%%%%%%%%%%%%%%%%%%%
\begin{lem}\label{lem:HLam}
	For any grid function $v \in \mathcal{V}_h^0$, it holds that  
	\begin{equation*}
		({\rm i})\  \left(  \mathcal{H}_h v, \Lambda_h v\right)  \leq 0;\qquad
		({\rm ii})\ \left(  \delta_x^2\delta_y^2 v, \mathcal{H}_h v\right)  \geq 0;\qquad
		({\rm iii})\  \left(  \delta_x^2\delta_y^2 v, \Lambda_h v\right)  \leq 0.
	\end{equation*}
\end{lem}
{\bf Proof.} First, by the definitions of $\mathcal{H}_h$ and $\Lambda_h$, the discrete inner product in (i) can be written as follows:
\begin{equation*}
	\begin{aligned}
		\left(  \mathcal{H}_h v, \Lambda_h v\right)  
		&=\left(  \mathcal{H}_{x} \mathcal{H}_{y} v, \mathcal{H}_{y} \delta_x^2 v\right) +\left(  \mathcal{H}_{x} \mathcal{H}_{y} v, \mathcal{H}_{x} \delta_y^2 v\right)  \\
		&=-\|\mathcal{H}_{y} \delta_x v\|_x^2
		+\frac{h_x^2}{12} \|\mathcal{H}_{y} \delta_x^2 v \|^2
		- \|\mathcal{H}_{x} \delta_y v \|_y^2
		+\frac{h_y^2}{12} \|\mathcal{H}_{x} \delta_y^2 v \|^2,
	\end{aligned}
\end{equation*}
where
\begin{align*}
	\|\mathcal{H}_{y} \delta_x^2 v\|^2 =\| \delta_x^2 \mathcal{H}_{y}v\|^2 
	&= h_x h_y \sum_{i=1}^{M_x-1}\sum_{j=1}^{M_y-1} \left| \delta_x^2 \mathcal{H}_{y} v_{i,j}\right|^2\\
	&= h_x h_y \sum_{i=1}^{M_x-1}\sum_{j=1}^{M_y-1} \left|\frac{ \delta_x \mathcal{H}_{y} v_{i+1/2,j} -  \delta_x \mathcal{H}_{y} v_{i-1/2,j}}{h_x}\right|^2\\
	&\leq \frac{4}{h_x^2} h_x h_y \sum_{i=1}^{M_x}\sum_{j=1}^{M_y-1}\left|\delta_x\mathcal{H}_{y} v_{i-1/2,j}\right|^2 = \frac{4}{h_x^2}  \|\delta_x\mathcal{H}_{y} v \|_x^2,
\end{align*}
and similarly
\begin{equation*}
	\|\mathcal{H}_{x} \delta_y^2 v\|^2 
	\leq   \frac{4}{h_y^2}  \|\delta_y\mathcal{H}_{x} v \|_y^2.
\end{equation*}
Thus, combining all the above results together, the proof of (i) is completed.

Second,  the discrete inner product in (ii) can be written as
\begin{equation*}
	\begin{aligned}
		\left( \delta_x^2\delta_y^2 v, \mathcal{H}_h v\right)
		&=\left( \delta_x^2\delta_y^2 v, \mathcal{H}_x \mathcal{H}_y v\right)\\
		&=\| \delta_x\delta_y v\|_{xy}^2 
		+\frac{h_x^2}{12}\left( \delta_x^2\delta_y^2 v, \delta_x^2 v\right)
		+\frac{h_y^2}{12}\left( \delta_x^2\delta_y^2 v, \delta_y^2 v\right)
		+\frac{h_x^2h_y^2}{144}\| \delta_x^2\delta_y^2 v\|^2.
	\end{aligned}
\end{equation*}
Using summation by parts and discrete inverse estimate, we have
\begin{equation}\label{inv:est}
	\begin{aligned}
		&\left( \delta_x^2\delta_y^2 v, \delta_x^2 v\right)
		=-\| \delta_x^2\delta_y v\|_{y}^2 
		\geq-\frac{4}{h_x^2}\| \delta_x\delta_y v\|_{xy}^2,\\
		&\left( \delta_x^2\delta_y^2 v, \delta_y^2 v\right)
		=-\| \delta_x \delta_y^2 v\|_{x}^2 
		\geq-\frac{4}{h_y^2}\| \delta_x\delta_y v\|_{xy}^2.
	\end{aligned}
\end{equation}
Then, the result of (ii) can be obtained by combining the above results together.

Third,  by the definition and a similar inverse estimate result in \eqref{inv:est} we have 
\begin{equation*}
	\begin{aligned}
		\left( \delta_x^2\delta_y^2 v, \Lambda_h v\right)
		&=\left( \delta_x^2\delta_y^2 v, \mathcal{H}_x\delta_y^2 v\right)
		+\left( \delta_x^2\delta_y^2 v, \mathcal{H}_y\delta_x^2 v\right) \\
		&=-\| \delta_x \delta_y^2 v\|_{x}^2 
		+\frac{h_x^2}{12}\| \delta_x^2\delta_y^2 v\|^2
		-\| \delta_x^2\delta_y v\|_{y}^2
		+\frac{h_y^2}{12}\| \delta_x^2\delta_y^2 v\|^2\\
		&\leq -\frac{2}{3}\| \delta_x \delta_y^2 v\|_{x}^2 -\frac{2}{3}\| \delta_x^2\delta_y v\|_{y}^2
		\leq 0,
	\end{aligned}
\end{equation*}
which proves the result in (iii). \qed
%%%%%%%%%%%%%%%%%%%%%%%%%%%%%%%%%%%%%%%%%%%%%%%%%%%%%%%%%%%%%%

\section{A fourth-order compact difference scheme and its unconditional analysis}\label{sec:compact}
\subsection{Construction of compact difference scheme}\label{subsec:compact}
To construct an unconditional maximum-norm optimal-order convergent HOC difference scheme, we shall introduce two temporal-derivative-dependent variables $p(\bm x, t) = u_t(\bm x, t) - a_{12}v(\bm x,t)$ and $q(\bm x, t) = v_t(\bm x, t)$. 
With these variables, the bacterial system \eqref{bact:eq} can be reformulated into the following equivalent four-variable coupled system:
\begin{subequations}\label{bact:eq:re}
	\begin{numcases}{}
		p = u_t - a_{12} v,                 &\qquad $(\bm x, t) \in \Omega \times J$,\label{bact:eq:re:a}\\
	    -d_{1} \Delta u + a_{11} u = -p,    &\qquad $(\bm x, t) \in \Omega \times J$, \label{bact:eq:re:b}\\
	    q = v_t ,                            &\qquad $(\bm x, t) \in \Omega \times J$,\label{bact:eq:re:c}\\
		-d_{2} \Delta v + a_{22} v - g(u) = -q,     & \qquad $(\bm x, t) \in \Omega \times J$, \label{bact:eq:re:d}
\end{numcases}
\end{subequations}
enclosed with homogeneous Dirichlet boundary conditions and initial values
\begin{equation}\label{bact:eq:bic}
	\begin{cases}{}
	u(\bm x, t)= v(\bm x, t)=p(\bm x, t)= q(\bm x, t)=0,   & (\bm x, t) \in \partial \Omega \times J,\\ 
	u(\bm x, 0)=u_{0}(\bm x), \quad v(\bm x, 0)=v_{0}(\bm x), &  \bm x \in \Omega,\\
	p(\bm x, 0)=d_{1} \Delta u_{0}(\bm x) - a_{11}u_{0}(\bm x)=:p_0(\bm x),  & \bm x \in \Omega, \\
	q(\bm x, 0)=d_{2} \Delta v_{0}(\bm x) - a_{22} v_{0}(\bm x) + g(u_{0}(\bm x))=:q_0(\bm x),	  &\bm x \in \Omega.
\end{cases}
\end{equation}
Below, we shall construct a linearized and decoupled fourth-order compact finite difference scheme for the bacterial system \eqref{bact:eq}   through the equivalent system \eqref{bact:eq:re}--\eqref{bact:eq:bic}. For simplicity, define grid functions
\begin{align*}
u_{i,j}^n=u (x_i, y_j, t_n ), ~~p_{i,j}^n=p (x_i, y_j, t_n ), ~~
v_{i,j}^n=v(x_i, y_j, t_n), ~~q_{i,j}^n=q(x_i, y_j, t_n),
\end{align*}
for $(i, j) \in \bar{\omega}$ and  $0 \leq n \leq N$.

We first consider the discretization of \eqref{bact:eq:re:a}, which only involves temporal derivative. For $2 \leq n \leq N$, it is approximated at $t=t_{n-1/2}$ by  the Crank-Nicholson  and second-order linear extrapolation formulas as 
\begin{align}\label{ext:ut}
\bar{p}_{i,j}^{n-1/2}= \delta_t u_{i,j}^{n-1/2} - a_{12} v_{i,j}^{*,n}+ R_{i,j}^{t,n-1/2}, 
\quad(i, j) \in \omega,
\end{align}
where $v^{*,n} = 3v^{n-1}/2 - v^{n-2}/2$. The truncation error $R_{i,j}^{t,n-1/2} = \mo(\tau^2)$, which can be verified from Lemmas \ref{lem:trunerr:t2}--\ref{lem:trunerr:t3}.
For $n=1$,  we can discretize \eqref{bact:eq:re:a} by  the Crank-Nicholson  and first-order linear extrapolation
formulas  at $t=t_{1/2}$ as
\begin{align}\label{ext:t1}
	\bar{p}_{i,j}^{1/2}= \delta_t u_{i,j}^{1/2} - a_{12} v_{i,j}^{*,1}+R_{i,j}^{t,1/2}, 
	\quad(i, j) \in \omega,
\end{align}
where $v^{*,1} = v^0$, and   Lemmas \ref{lem:trunerr:t1}--\ref{lem:trunerr:t2} show that the truncation error $R_{i,j}^{t,1/2} = \mo(\tau)$.

Next, to approximate \eqref{bact:eq:re:b} with fourth-order consistency, 
we apply the compact operator $\mathcal{H}_h$ to it
%it is easy to know that, for $1 \leq n \leq N$,
%\begin{equation*}
%	H_h p_{i,j}^n = d_1 H_h\Delta u_{i,j}^n - a_{11} H_h u_{i,j}^n,
%	\quad(i, j) \in \omega,
%\end{equation*}
and use Lemma \ref{lem:trunerr:t4} to get
\begin{equation}\label{ext:pxy}
	- d_1 \Lambda_h u_{i,j}^{n} + a_{11} \mathcal{H}_h u_{i,j}^{n} = -\mathcal{H}_h p_{i,j}^n + R_{i,j}^{s,n},
	\quad (i, j) \in \omega,
\end{equation}
for $1 \leq n \leq N$, where the  truncation error $R^{s,n}=\{R_{i,j}^{s,n}\}$ is estimated by Lemma \ref{lem:trunerr:t4} that
\begin{equation}\label{trun:Rs}
	\left\|R^{s,n}\right\| \leq C \left( h_x^4+h_y^4\right) \|u(\cdot,t_n)\|_{H^6(\Omega)}.
\end{equation}

Similarly,  \eqref{bact:eq:re:c} and \eqref{bact:eq:re:d} can be discretized as follows:
\begin{equation}\label{ext:vt}
 \bar{q}_{i,j}^{n-1/2}= \delta_t v_{i,j}^{n-1/2} + Z_{i,j}^{t,n-1/2}, \quad (i, j) \in \omega, 
\end{equation}
\begin{equation} \label{ext:qxy}
	 - d_2 \Lambda_h v_{i,j}^{n} + a_{22} \mathcal{H}_h v_{i,j}^{n} - \mathcal{H}_h g(u_{i,j}^n) = -\mathcal{H}_h q_{i,j}^n + Z_{i,j}^{s,n},
	 \quad(i, j) \in \omega,    
\end{equation}
for $1 \leq n \leq N$. Similarly, the truncation error $Z_{i,j}^{t,n-1/2} = \mo(\tau^2)$ and $Z^{s,n}=\{Z_{i,j}^{s,n}\}$ is estimated by 
\begin{equation}\label{trun:Zs}
	\left\|Z^{s,n}\right\| \leq C \left( h_x^4+h_y^4\right) \|v(\cdot,t_n)\|_{H^6(\Omega)}.
\end{equation}

Let $\{  P^n, U^n, Q^n, V^n\} =\{ P_{i,j}^n, U_{i,j}^n, Q_{i,j}^n, V_{i,j}^n \mid (i, j) \in \bar{\omega}\}  \in\mathcal{V}_h^0 $ be the difference approximations to the exact solutions $\{  p^n, u^n, q^n, v^n\} =\{  p_{i,j}^n, u_{i,j}^n, q_{i,j}^n, v_{i,j}^n \mid (i, j) \in \bar{\omega}\}  $. 
Omitting the small truncation errors in \eqref{ext:ut}--\eqref{ext:pxy} and \eqref{ext:vt}--\eqref{ext:qxy}, and replacing the exact grid solutions  with their numerical approximations,
we propose the following linearized Crank-Nicholson type compact finite difference  (CN-CFD) scheme as follows:
\begin{subequations}\label{bactre:cncfd}
	\begin{numcases}{}
	    \bar{P}_{i,j}^{n-1/2} = \delta_t U_{i,j}^{n-1/2} - a_{12} V_{i,j}^{*,n}, 	    &\qquad $ (i, j) \in \omega, $\label{bactre:cncfd:a}\\
		- d_1 \Lambda_h U_{i,j}^{n} + a_{11} \mathcal{H}_h U_{i,j}^{n} = - \mathcal{H}_h P_{i,j}^n,  	&\qquad $(i, j) \in \omega, $\label{bactre:cncfd:b} \\
		\bar{Q}_{i,j}^{n-1/2} = \delta_t V_{i,j}^{n-1/2} , 		&\qquad $(i, j) \in \omega, $   \label{bactre:cncfd:c}\\
		- d_2 \Lambda_h V_{i,j}^{n} + a_{22} \mathcal{H}_h V_{i,j}^{n} - \mathcal{H}_h g(U_{i,j}^{n})  = -\mathcal{H}_h Q_{i,j}^n,  		&\qquad $(i, j) \in \omega,$  \label{bactre:cncfd:d}
\end{numcases}
\end{subequations}
%where
%\begin{equation}\label{not:timF}
%	V^{*,n} =
%	\begin{cases}{}
%		V^0, &n = 1,\\
%		3V^{n-1}/2-V^{n-2}/2, &n\geq2,
%	\end{cases}
%\end{equation}
for $1 \leq n \leq N$,  enclosed with initial conditions
\begin{equation}\label{bactre:cncfd:ic}
	\begin{cases}{}
    	U_{i,j}^0 = u_0(x_i,y_j),\quad V_{i,j}^0 = v_0(x_i,y_j),		&\quad (i, j) \in \omega,\\
		P_{i,j}^0 = p_0(x_i,y_j),\quad Q_{i,j}^0 = q_0(x_i,y_j),		&\quad (i, j) \in \omega.
	\end{cases}
\end{equation}

%%%%%%%%%%%%%%%%%%%%%%%%
\begin{rem}\label{bact:dispq0}\rm In practical computation, the initial values $p^0$ and $q^0$ can also be approximated via \eqref{bactre:cncfd:b}  and \eqref{bactre:cncfd:d}  with fourth-order consistency, i.e.,
\begin{equation}\label{cncfd:iv}
	\begin{aligned}
	\mathcal{H}_h P_{i,j}^0 = d_1 \Lambda_h U_{i,j}^{0} - a_{11} \mathcal{H}_h U_{i,j}^{0},\quad
	\mathcal{H}_h Q_{i,j}^0 = d_2 \Lambda_h V_{i,j}^{0} - a_{22} \mathcal{H}_h V_{i,j}^{0} + \mathcal{H}_h g(U_{i,j}^0),
    \end{aligned}\end{equation}
for $(i, j) \in \omega$, where the truncation errors  $R_{i,j}^{s,0}$ and $Z_{i,j}^{s,0}$ are estimated by \eqref{trun:Rs} and \eqref{trun:Zs}, respectively.
\end{rem}
%%%%%%%%%%%%%%%%%%%%%%%%
\begin{rem}\label{rem:re=cfd} \rm By eliminating the auxiliary variables $P^n$ and $Q^n$,  the  difference system \eqref{bactre:cncfd} is equivalent to the following decoupled  CN-CFD scheme with respect to $\left\lbrace U^n, V^n\right\rbrace_{n \ge 2}  \in\mathcal{V}_h^0$
	\begin{equation}\label{bact:cncfd}
	\begin{cases}{}
		\mathcal{H}_h \delta_t U_{i,j}^{n-1/2}  = d_1 \Lambda_h \bar{U}_{i,j}^{n-1/2} - a_{11} \mathcal{H}_h \bar{U}_{i,j}^{n-1/2} + a_{12} \mathcal{H}_h V_{i,j}^{*,n},
		&(i, j) \in \omega, \\
		\mathcal{H}_h \delta_t V_{i,j}^{n-1/2}  = d_2 \Lambda_h \bar{V}_{i,j}^{n-1/2} - a_{22} \mathcal{H}_h \bar{V}_{i,j}^{n-1/2} +  \mathcal{H}_h \left[ \bar{g}(U)\right]_{i,j}^{n-1/2}, 
		&(i, j) \in \omega,
	\end{cases}
\end{equation}
 as well as the decoupled scheme for $\left\lbrace U^1, V^1\right\rbrace  \in\mathcal{V}_h^0$
\begin{equation}\label{bact:cncfd:t1}
	\left\lbrace 
	\begin{aligned}
		\mathcal{H}_h \delta_t U_{i,j}^{1/2}  &= \frac{d_1}{2} \Lambda_h U_{i,j}^{1} - \frac{a_{11}}{2} \mathcal{H}_h U_{i,j}^{1} + \frac{1}{2}\mathcal{H}_h P_{i,j}^0 + \frac{a_{12}}{2} \mathcal{H}_h V_{i,j}^{*,1}, &(i, j) \in \omega,\\
		\mathcal{H}_h \delta_t V_{i,j}^{1/2}  & = \frac{d_2}{2} \Lambda_h V_{i,j}^{1} - \frac{a_{22}}{2} \mathcal{H}_h V_{i,j}^{1} + \frac{1}{2}\mathcal{H}_h Q_{i,j}^0 + \frac{1}{2} \mathcal{H}_h [g(U)]_{i,j}^{1}, 		 & (i, j) \in \omega,
	\end{aligned}
	\right. 
\end{equation}
where initial conditions are given via \eqref{bactre:cncfd:ic}.
\end{rem}
%%%%%%%%%%%%%%%%%%%%%%%%
\begin{rem} \rm Below  the coupled CN-CFD difference system \eqref{bactre:cncfd} is used solely for the unconditional maximum-norm numerical analysis, in which the temporal-spatial error is split, and the unconditional convergence can be proved. However, in real numerical calculation, we use the decoupled  CN-CFD scheme \eqref{bact:cncfd}--\eqref{bact:cncfd:t1} to successively get the solutions $ U^n$ and $ V^n$, and then $ P^n$ and $ Q^n $  can be directly derived through \eqref{bactre:cncfd:a}  and \eqref{bactre:cncfd:c}, respectively. 
\end{rem}

%%%%%%%%%%%%%%%%%%%%%%%%%%%%%%%%%%%%%%%%%%%%%%%%%%%%%%%%%%%%%%
\subsection{Analysis of the  CN-CFD scheme}
To analyze the stability and convergence of the  CN-CFD scheme \eqref{bactre:cncfd}--\eqref{bactre:cncfd:ic}, we first prove the following   \textit{a priori}  estimate.
\begin{thm}\label{thm:priori}
Assume that grid functions $\{s^n, e^n, l^n, w^n \mid 1 \leq n \leq N\}\in\mathcal{V}_h^0$ are the solutions of the following CN-CFD system with given initial values $\{ s^0, e^0, l^0, w^0\}$ and data $\{ \xi^{n-1/2},\eta^{n},\zeta^{n-1/2},\lambda^{n}\} $:
\begin{subequations}\label{Priori:cncfd}
	\begin{numcases}{}
		\bar{e}_{i,j}^{n-1/2} = \delta_t s_{i,j}^{n-1/2} - a_{12} l_{i,j}^{*,n} + \xi_{i,j}^{n-1/2}, 
		&\qquad $ (i, j) \in \omega, $\label{Priori:cncfd:a}\\[0.05in]
		- d_1 \Lambda_h s_{i,j}^{n} + a_{11} \mathcal{H}_h s_{i,j}^{n} = - \mathcal{H}_h e_{i,j}^n + \eta_{i,j}^{n},  
		&\qquad $(i, j) \in \omega, $\label{Priori:cncfd:b} \\
		\bar{w}_{i,j}^{n-1/2} = \delta_t l_{i,j}^{n-1/2} + \zeta_{i,j}^{n-1/2} , 
		&\qquad $(i, j) \in \omega, $\label{Priori:cncfd:c}\\
		- d_2 \Lambda_h l_{i,j}^{n} + a_{22} \mathcal{H}_h l_{i,j}^{n} - \mathcal{H}_h g(s_{i,j}^{n})  = -\mathcal{H}_h w_{i,j}^n + \lambda_{i,j}^{n} ,  &\qquad $(i, j) \in \omega,$\label{Priori:cncfd:d}
	\end{numcases}
\end{subequations}
for $1 \leq n \leq N$, where $l^{*,n} $ defined as 
\begin{equation}\label{not:timF}
	l^{*,n} =
	\begin{cases}{}
			l^0, &n = 1,\\
			3l^{n-1}/2-l^{n-2}/2, &n\geq2,
		\end{cases}
\end{equation}
Then, there exists a positive constant $C$ which is only related to the coefficients $d_1$, $a_{11}$, $a_{12}$, $d_2$, $a_{22}$ and the Lipschitz constant $L$ such that
\begin{align*}
	\left\|  s^m \right\| _\infty^2 + \left\|  e^m \right\| ^2+ \left\|  l^m \right\| _\infty^2 + \left\|  w^m \right\|^2
	&\leq C\left(\Gamma^0 
	+  \left\|\eta^m\right\|^2 + \left\|\lambda^m\right\|^2
	+ \tau\sum_{n=1}^{m} \Psi^{n-1/2}\right), \quad 1 \leq m \leq N,
\end{align*}
where 
\begin{align}\label{apriori:Gamma}
\Gamma^0  = \|  s^0\|^2 + \|  e^0\|^2+ \|  l^0\|^2 + \| w^0\|^2,
\end{align}
\begin{align}\label{apriori:psi}
\Psi^{n-1/2} 
= \sum_{\chi=\{\xi, \zeta\}} \left( \|  \chi^{n-1/2}\|^2 + \| \Lambda_h \chi^{n-1/2}\|^2  \right) 
+ \sum_{\chi=\{\eta, \lambda\}} \left( \| \bar{\chi}^{n-1/2}\|^2 + \|\delta_t \chi^{n-1/2}\|^2\right) .
\end{align}
\end{thm}
{\bf Proof.}
%To begin with, we first estimate the maximum-norm of $  s^n$. 
From \eqref{Priori:cncfd:b}, it is seen that $\|\Lambda_h s^{n}\|$ can be controlled by $\|\mathcal{H}_h s^{n}\|$ and $\|\mathcal{H}_h e^{n}\|$, which means that once the estimates  $\|\mathcal{H}_h s^{n}\|$ and $\|\mathcal{H}_h e^{n}\|$ are obtained, by Lemma \ref{lem:infHL} we can easily get the estimate for $\left\|  s^n \right\| _\infty$. Also, $\left\|  l^n \right\| _\infty$ can be controlled by $\|\mathcal{H}_h l^{n}\|$ and $\|\mathcal{H}_h w^{n}\|$ through \eqref{Priori:cncfd:d}. In the following, we split the proof of  Theorem \ref{thm:priori}   in three main steps.

\paragraph{\bf Step I. Estimates for $\|\mathcal{H}_h s^{n}\|$ and $\|\mathcal{H}_h e^{n}\|$}

We shall first estimate $\|\mathcal{H}_h s^{n}\|$. By acting the compact operator $\mathcal{H}_h$ to \eqref{Priori:cncfd:a}, then taking the average of \eqref{Priori:cncfd:b} for the superscripts $n-1$ and $n$, and finally adding the two resulting equations, we have
\begin{align}\label{sta:uab:H}
	 \mathcal{H}_h \delta_t s_{i,j}^{n-1/2} 
	=d_1 \Lambda_h \bar{s}_{i,j}^{n-1/2} - a_{11} \mathcal{H}_h \bar{s}_{i,j}^{n-1/2}
	+ a_{12} \mathcal{H}_h l_{i,j}^{*,n} - \mathcal{H}_h \xi_{i,j}^{n-1/2}
	+\bar{\eta}_{i,j}^{n-1/2}, \quad (i, j) \in \omega.
\end{align}

Now, multiplying \eqref{sta:uab:H} by $h_x h_y \mathcal{H}_h \bar{s}_{i j}^{n-1/2}$ and then summing over all $(i, j) \in \omega$, i.e., taking discrete inner products with $ \mathcal{H}_h \bar{s}^{n-1/2}$, we  obtain
\begin{equation}\label{sta:uab:innH}
	\begin{aligned}
	& \left( \mathcal{H}_h \delta_t s^{n-1/2}, \mathcal{H}_h \bar{s}^{n-1/2}\right) \\
	& \quad = d_1 \left( \Lambda_h \bar{s}^{n-1/2}, \mathcal{H}_h \bar{s}^{n-1/2}\right) 
	- a_{11} \left(  \mathcal{H}_h \bar{s}^{n-1/2},\mathcal{H}_h \bar{s}^{n-1/2}\right) 
	+ a_{12}\left(   \mathcal{H}_h l^{*,n},\mathcal{H}_h \bar{s}^{n-1/2}\right)    \\
	&\qquad  
	    - \left(   \mathcal{H}_h \xi^{n-1/2},\mathcal{H}_h \bar{s}^{n-1/2}\right)     
	    + \left(  \bar{\eta}^{n-1/2},\mathcal{H}_h \bar{s}^{n-1/2}\right).  
\end{aligned}
\end{equation}

Noticing that the left-hand side of \eqref{sta:uab:innH} can be expressed as 
\begin{align}\label{sta:uab:innHL}
\left(  \mathcal{H}_h \delta_t s^{n-1/2}, \mathcal{H}_h \bar{s}^{n-1/2}\right) 
= \frac{ \| \mathcal{H}_h s^n \|^2 -  \| \mathcal{H}_h s^{n-1} \|^2}{2\tau},
\end{align}
and the first two right-hand side terms can be respectively estimated by Lemma \ref{lem:HLam} and the definition of discrete inner product that
\begin{align}\label{sta:uab:innHR12}
	d_1 \left(  \Lambda_h \bar{s}^{n-1/2}, \mathcal{H}_h \bar{s}^{n-1/2}\right) \leq 0,\quad
	- a_{11} \left(  \mathcal{H}_h \bar{s}^{n-1/2},\mathcal{H}_h \bar{s}^{n-1/2}\right) \leq 0.
\end{align}
While, the other three terms can be estimated  as follows
\begin{align}
	 \left(  \mathcal{H}_h l^{*,n},\mathcal{H}_h \bar{s}^{n-1/2}\right) 
	&   \le \frac{1}{2}  \| \mathcal{H}_h l^{*,n} \| ^2 + \frac{1}{4}\left(  \| \mathcal{H}_h s^n \|^2 +  \| \mathcal{H}_h s^{n-1} \|^2\right),        \label{sta:uab:innHR3}\\
	-\left(  \mathcal{H}_h \xi^{n-1/2},\mathcal{H}_h \bar{s}^{n-1/2}\right)  
	&\leq \frac{1}{2}  \| \mathcal{H}_h \xi^{n-1/2} \| ^2 + \frac{1}{4}\left(  \| \mathcal{H}_h s^n \|^2 +  \| \mathcal{H}_h s^{n-1} \|^2\right),       \label{sta:uab:innHR4}\\
	 \left( \bar{\eta}^{n-1/2}, \mathcal{H}_h \bar{s}^{n-1/2}\right)  
	&\leq  \frac{1}{2}\| \bar{\eta}^{n-1/2}\|^2
	+\frac{1}{4}\left(  \| \mathcal{H}_h s^n \|^2 +  \| \mathcal{H}_h s^{n-1} \|^2 \right).    \label{sta:uab:innHR5}
\end{align}
Therefore, inserting the estimates \eqref{sta:uab:innHL}--\eqref{sta:uab:innHR5} together into \eqref{sta:uab:innH},  and multiplying the resulting equation by $2\tau$, we have
\begin{equation}
	\begin{aligned}\label{sta:u:H}
	 & \| \mathcal{H}_h s^n \|^2 - \| \mathcal{H}_h s^{n-1} \|^2 \\
		&\quad \leq C\tau   \| \mathcal{H}_h l^{*,n} \| ^2  
		+ C\tau\left( \| \mathcal{H}_h \xi^{n-1/2}\|^2 
		+ \| \bar{\eta}^{n-1/2}\|^2\right) 
		+ C\tau \left(  \| \mathcal{H}_h s^n \|^2 +  \| \mathcal{H}_h s^{n-1} \|^2\right).
\end{aligned}
\end{equation}

Next, we shall estimate $\|\mathcal{H}_h e^{n}\|$. Making a difference for the superscripts $n$ and $n-1$ to \eqref{Priori:cncfd:b} and divided the resulting equation by $\tau$, further taking discrete inner product with $ \mathcal{H}_h \bar{e}^{n-1/2}$,  we have
 \begin{equation}\label{sta:ub:innL}
\begin{aligned}
&	%\frac{  \| \mathcal{H}_h e^{n} \| ^2 -   \| \mathcal{H}_h e^{n-1} \| ^2}{2\tau}=
\left( \mathcal{H}_h\delta_t e^{n-1/2}, \mathcal{H}_h \bar{e}^{n-1/2}\right)  \\
	&\quad = d_1 \left(  \Lambda_h\delta_t s^{n-1/2}, \mathcal{H}_h \bar{e}^{n-1/2}\right)   
	- a_{11} \left(  \mathcal{H}_h\delta_t s^{n-1/2},\mathcal{H}_h \bar{e}^{n-1/2}\right)   
	+ \left(  \delta_t \eta^{n-1/2}, \mathcal{H}_h \bar{e}^{n-1/2}\right).
\end{aligned}
\end{equation}

To estimate the first  term on the right-hand side of \eqref{sta:ub:innL}, we act the compact difference operator $\Lambda_h$ to \eqref{Priori:cncfd:a},  and then take discrete inner product with  $\mathcal{H}_h \bar{e}^{n-1/2}$ to get
 \begin{equation}\label{sta:ua:innL:e1} 
	\begin{aligned}
		\left( \Lambda_h \delta_t s^{n-1/2}, \mathcal{H}_h \bar{e}^{n-1/2}\right)   
	=\left( \Lambda_h \bar{e}^{n-1/2}, \mathcal{H}_h \bar{e}^{n-1/2} \right) 
	+ a_{12}\left(  \Lambda_h l^{*,n},\mathcal{H}_h \bar{e}^{n-1/2}\right)   
	- \left( \Lambda_h \xi^{n-1/2} , \mathcal{H}_h \bar{e}^{n-1/2}\right) . 
\end{aligned}
\end{equation}
Thus, using \eqref{sta:ua:innL:e1} and Lemma \ref{lem:HLam} we see
\begin{equation}\label{sta:ua:innL:est}
\begin{aligned}
 d_1 \left(  \Lambda_h \delta_t s^{n-1/2}, \mathcal{H}_h \bar{e}^{n-1/2}\right)   
& \leq C\left( \left\|\Lambda_h l^{*,n}\right\|^2 
+ \| \Lambda_h \xi^{n-1/2}\|^2
+   \| \mathcal{H}_h e^{n} \| ^2 
+   \| \mathcal{H}_h e^{n-1} \| ^2\right) .  
\end{aligned}
\end{equation}

Moreover, by applying Cauchy-Schwarz inequality and using \eqref{Priori:cncfd:a} again, one can estimate the second term on the right-hand side of \eqref{sta:ub:innL} by
\begin{equation}\label{sta:ub:innL:rhs2}
	\begin{aligned}
&- a_{11} \left(  \mathcal{H}_h\delta_t s^{n-1/2},\mathcal{H}_h \bar{e}^{n-1/2}\right)\\
&\quad =- a_{11} \left(  \mathcal{H}_h\bar{e}^{n-1/2},\mathcal{H}_h \bar{e}^{n-1/2}\right) 
  - a_{11}a_{12} \left(  \mathcal{H}_h l^{*,n},\mathcal{H}_h \bar{e}^{n-1/2}\right)      
  + a_{11} \left(  \mathcal{H}_h \xi^{n-1/2},\mathcal{H}_h \bar{e}^{n-1/2}\right)    \\
&\quad \leq C\left( \|\mathcal{H}_h l^{*,n}\|^2   
      + \|\mathcal{H}_h \xi^{n-1/2}\|^2
      + \| \mathcal{H}_h e^{n} \| ^2 +  \| \mathcal{H}_h e^{n-1} \| ^2\right) . 
\end{aligned}
\end{equation}

Finally, the last term on the right-hand side of \eqref{sta:ub:innL} can be estimated by  Cauchy-Schwarz inequality that
\begin{equation}\label{sta:ub:innL:rhs3}
	\begin{aligned}
	\left(  \delta_t \eta^{n-1/2}, \mathcal{H}_h \bar{e}^{n-1/2}\right) 
	\leq \frac{1}{2}\|\delta_t \eta^{n-1/2}\|^2
	+ \frac{1}{4}\left(  \| \mathcal{H}_h e^{n} \| ^2 +  \| \mathcal{H}_h e^{n-1} \| ^2  \right).
\end{aligned}
\end{equation}

Now, inserting the estimates \eqref{sta:ua:innL:est}--\eqref{sta:ub:innL:rhs3} together into \eqref{sta:ub:innL}, and using a similar lower bound \eqref{sta:uab:innHL} for the left-hand side of \eqref{sta:ub:innL}, we can derive
 \begin{equation}\label{sta:u:L}
	\begin{aligned}
	& \| \mathcal{H}_h e^n\|^2 -\| \mathcal{H}_h e^{n-1}\|^2  \\ 
		&\quad \leq C\tau\left( \|\mathcal{H}_h l^{*,n} \|^2 + \|\Lambda_h l^{*,n} \|^2 \right) 
		  + C\tau\left( \|\mathcal{H}_h \xi^{n-1/2} \|^2 
		  + \| \Lambda_h \xi^{n-1/2} \|^2     
		  + \|\delta_t \eta^{n-1/2} \|^2\right) \\
		  &\qquad 
		  + C\tau\left( \| \mathcal{H}_h e^n \|^2 + \| \mathcal{H}_h e^{n-1} \|^2  \right).
\end{aligned}
\end{equation}
%where $C_1 = (d_1a_{12} + a_{11}a_{12} + a_{11} + d_1 +1)/2$. 
%Also, we observe that  the estimate $\|\mathcal{H}_h e^{n}\|$ depends directly on the estimate $\left\|\Lambda_h l^{*,n}\right\|$ which be given in the last step.

Noting the definition of $l^{*,n}$ in \eqref{not:timF}  is different for $n=1$ and $n\geq2$.  Therefore, by adding \eqref{sta:u:H} and \eqref{sta:u:L} together for $n\ge 2$, and summing over $n$ from $2$ to $m $,   we  get
\begin{equation}\label{sta:u:sum}
\begin{aligned}
		 \| \mathcal{H}_h s^m\|^2 +\| \mathcal{H}_h e^m \|^2  
		& \le \| \mathcal{H}_h s^1 \|^2 + \| \mathcal{H}_h e^1 \|^2 
		+ C \tau \sum_{n=0}^{m-1}  \left( \|\mathcal{H}_h l^{n} \|^2  	 +  \|\Lambda_h l^{n}\|^2 \right) \\
		&\qquad	+ C\tau\sum_{n=2}^{m} \left(  \| \mathcal{H}_h \xi^{n-1/2}\|^2  
	     + \| \bar{\eta}^{n-1/2}\|^2+  \| \Lambda_h \xi^{n-1/2} \|^2 	+  \|\delta_t \eta^{n-1/2} \|^2 \right) \\
	   &\qquad + C \tau \sum_{n=1}^{m} \left(   \| \mathcal{H}_h s^n \|^2 	+\| \mathcal{H}_h e^n\|^2 \right).
\end{aligned}
\end{equation}
%where $C_2 = \max\left\{ 5(a_{12} + a_{11}a_{12}), 5d_1a_{12}\right\}$ and $C_3 = \max\left\{a_{12} + 2, 2C_1\right\}$. 

For $n=1$, the estimates $\| \mathcal{H}_h s^1\|$ and $\| \mathcal{H}_h e^1\|$ at the first time level can be obtained through   \eqref{not:timF}, \eqref{sta:u:H} and \eqref{sta:u:L} as
\begin{equation}\label{sta:u:sum:t1}
	\begin{aligned}
		  \| \mathcal{H}_h s^1 \|^2 + \| \mathcal{H}_h e^1 \|^2  
		& \leq  \| \mathcal{H}_h s^0 \|^2 + \| \mathcal{H}_h e^0 \|^2 
		+  C \tau \left(  \| \mathcal{H}_h l^{0} \|^2 
		+  \|\Lambda_h l^0 \|^2 \right)    \\
		&\qquad+ C\tau\left(  \| \mathcal{H}_h \xi^{1/2} \|^2 
		+ \| \bar{\eta}^{1/2}\|^2 
		+ \| \Lambda_h \xi^{1/2}\|^2 
		+ \|\delta_t \eta^{1/2}\|^2 \right)\\
		&\qquad+ C\tau\left(  \| \mathcal{H}_h s^1 \|^2 
		+  \| \mathcal{H}_h s^{0} \|^2 +  \| \mathcal{H}_h e^{1} \|^2 
		+ \| \mathcal{H}_h e^{0} \|^2\right). 
\end{aligned}
\end{equation}

Therefore, inserting \eqref{sta:u:sum:t1} into \eqref{sta:u:sum}, we have the following estimate 
\begin{equation}\label{sta:u:sum:tol}
\begin{aligned}
		  \| \mathcal{H}_h s^m \|^2 +  \| \mathcal{H}_h e^{m} \| ^2
		&  \leq  \| \mathcal{H}_h s^{0} \|^2 + \| \mathcal{H}_h e^{0} \|^2 
		+ C\tau\sum_{n=0}^{m-1}  \left( \left\|\mathcal{H}_h l^{n}\right\|^2  + \left\|\Lambda_h l^{n}\right\|^2\right)   \\
		&\qquad + C\tau\sum_{n=1}^{m}\left[  \| \mathcal{H}_h \xi^{n-1/2}\|^2  + \| \bar{\eta}^{n-1/2}\|^2
		+ \| \Lambda_h \xi^{n-1/2}\|^2 
		+ \|\delta_t \eta^{n-1/2}\|^2 \right]  \\
		&\qquad+ C\tau\sum_{n=0}^{m}\left(  \| \mathcal{H}_h s^n \|^2 
		+   \| \mathcal{H}_h e^{n} \| ^2 \right) . 
\end{aligned}
\end{equation}

It can be seen from \eqref{sta:u:sum:tol} that  the estimates $\|\mathcal{H}_h s^{m}\|$ and $  \| \mathcal{H}_h e^{m} \| $ depend directly on the estimates $\|\mathcal{H}_h l^{n}\|$ and $\|\Lambda_h l^{n}\|$   $(0\le n \le m-1)$ which will be given in the next step, and  once $\|\mathcal{H}_h l^{n}\|$ and $\|\Lambda_h l^{n}\|$ are analyzed, the estimates  $\| \mathcal{H}_h s^m\|$ and $\| \mathcal{H}_h e^m\|$ can be obtained immediately by virtue of discrete Gr\"{o}nwall's inequality. 
Furthermore, recalling that from \eqref{Priori:cncfd:d} one see $\|\Lambda_h l^{n}\|$ can be controlled by $\|\mathcal{H}_h l^{n}\|$ and $\|\mathcal{H}_h w^{n}\|$. Thus, to complete the proof, we finally have to  estimate  $\|\mathcal{H}_h l^{n}\|$ and $\|\mathcal{H}_h w^{n}\|$ that will be given as below.

\paragraph{\bf Step II. Estimates for $\|\mathcal{H}_h l^{n}\|$ and $\|\mathcal{H}_h w^{n}\|$}
 
We first estimate $\|\mathcal{H}_h l^{n}\|$  through equations \eqref{Priori:cncfd:c}--\eqref{Priori:cncfd:d}. 
Noticing that $g$ is  Lipschitz continuous, by Lemma \ref{lem:normEquivalent} one get
\begin{equation}\label{lips:g}
	\begin{aligned}
	\|\mathcal{H}_h g(s^{n}) - \mathcal{H}_h g(s^{n-1})\|
	\leq L \| s^n - s^{n-1}\|
	\leq \frac{9L}{4} \| \mathcal{H}_h (s^n - s^{n-1})\|,
\end{aligned}
\end{equation}
which is involved in the analysis below. 

Following the same line of proof for $\|\mathcal{H}_h s^{n}\|$, utilizing Lemma \ref{lem:HLam} and Cauchy-Schwarz inequality, as well as the assumption that $g(0)=0$, we obtain the following similar result
\begin{equation}\label{sta:v:H}
	\begin{aligned}
	& \| \mathcal{H}_h l^n\|^2 - \| \mathcal{H}_h l^{n-1} \|^2  \\
	&\quad \leq \tau\left(  \| \mathcal{H}_h l^n \|^2 +  \| \mathcal{H}_h l^{n-1} \|^2
	+  \| \mathcal{H}_h \zeta^{n-1/2}\|^2 
	+  \| \bar{\lambda}^{n-1/2}\|^2\right)    
    + 2\tau\left(  \mathcal{H}_h [\bar{g}(s)]^{n-1/2},\mathcal{H}_h \bar{l}^{n-1/2}\right)    \\
	& \quad \leq C\tau\left(  \| \mathcal{H}_h s^n \|^2 +  \| \mathcal{H}_h s^{n-1} \|^2\right)
	+ C\tau\left(\| \mathcal{H}_h \zeta^{n-1/2}\|^2 
	+  \| \bar{\lambda}^{n-1/2}\|^2\right)   \\
   	&\qquad+ C\tau\left(  \| \mathcal{H}_h l^n \|^2 +  \| \mathcal{H}_h l^{n-1} \|^2\right),
\end{aligned}
\end{equation}
where we have used the estimate
\begin{equation}\label{sta:v:H:g}
	\begin{aligned}
	\left(  \mathcal{H}_h [\bar{g}(s)]^{n-1/2},\mathcal{H}_h \bar{l}^{n-1/2}\right)  
	&\leq \|\mathcal{H}_h [\bar{g}(s)]^{n-1/2}\| \, \|\mathcal{H}_h \bar{l}^{n-1/2}\| \\
	&\leq  C\left(  \| \mathcal{H}_h s^n \|^2 +  \| \mathcal{H}_h s^{n-1} \|^2 +
	    \| \mathcal{H}_h l^n \|^2 +  \| \mathcal{H}_h l^{n-1} \|^2 \right).
\end{aligned}
\end{equation}

Next, we estimate $\|\mathcal{H}_h w^{n}\|$ also through equations \eqref{Priori:cncfd:c}--\eqref{Priori:cncfd:d}. 
Following the same estimate procedure for $\|\mathcal{H}_h e^{n}\|$, we obtain the following similar result
\begin{equation}\label{sta:v:Lmid}
	\begin{aligned}
	& \| \mathcal{H}_h w^{n} \|^2 - \| \mathcal{H}_h w^{n-1} \|^2   \\
	& \quad \leq   C\tau\left( \| \Lambda_h \zeta^{n-1/2}\|^2    
	+  \| \delta_t \lambda^{n-1/2}\|^2 \right) 
	+ C\tau\left(  \| \mathcal{H}_h w^{n} \|^2 + \| \mathcal{H}_h w^{n-1} \|^2 \right) \\
	&\qquad
	+ 2a_{22}\tau\, \left(  \mathcal{H}_h \delta_t l^{n-1/2},\mathcal{H}_h \bar{w}^{n-1/2}\right)    
	+ 2\tau\, \left(  \mathcal{H}_h \delta_t [g(s)]^{n-1/2},\mathcal{H}_h \bar{w}^{n-1/2}\right),
\end{aligned}
\end{equation}
where  Lemma \ref{lem:HLam} and Cauchy-Schwarz inequality have been applied.

Considering equation \eqref{Priori:cncfd:c},  the first discrete inner product term on the right-hand side of \eqref{sta:v:Lmid} can be estimated as
\begin{equation}\label{sta:v:Lmid:lt}
	\begin{aligned}
	2a_{22}\tau\left( \mathcal{H}_h \delta_t l^{n-1/2},\mathcal{H}_h \bar{w}^{n-1/2}\right) & =	2a_{22}\tau\left(\mathcal{H}_h \bar{w}^{n-1/2},\mathcal{H}_h \bar{w}^{n-1/2}\right)-	2a_{22}\tau\left( \mathcal{H}_h \zeta^{n-1/2},\mathcal{H}_h \bar{w}^{n-1/2}\right)\\
	&\leq C \tau\,\| \mathcal{H}_h \zeta^{n-1/2}\|^2
	+  C\tau\left(  \| \mathcal{H}_h w^{n} \|^2 + \| \mathcal{H}_h w^{n-1} \|^2 \right),
\end{aligned}
\end{equation}
and the last term of \eqref{sta:v:Lmid} can be estimated  using Cauchy-Schwarz inequality that
\begin{equation}\label{sta:v:Lmid:gt}
	\begin{aligned}
		2\tau( \mathcal{H}_h \delta_t [g(s)]^{n-1/2},\mathcal{H}_h \bar{w}^{n-1/2})
		&\leq 2\tau\|\mathcal{H}_h \delta_t[g(s)]^{n-1/2}\| \, \|\mathcal{H}_h \bar{w}^{n-1/2}\| \\
		&\leq C\tau\left(   \| \mathcal{H}_h e^{n} \| ^2 
		+  \| \mathcal{H}_h e^{n-1} \| ^2\right)  
		+ C\tau\,\| \mathcal{H}_h \xi^{n-1/2} \| ^2 \\
		&\qquad     + C\tau\left( \| \mathcal{H}_h l^{*,n} \| ^2+
		\| \mathcal{H}_h w^{n} \|^2 + \| \mathcal{H}_h w^{n-1} \|^2 \right) ,
	\end{aligned}
\end{equation}
in which \eqref{lips:g}, Lemma \ref{lem:normEquivalent} and  \eqref{Priori:cncfd:a} are utilized to show
\begin{equation}\label{sta:v:Lmid:gtl2}
	\begin{aligned}
		\|\mathcal{H}_h \delta_t[g(s)]^{n-1/2}\| 
		&\leq L\| \delta_t s^{n-1/2}\| 
		\leq L\, \left( \| \bar{e}^{n-1/2}\| + a_{12}\, \| l^{*n}\|+ \| \xi^{n-1/2}\|\right) \\
		&\leq C\left( \| \mathcal{H}_h\bar{e}^{n-1/2}\| + \|\mathcal{H}_h \xi^{n-1/2}\|+ \|\mathcal{H}_h l^{*n}\|\right).
	\end{aligned}
\end{equation}

Now, inserting the estimates \eqref{sta:v:Lmid:lt}--\eqref{sta:v:Lmid:gt} into \eqref{sta:v:Lmid}, we obtain  
\begin{equation}\label{sta:v:L}
\begin{aligned}
 &\| \mathcal{H}_h w^{n} \|^2 - \| \mathcal{H}_h w^{n-1} \|^2   \\
&\quad\leq C\tau\left(   \| \mathcal{H}_h e^{n} \| ^2 +  \| \mathcal{H}_h e^{n-1} \| ^2 \right)  
 +  C\tau \left( \| \Lambda_h \zeta^{n-1/2}\|^2 
+ \|\mathcal{H}_h \zeta^{n-1/2}\|^2   
+    \| \mathcal{H}_h \xi^{n-1/2} \| ^2 
+  \| \delta_t \lambda^{n-1/2}\|^2\right) \\
&\qquad+ C\tau\left(  \| \mathcal{H}_h l^{*,n} \| ^2 + \| \mathcal{H}_h w^{n} \|^2 + \| \mathcal{H}_h w^{n-1} \|^2 \right) . 
\end{aligned}
\end{equation}
%where $C_4= (d_2 + 3a_{22} + 1)/2 + 9L/8$.

Similar as the proof in Step I, we add \eqref{sta:v:H} and \eqref{sta:v:L} together,  and following the same discussions as in \eqref{sta:u:sum}--\eqref{sta:u:sum:t1} to obtain
\begin{equation}\label{sta:v:sum}
	\begin{aligned}
		&  \| \mathcal{H}_h l^m \|^2 + \| \mathcal{H}_h w^m \|^2     \\
		&\quad \leq  \| \mathcal{H}_h l^0 \|^2 + \| \mathcal{H}_h w^0\|^2  
		+ C  \tau\sum_{n=0}^{m} \left(   \| \mathcal{H}_h s^n \|^2 +    \| \mathcal{H}_h e^{n} \| ^2\right)\\
		&\qquad 
		+ C  \tau\sum_{n=1}^{m}\left( \| \mathcal{H}_h \xi^{n-1/2}\|^2   
		+ \| \mathcal{H}_h \zeta^{n-1/2}\|^2
		+ \|  \bar{\lambda}^{n-1/2}\|^2  +  \| \Lambda_h \zeta^{n-1/2}\|^2     
		+ \|\delta_t \lambda^{n-1/2}\|^2 \right)\\
		& 	\qquad	+ C  \tau\sum_{n=0}^{m} \left(  \| \mathcal{H}_h l^n \|^2 +  \| \mathcal{H}_h w^{n} \|^2 \right).
\end{aligned}
\end{equation}
%where $C_5 = \max\left\lbrace 27L/8, 1+a_{22}, d_2\right\rbrace $ and $C_6 = \max\lbrace (18L + 16 + 27a_{12}^2L)/8, 2C_4, 27L/8 \rbrace $.

In fact,  we can also observe from \eqref{sta:v:sum} that the estimates  $ \| \mathcal{H}_h l^m \|$ and $ \| \mathcal{H}_h w^m \|$ are in turn related to  $ \| \mathcal{H}_h s^n \|$ and $  \| \mathcal{H}_h e^{n} \| $ for $0\le n \le m$. Therefore, by adding up \eqref{sta:u:sum:tol} and \eqref{sta:v:sum}, we can obtain
\begin{equation}\label{sta:uv:sum}
	\begin{aligned}
		&  \| \mathcal{H}_h s^m \|^2 +  \| \mathcal{H}_h e^{m} \| ^2
		+  \| \mathcal{H}_h l^m \|^2 + \| \mathcal{H}_h w^m \|^2     \\
		&\quad\leq  \Gamma^0  + C\tau\sum_{n=1}^{m} \Psi^{n-1/2}
		+ C\tau\sum_{n=0}^{m} \left(  
		    \| \mathcal{H}_h s^n \|^2 +    \| \mathcal{H}_h e^{n}\|^2
		+ \| \mathcal{H}_h l^n \|^2 +  \| \mathcal{H}_h w^{n} \|^2   + \left\|\Lambda_h l^{n}\right\|^2 \right),
	\end{aligned}
\end{equation}
where $\Gamma^0$ and $\Psi^{n-1/2}$ are  respectively defined by \eqref{apriori:Gamma} and \eqref{apriori:psi}.
%where $C_7 = 2(C_5+a_{11}+1)\exp\left( 2T(2C_2 + 2 C_3 + C_6)\right) $.

\paragraph{\bf Step III. Estimates for  $\|s^{m}\|_\infty$ and $\|l^{m}\|_\infty$}

Obviously, the maximum-norm of $s^m$ and $l^m$ can be followed from \eqref{sta:uv:sum:tol} and Lemma \ref{lem:infHL}. Hence, in the last step, we turn to  the estimates of $ \| \Lambda_h s^m\|$ and $\left\| \Lambda_h l^m\right\|$.  
To this aim,
we rewrite equations \eqref{Priori:cncfd:b} and \eqref{Priori:cncfd:d}   as
\begin{align}\label{apriori:est:L0}
d_1 \Lambda_h s_{i,j}^{m} = \mathcal{H}_h e_{i,j}^m + a_{11} \mathcal{H}_h s_{i,j}^{m} - \eta_{i,j}^{m},\quad
d_2 \Lambda_h l_{i,j}^{m} = \mathcal{H}_h w_{i,j}^m + a_{22} \mathcal{H}_h l_{i,j}^{m} - \mathcal{H}_h g(s_{i,j}^{m}) - \lambda_{i,j}^{m}.
\end{align}

Then, taking discrete inner products with themselves, we  obtain   %and utilizing \eqref{sta:uv:sum:tol}
\begin{equation}\label{apriori:est:Ll}
 \begin{aligned}
d_1^2\left\| \Lambda_h s^{m}\right\|^2 
&\leq 3\left\|\mathcal{H}_h e^m\right\|^2 + 3a_{11}^2 \left\|\mathcal{H}_h s^m\right\|^2 + 3\left\|\eta^m\right\|^2,
% \le  C  \left( \Gamma^0 + \left\|\eta^m\right\|^2+ \tau\sum_{n=1}^{m} \Psi^{n-1/2}\right), 
\end{aligned}
\end{equation}
and
\begin{equation}\label{apriori:est:L2}
	\begin{aligned}
	d_2^2\left\| \Lambda_h l^{m}\right\|^2 
	&\leq 4\left\|\mathcal{H}_h w^m\right\|^2 + 4a_{22}^2 \left\|\mathcal{H}_h l^m\right\|^2 + \frac{81L^2}{4}\left\|\mathcal{H}_h s^m\right\|^2
	+ 4\left\|\lambda^m\right\|^2.
%	&\le C \left( \Gamma^0 + \left\|\lambda^m\right\|^2+ \tau\sum_{n=1}^{m} \Psi^{n-1/2}\right),
\end{aligned}
\end{equation}
%where $C_8 = 3(1+a_{11}^2)C_7???$ and $C_9 = 4(1+a_{22}^2 + 81L^2/16)C_7??$.

Therefore, for $\tau$ sufficiently small,  inserting \eqref{apriori:est:L2} into \eqref{sta:uv:sum}, and an application of discrete Gr\"{o}nwall's inequality immediately implies  
\begin{equation}\label{sta:uv:sum:tol}
	\begin{aligned}
		\| \mathcal{H}_h s^m \|^2 +  \| \mathcal{H}_h e^{m} \| ^2
		+ \| \mathcal{H}_h l^m \|^2 + \| \mathcal{H}_h w^m \|^2  
		\leq  C \left( \Gamma^0 + \tau\sum_{n=1}^{m} \Psi^{n-1/2}\right).
	\end{aligned}
\end{equation}

Finally, combinations of \eqref{sta:uv:sum:tol} with \eqref{apriori:est:Ll}--\eqref{apriori:est:L2}, we get
\begin{equation}\label{apriori:result}
\begin{aligned}
		&\left\| \mathcal{H}_h s^m \right\|^2 + \left\| \Lambda_h s^m \right\|^2 + \left\|  e^m \right\| ^2
		+ \left\| \mathcal{H}_h l^m \right\|^2 + \left\| \Lambda_h l^m \right\|^2 + \left\|  w^m \right\|^2\\		
		&\quad \leq C \left(\Gamma^0 
		+  \left\|\eta^m\right\|^2 + \left\|\lambda^m\right\|^2		+ \tau\sum_{n=1}^{m} \Psi^{n-1/2}\right).
\end{aligned}
\end{equation}
%where $C_{10} = C_7 + C_8/d_1^2 + C_9/d_2^2$. 
Then, applying Lemma \ref{lem:infHL} to \eqref{apriori:result}, the proof of the theorem is completed.
\qed

Following the \textit{a priori} estimate, we can easily establish the stability result for the linearized CN-CFD scheme \eqref{bactre:cncfd}--\eqref{bactre:cncfd:ic}.
\begin{thm}\label{stability}
The linearized CN-CFD scheme \eqref{bactre:cncfd}--\eqref{bactre:cncfd:ic} is unconditionally stable with respect to the initial values in the sense that
\begin{align*}
		\|  U^m \| _\infty +\|  P^m \|+\|  V^m \| _\infty  + \|  Q^m \|		\leq C \left( \| U^0\| + \|  P^0\|	+ \|V^0\| + \|  Q^0\|\right),\quad 1\leq m\leq N,
\end{align*}
 where the constant $C$ is only related to the coefficients $d_1$, $a_{11}$, $a_{12}$, $d_2$, $a_{22}$, and the Lipschitz constant $L$.
\end{thm}
{\bf Proof.}  The  scheme \eqref{bactre:cncfd} can be viewed as taking $\{s^n, e^n, l^n, w^n\}=\{U^n, P^n, V^n, Q^n\}$ in \eqref{Priori:cncfd:a}--\eqref{Priori:cncfd:d} with $\xi^{n-1/2} = \zeta^{n-1/2} = \eta^{n}= \lambda^{n} = 0$.
Then, the conclusion can be directly derived from Theorem \ref{thm:priori}.\qed

\begin{rem}\rm
Actually, it is inferred from \eqref{bact:eq:bic} or \eqref{cncfd:iv} that  $\lbrace P^0, Q^0\rbrace $   can be controlled by initial conditions $\lbrace U^0, V^0\rbrace $. Thus, Theorem \ref{stability} indicates that the maximum-norm of $U^m$ and $V^m$  can be controlled by the discrete norms of the initial values $  U^0$ and $V^0 $.
\end{rem}

%Remark \ref{rem:re=cfd} and Theorem \ref{thm:priori} also imply the following convergence theorem, by giving that $g$ is twice continuously differentiable.
\begin{thm}\label{convergence}
 Let  $u(x, y, t), v(x, y, t) \in C^{3}( J;W^{2,\infty}(\Omega)) \cap C^{1}( J;H^6(\Omega)) $ be the exact solutions of the bacterial system \eqref{bact:eq} and $\left\{U^n, V^n \mid 0 \leq n \leq N\right\} \in\mathcal{V}_h^0$ be  the numerical solutions  of the linearized CN-CFD scheme \eqref{bactre:cncfd}--\eqref{bactre:cncfd:ic}. Then, the maximum-norm error estimates hold unconditionally  in the sense that
 \begin{align*}
 	\| u^m- U^m \| _\infty  +\|  v^m-V^m \| _\infty  	\leq \tilde{C}\left(\tau^2+h_x^4+h_y^4\right),\quad 1\leq m\leq N,
 \end{align*}
where the constant $\tilde{C}$ is related to the coefficients $d_1$, $a_{11}$, $a_{12}$, $d_2$, $a_{22}$, the Lipschitz constant $L$ and the bounds of the exact solutions $u$ and $v$.
\end{thm}
{\bf Proof.}
Let $\theta^n:= u^n - U^n$, $\sigma^n := p^n - P^n$, $\varphi^n := v^n - V^n$, and $\mu^n := q^n - Q^n$ for $0 \leq n \leq N$.
By subtracting those equations in \eqref{bactre:cncfd} from corresponding \eqref{ext:ut}--\eqref{ext:pxy} and \eqref{ext:vt}--\eqref{ext:qxy}, we derive the following error equations
\begin{subequations}\label{bactre:err}
	\begin{numcases}{}
	\bar{\sigma}_{i,j}^{n-1/2} = \delta_t \theta_{i,j}^{n-1/2} - a_{12} \varphi_{i,j}^{*,n} + R_{i,j}^{t,n-1/2}, 
	&$ (i, j) \in \omega,$\label{bactre:err:a}\\
	 - d_1 \Lambda_h \theta_{i,j}^{n} + a_{11} \mathcal{H}_h \theta_{i,j}^{n} 
	 = -\mathcal{H}_h \sigma_{i,j}^n + R_{i,j}^{s,n},
	&$ (i, j) \in \omega,$\label{bactre:err:b} \\
	\bar{\mu}_{i,j}^{n-1/2} = \delta_t \varphi_{i,j}^{n-1/2} + Z_{i,j}^{t,n-1/2}, 
	&$ (i, j) \in \omega,$\label{bactre:err:c}\\
	 -d_2 \Lambda_h \varphi_{i,j}^{n} + a_{22} \mathcal{H}_h \varphi_{i,j}^{n} 
	 - \mathcal{H}_h\left[ g(u_{i,j}^{n}) - g(U_{i,j}^{n})\right] = -\mathcal{H}_h \mu_{i,j}^n + Z_{i,j}^{s,n}, 
	&$ (i, j) \in \omega,$\label{bactre:err:d}
	\end{numcases}
\end{subequations}
for $1\leq n\leq N$. It is easy to check that $\theta_{i,j}^0= \sigma_{i,j}^0= \varphi_{i,j}^0= \mu_{i,j}^0=0$ for $(i,j)\in \bar{\omega}$.

Noting that \eqref{bactre:err} has basically the same form as \eqref{Priori:cncfd} with $\{s^n, e^n, l^n, w^n\}=\{\theta^n, \sigma^n, \varphi^n, \mu^n\}$, $\{ \xi^{n-1/2}, \eta^{n}, \zeta^{n-1/2}, \lambda^{n}\} = \{R^{t,n-1/2}, R^{s,n}, Z^{t,n-1/2}, Z^{s,n}\}$.
The only difference lies in the error term $\mathcal{H}_h\left[ g(u_{i,j}^{n})-g(U_{i,j}^{n})\right]$ in \eqref{bactre:err:d}, which is involved in the analysis of Steps II--III of Theorem \ref{thm:priori}. 
Therefore, we only need to pay special attention on the analysis in which this term is involved.

\paragraph{\bf Result of Step I}
   Following Step I of Theorem \ref{thm:priori}, a similar result as \eqref{sta:u:sum} can be obtained
\begin{equation}\label{cov:u:sum:tol}
	\begin{aligned}
		&  \| \mathcal{H}_h \theta^m \|^2 + \| \mathcal{H}_h  \sigma^m \|^2\\
		&\quad\leq  \| \mathcal{H}_h \theta^0 \|^2 + \| \mathcal{H}_h  \sigma^0 \|^2 
		+ C\tau\sum_{n=0}^{m-1}  \left(  \| \mathcal{H}_h \varphi^{n} \|^2  +  \|\Lambda_h \varphi^{n}\|^2\right) 
		+ C\tau\sum_{n=0}^{m}\left(  \| \mathcal{H}_h \theta^{n} \|^2 +  \| \mathcal{H}_h  \sigma^n \|^2 \right)  \\
		&\qquad + 	C\tau\sum_{n=1}^{m}\left( \| \mathcal{H}_h R^{t,n-1/2}\|^2  + \| \bar{R}^{s,n-1/2}\|^2
		+ \| \Lambda_h R^{t,n-1/2}\|^2 
		+ \|\delta_t R^{s,n-1/2}\|^2 \right). 
	\end{aligned}
\end{equation}
%where $C_2$ and $C_3$ are defined in \eqref{sta:u:sum}.

\paragraph{\bf Result of Step II} 
   
Given that $g$ is twice continuously differentiable and note that the numerical solution $\{U^n\mid 0\leq n\leq N\}$ is bounded (See Theorem \ref{stability}), we obtain 
   \begin{align}\label{}
   	\left\| \mathcal{H}_h (g(u^n) - g(U^n))\right\| 
   	&\leq C\,\| \mathcal{H}_h \theta^{n} \|, \label{est:g:err}\\
   	\| \mathcal{H}_h \delta_t \left( g(u) - g(U)\right)^{n-1/2}\| 
   	&\leq C\left( \| \mathcal{H}_h \theta^n \|  + \| \mathcal{H}_h \theta^{n-1} \| + \| \mathcal{H}_h \delta_t\theta^{n-1/2}\| \right). \label{est:gt:err}
   \end{align}
%where $C_g$ is related to the bound of twice derivative of $g$.
   
Throughout the analysis of Step II in Theorem \ref{thm:priori}, we observe that the estimate \eqref{sta:v:H:g} 
is now replaced by
\begin{equation}\label{cov:v:H:g}
	\begin{aligned}
		& \left(  \mathcal{H}_h [\bar{g}(u) - \bar{g}(U)]^{n-1/2},\mathcal{H}_h \bar{\varphi}^{n-1/2}\right)  
		\leq \|\mathcal{H}_h [\bar{g}(u) - \bar{g}(U)]^{n-1/2}\| \, \|\mathcal{H}_h \bar{\varphi}^{n-1/2}\| \\
		&\quad \leq  C\left(  \| \mathcal{H}_h \theta^{n} \|^2 +  \| \mathcal{H}_h \theta^{n-1} \|^2 +
		 \| \mathcal{H}_h \varphi^{n} \|^2 +  \| \mathcal{H}_h \varphi^{n-1} \|^2 \right),
	\end{aligned}
\end{equation}
where \eqref{est:g:err} is applied. Meanwhile, the estimate \eqref{sta:v:Lmid:gt} is changed to
\begin{equation}\label{cov:v:Lmid:gt}
	\begin{aligned}
		&  \left(  \mathcal{H}_h \delta_t [g(u) - g(U)]^{n-1/2},\mathcal{H}_h \bar{\mu}^{n-1/2}\right)   
		\leq \|\mathcal{H}_h \delta_t[g(u) - g(U)]^{n-1/2}\| \|\mathcal{H}_h \bar{\mu}^{n-1/2}\| \\
		&\quad \leq C\left( \sum_{\chi=\theta, \sigma, \mu} \left( \| \mathcal{H}_h \chi^n \|^2  + \| \mathcal{H}_h \chi^{n-1} \|^2\right) 
        +  \| \mathcal{H}_h \varphi^{*,n}\|^2
		+ \| \mathcal{H}_h R^{t,n-1/2}\|^2\right),
	\end{aligned}
\end{equation}
where \eqref{est:gt:err} is used.

Thus, the result corresponding to \eqref{sta:v:sum} in Step II of  Theorem \ref{thm:priori} is now as follows:  
\begin{equation}\label{cov:v:sum}
	\begin{aligned}
		 & \| \mathcal{H}_h \varphi^m \|^2 + \| \mathcal{H}_h \mu^m \|^2     \\
		&\quad \leq \| \mathcal{H}_h \varphi^0 \|^2 + \| \mathcal{H}_h \mu^0\|^2
		+ C\tau\sum_{n=1}^{m}\left(\| \mathcal{H}_h R^{t,n-1/2}\|^2   
		+ \| \mathcal{H}_h Z^{t,n-1/2}\|^2
		+ \|  \bar{Z}^{s,n-1/2}\|^2 \right.\\
		& \qquad \qquad \left.  +  \| \Lambda_h Z^{t,n-1/2}\|^2     
		+ \|\delta_t Z^{s,n-1/2}\|^2 \right)\\
		& 	\qquad	+ C \tau\sum_{n=0}^{m} \left(  \| \mathcal{H}_h \varphi^{n} \|^2 + \left\| \mathcal{H}_h \mu^n\right\|^2 
		+   \| \mathcal{H}_h \theta^{n} \|^2 +   \| \mathcal{H}_h  \sigma^n \|^2\right),
	\end{aligned}
\end{equation}
%where $C_{11} = \max\left\lbrace 5C_g, 1+a_{22}, d_2\right\rbrace $ and $C_{12} = \max\lbrace (9L/4 + 2+ 5a_{12}^2C_g, 2C_4, 5C_g \rbrace $. 

Now, adding up \eqref{cov:u:sum:tol} and \eqref{cov:v:sum}, we have the following result
\begin{equation}\label{cov:uv:sum:tolm}
	\begin{aligned}
		&\| \mathcal{H}_h \theta^m \|^2 + \| \mathcal{H}_h  \sigma^m \|^2
		+ \| \mathcal{H}_h \varphi^m \|^2 + \| \mathcal{H}_h \mu^m \|^2  \\
		&\quad\leq  \tilde \Gamma^0 + C\tau\sum_{n=1}^{m} \tilde \Psi^{n-1/2}
		+ C\tau\sum_{n=0}^{m} \left(  
		\| \mathcal{H}_h \theta^n \|^2 +    \| \mathcal{H}_h \sigma^{n}\|^2
		+ \| \mathcal{H}_h \varphi^n \|^2 +  \| \mathcal{H}_h \mu^{n} \|^2   + \left\|\Lambda_h \varphi^{n}\right\|^2 \right),
	\end{aligned}
\end{equation}
where 
%$C_{13} = 2(C_{11}+a_{11}+1)\exp\left( 2T(2C_2 + 2 C_3 + C_{12})\right) $, and 
\[\tilde\Gamma^0= \|  \theta^0\|^2 + \|  \sigma^0\|^2+ \|  \varphi^0\|^2 + \| \mu^0\|^2,
\]  %and 
\[\tilde\Psi^{n-1/2}= \sum_{\chi=\{R^t, Z^t\}} \left( \|  \chi^{n-1/2}\|^2 + \| \Lambda_h \chi^{n-1/2}\|^2  \right) 
+ \sum_{\chi=\{R^s, Z^s\}} \left( \| \bar{\chi}^{n-1/2}\|^2 + \|\delta_t \chi^{n-1/2}\|^2\right) .
\]

\paragraph{\bf Result of Step III} 
   In Step III of Theorem \ref{thm:priori}, the second equation of \eqref{apriori:est:L0} is changed to the following form 
\begin{align}
	d_2 \Lambda_h \varphi_{i,j}^{m} = \mathcal{H}_h \mu_{i,j}^m + a_{22} \mathcal{H}_h \varphi_{i,j}^{m} 
	- \mathcal{H}_h\left[ g(u_{i,j}^{m}) - g(U_{i,j}^{m})\right] - Z_{i,j}^{s,m}.
\end{align}
By a similar treatment as \eqref{apriori:est:L2} and using the inequality \eqref{est:g:err} we obtain
\begin{align*}
	d_2^2\left\| \Lambda_h \varphi^{m}\right\|^2 
	\leq 4\left\|\mathcal{H}_h \mu^m\right\|^2 + 4a_{22}^2 \left\|\mathcal{H}_h \varphi^m\right\|^2 + \frac{81L^2}{4}\left\|\mathcal{H}_h \theta^m\right\|^2
	+ 4\left\|Z^{s,m}\right\|^2.
\end{align*}
Inserting the above result into \eqref{cov:uv:sum:tolm}, for $\tau $ sufficiently small, the application of discrete Gr\"{o}nwall's inequality shows
\begin{equation}\label{cov:uv:sum:tol}
	\begin{aligned}
		\| \mathcal{H}_h \theta^m \|^2 + \| \mathcal{H}_h  \sigma^m \|^2
		+ \| \mathcal{H}_h \varphi^m \|^2 + \| \mathcal{H}_h \mu^m \|^2  
		\leq  C \left(\tilde \Gamma^0 + \tau\sum_{n=1}^{m} \tilde \Psi^{n-1/2}\right).  
	\end{aligned}
\end{equation}
Therefore, one can come to a similar result as \eqref{apriori:result}  just with different constant.
%Following the same procedure in the \textit{a priori} estimate, one can arrive at a similar conclusion as in Theorem \ref{thm:priori}.

Finally, applying the truncation errors \eqref{trun:Rs}, \eqref{trun:Zs} and the zero-valued initial errors to the similar result of \eqref{apriori:result}, we have
\begin{equation*}
\begin{aligned}
	\left\| \mathcal{H}_h \theta^m \right\|^2 + \left\| \Lambda_h \theta^m \right\|^2 
	+ \left\| \mathcal{H}_h \varphi^m \right\|^2 + \left\| \Lambda_h \varphi^m \right\|^2 
%\left\|\theta^n\right\|_{\infty}^2 + \left\|\varphi^n\right\|_{\infty}^2 
&\leq  C \left(\tilde{\Gamma}^0 +  \left\|R^{s,m}\right\|^2 + \left\|Z^{s,m}\right\|^2
	+ \tau\sum_{n=1}^{m} \tilde{\Psi}^{n-1/2}\right)\\
 &=\mo\left(\tau^4+h_x^8+h_y^8\right), \quad 1\leq m\leq N,
\end{aligned}
\end{equation*}
%where $C_{14} = C_{13} + C_8/d_1^2 + C_9/d_2^2$, 
which directly proves Theorem \ref{convergence} through Lemma \ref{lem:infHL}.  
\qed

%%%%%%%%%%%%%%%%%%%%%%%%%%%%%%%%%%%%%%%%%%%%%%%%%%%%%%%%%%%%%%%%
\section{A compact ADI algorithm and its numerical analysis}

We mention that although the developed linearized fourth-order compact difference scheme \eqref{bactre:cncfd} can be decoupled in practical computation via \eqref{bact:cncfd}--\eqref{bact:cncfd:t1}, at each time level $t_n$ ($1\leq n\leq N$), we still have to solve two linear algebra systems with dimensions $M:=M_x M_y$, which is computationally expensive for large-scale modeling and simulations. This motivates us to develop a linearized and decoupled fourth-order compact ADI scheme for the approximation of the bacterial system \eqref{bact:eq}, which can efficiently reduce the solution of the large-scale two-dimension problem into a series of small-scale one-dimensional problems. 
Actually, at each time level, only some tridiagonal  linear algebra systems with dimensions $M_x$ or  $M_y$ need to be solved.
Thus, it is much cheaper and more efficient, especially suitable for large-scale modeling and simulation of multi-dimensional problems.
Meanwhile, we shall also discuss the corresponding unconditional stability and error estimate via the discrete energy method and temporal-spatial error splitting technique. 

\subsection{ Derivation of the compact ADI scheme}
To construct a compact ADI scheme, we first recall from \eqref{ext:ut}--\eqref{ext:vt} by eliminating the auxiliary variables  $\left\lbrace p^n, q^n \right\rbrace $ that the original variables $\left\lbrace u^n, v^n\mid 1\leq n\leq N \right\rbrace \in\mathcal{V}_h^0$ satisfy  the following equations
\begin{subequations}\label{bact:cncfd+trun}
	\begin{numcases}{}
		\mathcal{H}_h \delta_t u_{i,j}^{n-1/2}  = d_1 \Lambda_h \bar{u}_{i,j}^{n-1/2} - a_{11} \mathcal{H}_h \bar{u}_{i,j}^{n-1/2} 
		+ a_{12} \mathcal{H}_h v_{i,j}^{*,n} + \mathcal{H}_h R_{i,j}^{t,n-1/2} + \bar{R}_{i,j}^{s,n-1/2},\label{bact:cncfd+trun:a}\\
		\mathcal{H}_h \delta_t v_{i,j}^{n-1/2}  = d_2 \Lambda_h \bar{v}_{i,j}^{n-1/2} - a_{22} \mathcal{H}_h \bar{v}_{i,j}^{n-1/2} 
		+  \mathcal{H}_h \left[ \bar{g}(u)\right]_{i,j}^{n-1/2} + \mathcal{H}_h Z_{i,j}^{t,n-1/2} + \bar{Z}_{i,j}^{s,n-1/2}, \label{bact:cncfd+trun:b}
	\end{numcases}
\end{subequations}
for $(i, j) \in \omega$ and $1\leq n\leq N$, where $\bar{R}^{s,n-1/2} = ( \bar{R}^{s,n}+ \bar{R}^{s,n-1}) /2 $, $\bar{Z}^{s,n-1/2} = ( \bar{Z}^{s,n}+ \bar{Z}^{s,n-1}) /2$, and $R^{t,n-1/2}$, $R^{s,n}$, $Z^{t, n-1/2}$ and $Z^{s,n}$ are defined and estimated in subsection \ref{subsec:compact}. 

Define two temporal second-order perturbation terms:
\begin{equation}\label{pertub:term}
	\begin{aligned}
		&R_{i,j}^{P, n-1/2}:=\frac{d_1^2 \tau^2}{4} \delta_x^2 \delta_y^2 \delta_t u_{i, j}^{n-1/2} 
		- \frac{d_1 a_{11}\tau^2}{8} \Lambda_h \delta_t u_{i, j}^{n-1/2}
		+ \frac{a_{11}^2\tau^2}{16} \mathcal{H}_h \delta_t u_{i, j}^{n-1/2},\\
		&Z_{i,j}^{P, n-1/2} := \frac{d_2^2 \tau^2}{4} \delta_x^2 \delta_y^2 \delta_t v_{i, j}^{n-1/2} 
		- \frac{d_2 a_{22}\tau^2}{8} \Lambda_h \delta_t v_{i, j}^{n-1/2}
		+ \frac{a_{22}^2\tau^2}{16} \mathcal{H}_h \delta_t v_{i, j}^{n-1/2}.
	\end{aligned}
\end{equation}
Then, by adding these small perturbation terms respectively to both sides of \eqref{bact:cncfd+trun:a} and \eqref{bact:cncfd+trun:b}, we further have
\begin{equation}\label{ADI:extu}
\begin{aligned}
	&\mathcal{H}_h \delta_t u_{i,j}^{n-1/2}  
	+ \frac{d_1 \tau^2}{4} \delta_x^2 \delta_y^2 \delta_t u_{i, j}^{n-1/2} 
	- \frac{d_1 a_{11}\tau^2}{8} \Lambda_h \delta_t u_{i, j}^{n-1/2}
	+ \frac{a_{11}^2\tau^2}{16} \mathcal{H}_h \delta_t u_{i, j}^{n-1/2} \\
	&\quad= d_1 \Lambda_h \bar{u}_{i,j}^{n-1/2} - a_{11} \mathcal{H}_h \bar{u}_{i,j}^{n-1/2} + a_{12} \mathcal{H}_h v_{i,j}^{*,n}  + R_{i,j}^{I, n-1/2},
\end{aligned}
\end{equation}
\begin{equation}\label{ADI:extv}
	\begin{aligned}
		&\mathcal{H}_h \delta_t v_{i,j}^{n-1/2}  
		+ \frac{d_2 \tau^2}{4} \delta_x^2 \delta_y^2 \delta_t v_{i, j}^{n-1/2} 
		- \frac{d_2 a_{22}\tau^2}{8} \Lambda_h \delta_t v_{i, j}^{n-1/2}
		+ \frac{a_{22}^2\tau^2}{16} \mathcal{H}_h \delta_t v_{i, j}^{n-1/2}   \\
		&\quad = d_2 \Lambda_h \bar{v}_{i,j}^{n-1/2} - a_{22} \mathcal{H}_h \bar{v}_{i,j}^{n-1/2} +  \mathcal{H}_h \left[ \bar{g}(u)\right]_{i,j}^{n-1/2} +Z_{i,j}^{I, n-1/2},
	\end{aligned}
\end{equation}
with truncation errors
\begin{equation}
	\begin{aligned}
	& R_{i,j}^{I, n-1/2} = \mathcal{H}_h R_{i,j}^{t,n-1/2} + \bar{R}_{i,j}^{s,n-1/2} + R_{i,j}^{P, n-1/2},
	&Z_{i,j}^{I, n-1/2} = \mathcal{H}_h Z_{i,j}^{t,n-1/2} + \bar{Z}_{i,j}^{s,n-1/2} + Z_{i,j}^{P, n-1/2}.
\end{aligned}
\end{equation}

Now,  neglecting the small truncation errors $R_{i,j}^{I, n-1/2}$ and $Z_{i,j}^{I, n-1/2}$ in \eqref{ADI:extu}--\eqref{ADI:extv}, and rearranging the resulting equations, a Crank-Nicolson type compact ADI finite difference (CN-ADI-CFD) method for system \eqref{bact:eq} can be defined as finding $ \{U^n, V^n \mid 1\leq n\leq N\} \in\mathcal{V}_h^0$ such that
\begin{subequations}\label{ADI:cfd}
	\begin{numcases}{}
		\Big( \mathcal{H}_{x} - \frac{d_1\tau}{2} \delta_x^2 + \frac{ a_{11}\tau}{4} \mathcal{H}_{x}\Big)
		\Big( \mathcal{H}_{y} - \frac{d_1\tau}{2} \delta_y^2 + \frac{a_{11}\tau}{4} \mathcal{H}_{y}\Big)U_{i,j}^{n}\nonumber\\
		\quad\quad= \Big( \mathcal{H}_{x} - \frac{d_1\tau}{2} \delta_x^2 + \frac{ a_{11}\tau}{4} \mathcal{H}_{x}\Big)
		\Big( \mathcal{H}_{y} - \frac{d_1\tau}{2} \delta_y^2 + \frac{a_{11}\tau}{4} \mathcal{H}_{y}\Big)U_{i,j}^{n-1} \label{ADI:cfd:a} \\
		\qquad\qquad\qquad +~d_1\tau \Lambda_h U_{i,j}^{n-1} - a_{11}\tau \mathcal{H}_h U_{i,j}^{n-1} + a_{12}\tau \mathcal{H}_h V_{i,j}^{*,n}
		=: F_{i,j}^{n-1}, &\quad $(i, j) \in \omega, $\nonumber\\
		\Big( \mathcal{H}_{x} - \frac{d_2\tau}{2} \delta_x^2 + \frac{a_{22}\tau}{4} \mathcal{H}_{x}\Big)
		\Big( \mathcal{H}_{y} - \frac{d_2\tau}{2} \delta_y^2 + \frac{a_{22}\tau}{4} \mathcal{H}_{y}\Big)V_{i,j}^{n}\nonumber\\
		\quad\quad= \Big( \mathcal{H}_{x} - \frac{d_2\tau}{2} \delta_x^2 + \frac{a_{22}\tau}{4} \mathcal{H}_{x}\Big)
		\Big( \mathcal{H}_{y} - \frac{d_2\tau}{2} \delta_y^2 + \frac{a_{22}\tau}{4} \mathcal{H}_{y}\Big)V_{i,j}^{n-1}\label{ADI:cfd:b}\\
		\qquad\qquad\qquad +~d_2 \tau\Lambda_h V_{i,j}^{n-1} - a_{22}\tau \mathcal{H}_h V_{i,j}^{n-1} +  \tau \mathcal{H}_h \bar{g}(U)_{i,j}^{n-1/2}=: G_{i,j}^{n-1}, &\quad $(i, j) \in \omega, $\nonumber
	\end{numcases}
\end{subequations}  
enclosed with initial conditions $ \{U^0, V^0\}$:
\begin{equation}\label{ADI:cncfd:ic}
		U_{i,j}^0 = u_0(x_i,y_j),\quad V_{i,j}^0 = v_0(x_i,y_j),	\quad (i, j) \in \omega.
\end{equation}

By introducing intermediate variables
\begin{equation*}
\begin{aligned}
	U_{i,j}^\star =	\Big( \mathcal{H}_{y} - \frac{d_1\tau}{2} \delta_y^2 + \frac{a_{11}\tau}{4} \mathcal{H}_{y}\Big)U_{i,j}^{n}, \quad
	V_{i,j}^\star =	\Big( \mathcal{H}_{y} - \frac{d_2\tau}{2} \delta_y^2 + \frac{a_{22}\tau}{4} \mathcal{H}_{y}\Big)V_{i,j}^{n},
\end{aligned}
\end{equation*}
the solutions $ \{U^n, V^n\} \in\mathcal{V}_h^0$ of the CN-ADI-CFD scheme \eqref{ADI:cfd} are reduced to solve two sequences of one-dimensional sub-problems along each spatial direction. To be specific, at each time level $t=t_n$,  the solutions $ \{U^n, V^n\}$ are determined as follows:

{\bf Step 1}: for each fixed $j \in \{1,2,\ldots M_y-1\}$, compute the intermediate solution $ U^\star=\{U_{i,j}^\star\}$ along  $x$-direction via
\begin{equation}\label{ADI:cfd:uxd}
		\Big( \mathcal{H}_{x} - \frac{d_1\tau}{2} \delta_x^2 + \frac{ a_{11}\tau}{4} \mathcal{H}_{x}\Big)U_{i,j}^\star
		= F_{i,j}^{n-1},\quad	U_{0, j}^\star = U_{M_x, j}^\star=0;		
\end{equation} 

{\bf Step 2}: for each fixed $ i \in \{1,2,\ldots M_x-1\}$, compute the solution $ U^n \in\mathcal{V}_h^0$ along $y$-direction via
\begin{equation}\label{ADI:cfd:uyd}
		\Big( \mathcal{H}_{y} - \frac{d_1\tau}{2} \delta_y^2 + \frac{a_{11}\tau}{4} \mathcal{H}_{y}\Big)U_{i,j}^{n} = U_{i,j}^\star;
\end{equation} 

{\bf Step 3}: once $ U^n$ is available, for each fixed $j \in \{1,2,\ldots M_y-1\}$, we compute the intermediate solution $ V^\star=\{V_{i,j}^\star\}$  along $x$-direction via
\begin{equation}\label{ADI:cfd:vxd}
	\Big( \mathcal{H}_{x} - \frac{d_2\tau}{2} \delta_x^2 + \frac{a_{22}\tau}{4} \mathcal{H}_{x}\Big)V_{i,j}^\star
	= G_{i,j}^{n-1}, \quad 		V_{0, j}^\star =  V_{M_x, j}^\star=0;		
\end{equation} 

{\bf Step 4}: finally, for each fixed $ i \in \{1,2,\ldots M_x-1\}$, compute the solution $ V^n \in\mathcal{V}_h^0$ along $y$-direction via
\begin{equation}\label{ADI:cfd:vyd}
		\Big( \mathcal{H}_{y} - \frac{d_2\tau}{2} \delta_y^2 + \frac{a_{22}\tau}{4} \mathcal{H}_{y}\Big)V_{i,j}^{n} = V_{i,j}^\star.
\end{equation} 

\begin{rem}\rm
	If the boundary conditions $u({\bm x}, t)=\phi_1({\bm x}, t) \neq 0$, $v({\bm x}, t)=\phi_2({\bm x}, t) \neq 0$ for ${\bm x} \in \partial \Omega$, then the boundary equations for \eqref{ADI:cfd:uxd} and \eqref{ADI:cfd:vxd} are respectively replaced by
	\begin{equation*}
		\begin{aligned}
		U_{0, j}^\star  = \Big( \mathcal{H}_{y} - \frac{d_1\tau}{2} \delta_y^2 + \frac{a_{11}\tau}{4} \mathcal{H}_{y}\Big) \phi_1\left(0, y_j, t_n\right), \quad
		U_{M_x, j}^\star   =\Big( \mathcal{H}_{y} - \frac{d_1\tau}{2} \delta_y^2 + \frac{a_{11}\tau}{4} \mathcal{H}_{y}\Big) \phi_1\left(1, y_j, t_n\right) ;
\\
		V_{0, j}^\star  = \Big( \mathcal{H}_{x} - \frac{d_2\tau}{2} \delta_x^2 + \frac{a_{22}\tau}{4} \mathcal{H}_{x}\Big) \phi_2\left(0, y_j, t_n\right), \quad
		V_{M_x, j}^\star   =\Big( \mathcal{H}_{x} - \frac{d_2\tau}{2} \delta_x^2 + \frac{a_{22}\tau}{4} \mathcal{H}_{x}\Big) \phi_2\left(1, y_j, t_n\right) .
	\end{aligned}
	\end{equation*} 
\end{rem}

\begin{rem}\rm
For the purpose of theoretical analysis below, we also introduce two auxiliary variables $\left\lbrace P^n,  Q^n \right\rbrace \in\mathcal{V}_h^0 $ such that  \eqref{ADI:cfd} can be rewritten into the following equivalent form 
\begin{subequations}\label{bactre:ADIcncfd}
	\begin{numcases}{}
		\bar{P}_{i,j}^{n-1/2} = \delta_t U_{i,j}^{n-1/2} - a_{12} V_{i,j}^{*,n},   &\qquad $(i, j) \in \omega, $ \label{bactre:ADIcncfd:a}\\
		- d_1 \Lambda_h \bar{U}_{i,j}^{n-1/2} + a_{11} \mathcal{H}_h \bar{U}_{i,j}^{n-1/2} 
		+ \frac{d_1 \tau^2}{4} \delta_x^2 \delta_y^2 \delta_t U_{i, j}^{n-1/2}   \nonumber\\
		\qquad- \frac{d_1 a_{11}\tau^2}{8} \Lambda_h \delta_t U_{i, j}^{n-1/2} 
		+ \frac{a_{11}^2\tau^2}{16} \mathcal{H}_h \delta_t U_{i, j}^{n-1/2}
		= - \mathcal{H}_h \bar{P}_{i,j}^{n-1/2},   &\qquad $(i, j) \in \omega, $  \label{bactre:ADIcncfd:b} \\
		\bar{Q}_{i,j}^{n-1/2} = \delta_t V_{i,j}^{n-1/2} ,    &\qquad $(i, j) \in \omega, $ \label{bactre:ADIcncfd:c}\\
		- d_2 \Lambda_h \bar{V}_{i,j}^{n-1/2} + a_{22} \mathcal{H}_h \bar{V}_{i,j}^{n-1/2} - \mathcal{H}_h [\bar{g}(U)]_{i,j}^{n-1/2}      
		+ \frac{d_2 \tau^2}{4} \delta_x^2 \delta_y^2 \delta_t V_{i, j}^{n-1/2}  \nonumber\\
		\qquad- \frac{d_2 a_{22}\tau^2}{8} \Lambda_h \delta_t V_{i, j}^{n-1/2}
		+ \frac{a_{22}^2\tau^2}{16} \mathcal{H}_h \delta_t V_{i, j}^{n-1/2}
		= -\mathcal{H}_h \bar{Q}_{i,j}^{n-1/2},    &\qquad $(i, j) \in \omega, $  \label{bactre:ADIcncfd:d}
	\end{numcases}
\end{subequations}
for  $1\leq n\leq N$,  enclosed with initial conditions 
\begin{equation}\label{bactre:ADIcncfd:ic}
	\begin{cases}{}
		U_{i,j}^0 = u_0(x_i,y_j),\quad V_{i,j}^0 = v_0(x_i,y_j),		&(i, j) \in \omega,\\
		\mathcal{H}_h P_{i,j}^0 = d_1 \Lambda_h U_{i,j}^{0} - a_{11} \mathcal{H}_h U_{i,j}^{0},   &(i, j) \in \omega,\\
		\mathcal{H}_h Q_{i,j}^0 = d_2 \Lambda_h V_{i,j}^{0} - a_{22} \mathcal{H}_h V_{i,j}^{0} + \mathcal{H}_h [g(U)]_{i,j}^{0},   &(i, j) \in \omega.
	\end{cases}
\end{equation}
\end{rem}

%%%%%%%%%%%%%%%%%%%%%%%%%%%%%%%%%%%
\subsection{Analysis of the compact ADI scheme}

Firstly, we prove the following \textit{a priori} estimate for the equivalent  form \eqref{bactre:ADIcncfd}--\eqref{bactre:ADIcncfd:ic} of the compact ADI scheme  \eqref{ADI:cfd}--\eqref{ADI:cncfd:ic}. 
\begin{thm}\label{a priori:ADI}
Assume that grid functions $\{s^n, e^n, l^n, w^n \mid 1 \leq n \leq N\}\in\mathcal{V}_h^0$ are the solutions of the following system with given initial values $\{ s^0, e^0, l^0, w^0\} $ and data $\{ \xi^{n-1/2},\eta^{n-1/2},\zeta^{n-1/2},\lambda^{n-1/2}\} $: 
\begin{subequations}\label{Priori:ADIcncfd}
	\begin{numcases}{}
		\bar{e}_{i,j}^{n-1/2} = \delta_t s_{i,j}^{n-1/2} - a_{12} l_{i,j}^{*,n} + \xi_{i,j}^{n-1/2},    &\qquad $ (i, j) \in \omega, $ \label{Priori:ADIcncfd:a}\\
		- d_1 \Lambda_h \bar{s}_{i,j}^{n-1/2} + a_{11} \mathcal{H}_h \bar{s}_{i,j}^{n-1/2} 
		 + \frac{d_1 \tau^2}{4} \delta_x^2 \delta_y^2 \delta_t s_{i, j}^{n-1/2} &\nonumber\\
		\qquad - \frac{d_1 a_{11}\tau^2}{8} \Lambda_h \delta_t s_{i, j}^{n-1/2}   
		+ \frac{a_{11}^2\tau^2}{16} \mathcal{H}_h \delta_t s_{i, j}^{n-1/2}
		= - \mathcal{H}_h \bar{e}_{i,j}^{n-1/2} + \eta_{i,j}^{n-1/2},    &\qquad $ (i, j) \in \omega, $ \label{Priori:ADIcncfd:b} \\
		\bar{w}_{i,j}^{n-1/2} = \delta_t l_{i,j}^{n-1/2} +\zeta_{i,j}^{n-1/2},   &\qquad $ (i, j) \in \omega, $  \label{Priori:ADIcncfd:c}\\
		- d_2 \Lambda_h \bar{l}_{i,j}^{n-1/2} + a_{22} \mathcal{H}_h \bar{l}_{i,j}^{n-1/2} - \mathcal{H}_h [{g}(s)]_{i,j}^{n-1/2}  
		 + \frac{d_2 \tau^2}{4} \delta_x^2 \delta_y^2 \delta_t l_{i, j}^{n-1/2}  & \nonumber\\
		\qquad- \frac{d_2 a_{22}\tau^2}{8} \Lambda_h \delta_t l_{i, j}^{n-1/2}    
		 + \frac{a_{22}^2\tau^2}{16} \mathcal{H}_h \delta_t l_{i, j}^{n-1/2}
		= -\mathcal{H}_h \bar{w}_{i,j}^{n-1/2} + \lambda_{i,j}^{n-1/2},   &\qquad $ (i, j) \in \omega, $ \label{Priori:ADIcncfd:d}
	\end{numcases}
\end{subequations}
for $1\leq n\leq N$, where $l^{*,n} $ is defined by  \eqref{not:timF}.
Then, there exists a positive constant $C$ which is only related to the coefficients $d_1$, $a_{11}$, $a_{12}$, $d_2$, $a_{22}$ and the Lipschitz constant $L$ such that
\begin{equation*}
	\left\|  s^m \right\| _\infty^2 + \left\|  l^m \right\| _\infty^2
	\leq C \left( \Gamma^0
	+ \tau\sum_{n=1}^{m} \Psi^{n-1/2} 
	+  \sum_{\chi=\{\eta, \lambda\}}\left(\|\chi^{1/2} \|^2 + \|\chi^{m-1/2} \|^2  + \tau\sum_{n=1}^{m-1} \| \triangle_t \chi^{n}\|^2\right)   \right) ,\quad 1\leq m\leq N,
\end{equation*}
where 
\begin{align}\label{aprioro:adi:Gamma}
\Gamma^0 =  \| \mathcal{H}_h s^0 \|^2 + \|\Lambda_h s^0\|^2 + \|\mathcal{H}_h l^0\|^2 + \|\Lambda_h l^0\|^2,
\end{align}
\begin{align}\label{aprioro:adi:Psi}
\Psi^{n-1/2} 
= \sum_{\chi=\{\xi, \zeta\}}\left( \| \mathcal{H}_h \chi^{n-1/2}\|^2 + \| \Lambda_h \chi^{n-1/2}\|^2\right)   
+ \|  \eta^{n-1/2}\|^2  
+ \|  \lambda^{n-1/2}\|^2 .
\end{align}
\end{thm}
{\bf Proof.} %It is seen from Lemma \ref{lem:infHL} that the maximum-norm estimate of $v$ can be obtained through analyzing $\left\| \mathcal{H}_h v\right\|$ and $\left\| \Lambda_h v\right\|$.  Thus, 
The proof procedure of this theorem is basically the same as that of Theorem \ref{thm:priori}, but is much more complicated. To make it more clear, we also divide the proof into three main steps.

\paragraph{\bf Step I. Estimates for $ \| \mathcal{H}_h s^n \|$ and $\left\| \Lambda_h s^n\right\|$} 

We shall first estimate $ \| \mathcal{H}_h s^n \|$. Taking discrete inner products with $ \mathcal{H}_h \delta_t s^{n-1/2}$ for \eqref{Priori:ADIcncfd:b}, we  obtain
\begin{equation}\label{Stab:ADI:e1}
\begin{aligned}
	& a_{11}\frac{ \| \mathcal{H}_h s^n \|^2 -  \| \mathcal{H}_h s^{n-1} \|^2}{2\tau}
	+ \frac{d_1 \tau^2}{4} \left( \delta_x^2 \delta_y^2 \delta_t s^{n-1/2}, \mathcal{H}_h \delta_t s^{n-1/2}\right) 
	+ \frac{a_{11}^2\tau^2}{16} \| \mathcal{H}_h \delta_t s^{n-1/2}\|^2  \\
	&\quad = d_1\left( \Lambda_h \bar{s}^{n-1/2}, \mathcal{H}_h \delta_t s^{n-1/2}\right) 
	+ \left(  \eta^{n-1/2}, \mathcal{H}_h \delta_t s^{n-1/2}\right) 
	- \left(  \mathcal{H}_h \bar{e}^{n-1/2}, \mathcal{H}_h \delta_t s^{n-1/2}\right) \\
	 & \qquad + \frac{d_1 a_{11}\tau^2}{8} \left(  \mathcal{H}_h \delta_t s^{n-1/2}, \Lambda_h \delta_t s^{n-1/2}\right). 
\end{aligned}
\end{equation}
Below, we also estimate \eqref{Stab:ADI:e1} term-by-term. First, for the second term on the left-hand side of \eqref{Stab:ADI:e1} and the last term on the right-hand side of \eqref{Stab:ADI:e1}, we use (i)--(ii) of Lemma \ref{lem:HLam} to arrive at
\begin{equation}\label{Stab:ADI:e2}
	\begin{aligned}
	\left( \delta_x^2 \delta_y^2 \delta_t s^{n-1/2}, \mathcal{H}_h \delta_t s^{n-1/2}\right) \ge 0,
%	\geq \frac{1}{3}\| \delta_x \delta_y \delta_t s^{n-1/2}\|_{xy}^2
%	     + \frac{h_x^2 h_y^2}{144}\| \delta_x^2 \delta_y^2 \delta_t s^{n-1/2}\|^2.
\quad 
	\left(  \mathcal{H}_h \delta_t s^{n-1/2}, \Lambda_h \delta_t s^{n-1/2}\right)  \leq 0. 
\end{aligned}
\end{equation}

Second, by acting the compact operator $ \mathcal{H}_h$ to  \eqref{Priori:ADIcncfd:a} to obtain  $ \mathcal{H}_h \delta_t s^{n-1/2}$, the first three terms on the right-hand side of \eqref{Stab:ADI:e1} can be estimated
\begin{equation}\label{Stab:ADI:e3}
	\begin{aligned}
	&d_1\left( \Lambda_h \bar{s}^{n-1/2}, \mathcal{H}_h \delta_t s^{n-1/2}\right) 
	+ \left(  \eta^{n-1/2}, \mathcal{H}_h \delta_t s^{n-1/2}\right) 
	- \left(  \mathcal{H}_h \bar{e}^{n-1/2}, \mathcal{H}_h \delta_t s^{n-1/2}\right)\\
	&\quad = d_1\left( \Lambda_h \bar{s}^{n-1/2}, \mathcal{H}_h \bar{e}^{n-1/2}\right) 
	+ d_1a_{12}\left( \Lambda_h \bar{s}^{n-1/2}, \mathcal{H}_h l^{*,n}\right) 
	- d_1\left( \Lambda_h \bar{s}^{n-1/2}, \mathcal{H}_h \xi^{n-1/2}\right)     \\
	&\qquad+\left( \eta^{n-1/2}, \mathcal{H}_h \bar{e}^{n-1/2}\right) 
	+ a_{12}\left( \eta^{n-1/2}, \mathcal{H}_h l^{*,n}\right) 
	- \left( \eta^{n-1/2}, \mathcal{H}_h \xi^{n-1/2}\right) \\
	&\qquad-\left( \mathcal{H}_h \bar{e}^{n-1/2}, \mathcal{H}_h \bar{e}^{n-1/2}\right) 
	- a_{12}\left( \mathcal{H}_h \bar{e}^{n-1/2}, \mathcal{H}_h l^{*,n}\right)     
	+ \left( \mathcal{H}_h \bar{e}^{n-1/2}, \mathcal{H}_h \xi^{n-1/2}\right) \\
	&\quad \leq -\frac{1}{2}\|  \mathcal{H}_h \bar{e}^{n-1/2}\|^2
		+C \left( \|  \Lambda_h \bar{s}^{n-1/2}\|^2
	+   \|  \mathcal{H}_h l^{*,n}\|^2
	+  \|  \mathcal{H}_h \xi^{n-1/2}\|^2
	+  \|  \eta^{n-1/2}\|^2\right). 
\end{aligned}
\end{equation}

Now, inserting \eqref{Stab:ADI:e2}--\eqref{Stab:ADI:e3} into \eqref{Stab:ADI:e1}, multiplying the resulting equation by $2\tau$, and then summing over $n$ from $1$ to $m$ (noting the definition of $l^{*,n}$ in \eqref{not:timF}), we reach the estimate for $ \| \mathcal{H}_h s^m\|$:
\begin{equation}\label{sta:ADIuab:innH:tol}
	\begin{aligned}
		&a_{11} \| \mathcal{H}_h s^m \|^2 + \tau\sum_{n=1}^{m}\|  \mathcal{H}_h \bar{e}^{n-1/2}\|^2  \\
		&  \le a_{11} \| \mathcal{H}_h s^0 \|^2
		 + C\tau\sum_{n=0}^{m-1}   \| \mathcal{H}_h l^{n} \| ^2 
		+ C\tau\sum_{n=1}^{m}\left(\|  \mathcal{H}_h \xi^{n-1/2}\|^2 + \|  \eta^{n-1/2}\|^2\right)
		+ C\tau\sum_{n=0}^{m}  \| \Lambda_h s^n \|^2,
	\end{aligned}
\end{equation}
which is similar to \eqref{sta:u:sum:tol}.

Next, by taking discrete inner products with $\Lambda_h \delta_t s^{n-1/2}$ for \eqref{Priori:ADIcncfd:b}, we estimate $\left\| \Lambda_h s^n\right\|$ as follows
\begin{equation}\label{sta:ADIub:innL}
	\begin{aligned}
		&d_1\frac{\| \Lambda_h s^n\|^2 -  \| \Lambda_h s^{n-1}\|^2}{2\tau}
		+ \frac{d_1a_{11}\tau^2}{8} \| \Lambda_h \delta_t s^{n-1/2}\|^2    \\
		&\quad = \left(  \mathcal{H}_h \bar{e}^{n-1/2}, \Lambda_h \delta_t s^{n-1/2}\right) 
		+ \frac{d_1 \tau^2}{4} \left(  \delta_x^2 \delta_y^2 \delta_t s^{n-1/2}, \Lambda_h \delta_t s^{n-1/2}\right)    
	    + \frac{a_{11}^2\tau^2}{16}\left( \mathcal{H}_h \delta_t s^{n-1/2}, \Lambda_h \delta_t s^{n-1/2}\right)  \\
	    &\qquad+a_{11}\left( \mathcal{H}_h \bar{s}^{n-1/2}, \Lambda_h \delta_t s^{n-1/2}\right)  
		- \left( \eta^{n-1/2}, \Lambda_h \delta_t s^{n-1/2}\right) \\
		&\quad \le \left(  \mathcal{H}_h \bar{e}^{n-1/2}, \Lambda_h \delta_t s^{n-1/2}\right) 
   +a_{11}\left( \mathcal{H}_h \bar{s}^{n-1/2}, \Lambda_h \delta_t s^{n-1/2}\right)  
- \left( \eta^{n-1/2}, \Lambda_h \delta_t s^{n-1/2}\right), 		
\end{aligned}
\end{equation}
where  (iii) of Lemma \ref{lem:HLam}  and \eqref{Stab:ADI:e2} are utilized to show that  the second and third terms on the right-hand side of
equality sign are less than zero.
%\begin{align}\label{sta:ADIub:innL:2r}
%	\left(  \delta_x^2 \delta_y^2 \delta_t s^{n-1/2}, \Lambda_h \delta_t s^{n-1/2}\right) 
%\leq  0.
%%-\frac{2}{3} \| \delta_x \delta_y^2 \delta_t s^{n-1/2} \|_x^2
%%	        -\frac{2}{3} \| \delta_x^2 \delta_y \delta_t s^{n-1/2} \|_y^2.
%\end{align}

Moreover,   by acting the compact difference operator $\Lambda_h$ to \eqref{Priori:ADIcncfd:a}, and taking discrete inner product with $\mathcal{H}_h \bar{e}^{n-1/2}$, we conclude from Lemma   \ref{lem:HLam}  that
\begin{equation}\label{sta:ADIua:innL}
	\begin{aligned}
 &	\left(  \Lambda_h \delta_t s^{n-1/2}, \mathcal{H}_h \bar{e}^{n-1/2}\right)  \\
	& \quad = \left(  \Lambda_h \bar{e}^{n-1/2} , \mathcal{H}_h \bar{e}^{n-1/2} \right) 
	+a_{12}\left(   \Lambda_h l^{*,n},\mathcal{H}_h \bar{e}^{n-1/2}\right)   
	- \left(  \Lambda_h \xi^{n-1/2} , \mathcal{H}_h \bar{e}^{n-1/2} \right) \\
	& \quad  \le a_{12}\left(   \Lambda_h l^{*,n},\mathcal{H}_h \bar{e}^{n-1/2}\right)   
	- \left(  \Lambda_h \xi^{n-1/2} , \mathcal{H}_h \bar{e}^{n-1/2} \right). 
	\end{aligned}
\end{equation}

Therefore, inserting \eqref{sta:ADIua:innL} into \eqref{sta:ADIub:innL}, multiplying the resulting equation by $2\tau$, and then summing over $n$ from $1$ to $m$, we can find
\begin{equation}\label{sta:ADIuab:innL:tol}
	\begin{aligned}
		d_1\| \Lambda_h s^m\|^2  
		&\leq  d_1\| \Lambda_h s^{0}\|^2	
		 + 2a_{12}\tau\sum_{n=1}^{m}\left(   \Lambda_h l^{*,n},\mathcal{H}_h\bar{e}^{n-1/2}\right) 
		 - 2\tau\sum_{n=1}^{m}\left(  \Lambda_h \xi^{n-1/2} , \mathcal{H}_h \bar{e}^{n-1/2} \right) \\
		&\qquad+ 2a_{11}\tau\sum_{n=1}^{m}\left( \mathcal{H}_h \bar{s}^{n-1/2}, \Lambda_h \delta_t s^{n-1/2}\right) 
		- 2\tau\sum_{n=1}^{m}\left( \eta^{n-1/2}, \Lambda_h \delta_t s^{n-1/2}\right) \\
		&\leq d_1\| \Lambda_h s^{0}\|^2 
		 + C\tau\sum_{n=1}^{m}\left( \| \Lambda_h l^{n-1}\| ^2 + \| \Lambda_h \xi^{n-1/2}\| ^2\right) 
		 + \frac{\tau}{8}\sum_{n=1}^{m}\|  \mathcal{H}_h \bar{e}^{n-1/2}\|^2\\
		&\qquad+ 2a_{11}\tau\sum_{n=1}^{m}\left( \mathcal{H}_h \bar{s}^{n-1/2}, \Lambda_h \delta_t s^{n-1/2}\right) 
		- 2\tau\sum_{n=1}^{m}\left( \eta^{n-1/2}, \Lambda_h \delta_t s^{n-1/2}\right) .
	\end{aligned}
\end{equation}

We still have to analyze the last two inner product terms in \eqref{sta:ADIuab:innL:tol}.  Following summation by parts in temporal direction, we have
\begin{equation}\label{sta:ADIuab:innL:1}
	\begin{aligned}
 &2a_{11}\tau\sum_{n=1}^{m} \left( \mathcal{H}_h \bar{s}^{n-1/2}, \Lambda_h \delta_t s^{n-1/2}\right) - 2\tau \sum_{n=1}^{m}\left( \eta^{n-1/2}, \Lambda_h \delta_t s^{n-1/2}\right)  \\
  &\quad =  -a_{11}\tau \sum_{n=1}^{m-1}\left( \mathcal{H}_h \delta_t {s}^{n+1/2} + \mathcal{H}_h \delta_t {s}^{n-1/2},  \Lambda_h  s^{n}\right) 
             - 2a_{11}\left( \mathcal{H}_h \bar{s}^{1/2}, \Lambda_h s^0\right)      
             + 2a_{11}\left(  \mathcal{H}_h \bar{s}^{m-1/2}, \Lambda_h s^{m} \right) \\
             & \qquad 
    +\tau \sum_{n=1}^{m-1}\left(  \triangle_t {\eta}^{n},  \Lambda_h  s^{n}\right) 
   +2\left(  \eta^{1/2}, \Lambda_h s^0\right)
     - 2\left(  \eta^{m-1/2}, \Lambda_h s^{m}\right)=:\sum_{i=1}^6 I_i.
	\end{aligned}
\end{equation}

For the $I_1$ term,  we get  $\mathcal{H}_h \delta_t {s}^{n-1/2} $ from  \eqref{Priori:ADIcncfd:a}, and then Cauchy-Schwarz inequality implies that
\begin{equation}\label{sta:ADIuab:innL:1:rhs1}
	\begin{aligned}
 I_1&\leq  a_{11}\tau \sum_{n=1}^{m}\|\mathcal{H}_h \delta_t {s}^{n-1/2} \| \left( \|\Lambda_h  s^{n-1} \| +\|\Lambda_h  s^{n} \| \right) \\
       &\leq a_{11}\tau \sum_{n=1}^{m}\left( \|\mathcal{H}_h \bar{e}^{n-1/2} \| 
                        + \|\mathcal{H}_h l^{*,n} \| 
                        + \|\mathcal{H}_h \xi^{n-1/2} \|\right)  \left( \|\Lambda_h  s^{n-1} \| +\|\Lambda_h  s^{n} \| \right)   \\
 &\leq C\tau \sum_{n=1}^{m}\|\mathcal{H}_h l^{n-1} \|^2
       + C\tau \sum_{n=1}^{m}\|\mathcal{H}_h \xi^{n-1/2} \|^2
       + C\tau \sum_{n=0}^{m}\|\Lambda_h s^{n} \|^2
       + \frac{\tau}{8}\sum_{n=1}^{m}\|\mathcal{H}_h \bar{e}^{n-1/2} \|^2.
	\end{aligned}
\end{equation}

For the $I_2$ term, we utilize the relation $\bar{s}^{1/2}=s^0+\frac{\tau}{2}\delta_t s^{1/2}$ as well as \eqref{Priori:ADIcncfd:a} for $\mathcal{H}_h \delta_t {s}^{1/2} $ again to see
\begin{equation}\label{sta:ADIuab:innL:1:rhs2}
	\begin{aligned}
		I_2  &= - 2a_{11}\left( \mathcal{H}_h s^0, \Lambda_h s^0\right) 
	   - a_{11}\tau\left( \mathcal{H}_h \delta_t s^{1/2}, \Lambda_h s^0\right) \\
	 &\leq a_{11}\left( \|\mathcal{H}_h s^{0} \|^2 + \|\Lambda_h s^{0} \|^2\right) 
	     + C\tau  \|\mathcal{H}_h l^{0} \|^2 
	     + C\tau \|\mathcal{H}_h \xi^{1/2} \|^2 
	     + C\tau \|\Lambda_h s^{0}\|^2
	     + \frac{\tau}{8} \|\mathcal{H}_h \bar{e}^{1/2} \|^2.
	\end{aligned}
\end{equation}

For the $I_3$ term,  we utilize the relation $\bar{s}^{m-1/2}=s^m-\frac{\tau}{2}\delta_t s^{m-1/2}$,  Lemma \ref{lem:HLam} and  \eqref{Priori:ADIcncfd:a} for $\mathcal{H}_h \delta_t {s}^{m-1/2} $ again to see
\begin{equation}\label{sta:ADIuab:innL:1:rhs3}
	\begin{aligned}
			I_3&
		=   2a_{11}\left(  \mathcal{H}_h s^m, \Lambda_h s^m\right) 
		- a_{11}\tau \left( \mathcal{H}_h \delta_t s^{m-1/2}, \Lambda_h s^m\right)  \\
		&\leq   C\tau   \|\mathcal{H}_h l^{*,m} \|^2 
		+ C\tau \|\mathcal{H}_h \xi^{m-1/2} \|^2 
		+ C\tau \|\Lambda_h s^{m} \|^2 
		+ \frac{\tau}{8} \|\mathcal{H}_h \bar{e}^{m-1/2} \|^2.
	\end{aligned}
\end{equation}

For the remainder terms,  Cauchy-Schwarz inequality directly implies that
\begin{equation}\label{sta:ADIuab:innL:1:rhs4}
	\begin{aligned}
		I_4+I_5+I_6	& \leq
		 \frac{d_1}{2}\| \Lambda_h s^m\|^2  
		+ \| \Lambda_h s^0\|^2 
		+\| \eta^{1/2}\|^2 
		+\frac{1}{d_1}\| \eta^{m-1/2}\|^2 
		+	2\tau\sum_{n=1}^{m-1} \left( \| \triangle_t\eta^{n}\|^2  
		+   \| \Lambda_h s^{n}\|^2\right) .
	\end{aligned}
\end{equation}

Therefore,  by combining \eqref{sta:ADIuab:innL:1:rhs1}--\eqref{sta:ADIuab:innL:1:rhs4}  with \eqref{sta:ADIuab:innL:1}, and then inserting them into  \eqref{sta:ADIuab:innL:tol}, we  get the estimate for $ \| \Lambda_h s^m\|$:
\begin{equation}\label{sta:ADIuab:innL:toltol}
	\begin{aligned}
		\frac{d_1}{2} \| \Lambda_h s^m\|^2  
		&\leq (d_1 +a_{11} + 1) \| \Lambda_h s^0\|^2 
		+ C\tau\sum_{n=0}^{m-1} \left(  \| \mathcal{H}_h l^{n} \| ^2 +\| \Lambda_h l^{n}\| ^2\right)		\\
		&\qquad 
		  + C\tau\sum_{n=1}^{m}\| \Lambda_h \xi^{n-1/2}\| ^2 
		  + C \left( \|\eta^{1/2} \|^2 + \|\eta^{m-1/2} \|^2  + \tau\sum_{n=1}^{m-1} \| \triangle_t \eta^{n}\| ^2\right) \\
		&\qquad 
		+ C \tau\sum_{n=0}^{m}\| \Lambda_h s^{n}\|^2
		+ \frac{\tau}{2}\sum_{n=1}^{m}\|  \mathcal{H}_h \bar{e}^{n-1/2}\|^2.
	\end{aligned}
\end{equation} 

As before, by combining \eqref{sta:ADIuab:innL:toltol} and \eqref{sta:ADIuab:innH:tol} together, we obtain
\begin{equation}\label{sta:ADIuab:innHL}
	\begin{aligned}
		&	a_{11} \| \mathcal{H}_h s^m \|^2 + \frac{d_1}{2}\| \Lambda_h s^m\|^2   
		 + \frac{\tau}{2}\sum_{n=1}^{m}\|  \mathcal{H}_h \bar{e}^{n-1/2}\|^2   \\
		&\quad \leq a_{11} \| \mathcal{H}_h s^{0} \|^2 + (d_1 +a_{11} + 1)\| \Lambda_h s^0\|^2
		+ C\tau\sum_{n=0}^{m-1} \left(  \| \mathcal{H}_h l^{n} \| ^2 +\| \Lambda_h l^{n}\| ^2\right)		 \\
		&\qquad + C\tau\sum_{n=1}^{m}\left( \| \mathcal{H}_h \xi^{n-1/2}\|^2+ \| \Lambda_h \xi^{n-1/2}\|^2   \right)\\
		&\qquad
		 +C\left( \|\eta^{1/2} \|^2 + \|\eta^{m-1/2} \|^2 
		+  \tau \sum_{n=1}^{m} \| \eta^{n-1/2}\|^2   +   \tau \sum_{n=1}^{m-1} \| \triangle_t \eta^{n}\| ^2\right)  
	 		+   C\tau\sum_{n=0}^{m}\| \Lambda_h s^{n}\|^2.  
	\end{aligned}
\end{equation}

From \eqref{sta:ADIuab:innHL}, we  can see that to complete  the maximum-norm estimate of $s^m$, i.e., the estimates $\|\mathcal{H}_h s^{m}\|$ and $  \| \Lambda_h s^{m} \| $, we still  need to analyze $\|\mathcal{H}_h l^{n}\|$ and $\|\Lambda_h l^{n}\|$  for $0\le n \le m-1$, which is also necessary for the  maximum-norm estimate of $l^m$.

\paragraph{\bf Step II. Estimates for  $ \| \mathcal{H}_h l^n \|$ and $\| \Lambda_h l^n\|$} 
First, by taking discrete inner products with $ \mathcal{H}_h \delta_t l^{n-1/2}$ for \eqref{Priori:ADIcncfd:d}, we  estimate $ \| \mathcal{H}_h l^n \|$ as follows
\begin{equation}\label{sta:ADIvb:innH}
	\begin{aligned}
	& a_{22}\frac{ \| \mathcal{H}_h l^n \|^2 -  \| \mathcal{H}_h l^{n-1} \|^2}{2\tau}
	+ \frac{d_2 \tau^2}{4}\left( \delta_x^2 \delta_y^2 \delta_t l^{n-1/2}, \mathcal{H}_h \delta_t l^{n-1/2}\right) 
	+ \frac{a_{22}^2\tau^2}{16} \| \mathcal{H}_h \delta_t l^{n-1/2}\|^2  \\
	&\quad = d_2\left( \Lambda_h \bar{l}^{n-1/2}, \mathcal{H}_h \delta_t l^{n-1/2}\right) 
	+ \left( \lambda^{n-1/2}, \mathcal{H}_h \delta_t l^{n-1/2}\right) 
	- \left( \mathcal{H}_h \bar{w}^{n-1/2}, \mathcal{H}_h \delta_t l^{n-1/2}\right)   \\
	& \qquad  
	+ \frac{d_2 a_{22}\tau^2}{8}\left( \Lambda_h \delta_t l^{n-1/2}, \mathcal{H}_h \delta_t l^{n-1/2}\right)
	+ \left( \mathcal{H}_h [\bar{g}(s)]^{n-1/2}, \mathcal{H}_h \delta_t l^{n-1/2}\right).
	\end{aligned}
\end{equation}

Compared to the estimates in Step I, we only need to analyze the last term in \eqref{sta:ADIvb:innH}.  
By Cauchy-Schwarz inequality and utilizing \eqref{Priori:ADIcncfd:c} as well as the assumption that $g(0)=0$, we   find
\begin{equation}\label{sta:ADIvb:innH:g}
	\begin{aligned}
	&\left( \mathcal{H}_h [\bar{g}(s)]^{n-1/2}, \mathcal{H}_h \delta_t l^{n-1/2}\right) 
	= \left( \mathcal{H}_h [\bar{g}(s)]^{n-1/2}, \mathcal{H}_h \bar{w}^{n-1/2}\right)  - \left( \mathcal{H}_h [\bar{g}(s)]^{n-1/2}, \mathcal{H}_h \zeta^{n-1/2}\right) \\
	&\quad \leq C \left(  \| \mathcal{H}_h s^n \|^2 +  \| \mathcal{H}_h s^{n-1} \|^2 \right) 
	+ C\| \mathcal{H}_h \zeta^{n-\frac{1}{2}}\|^2
	+ \frac{1}{8} \| \mathcal{H}_h \bar{w}^{n-\frac{1}{2}}\|^2 .
	\end{aligned}
\end{equation}
Then, following the same procedure of the estimate for $\left\| \mathcal{H}_h s^{m}\right\|$ in Step I, we  obtain 
\begin{equation}\label{sta:ADIvab:innH:tol}
	\begin{aligned}
	&a_{22} \| \mathcal{H}_h l^m \|^2 + \tau\sum_{n=1}^{m}\|  \mathcal{H}_h \bar{w}^{n-1/2}\|^2 \\
	&\quad \leq a_{22} \| \mathcal{H}_h l^{0}\|^2
	+ C\tau\sum_{n=0}^{m}   \|  \mathcal{H}_h s^{n}\|^2
	+ C\tau\sum_{n=1}^{m}\left( \|  \mathcal{H}_h \zeta^{n-1/2}\|^2 + \|  \lambda^{n-1/2}\|^2\right) 
	+  C\tau\sum_{n=0}^{m}\|  \Lambda_h l^{n}\|^2 .
	\end{aligned}
\end{equation}

Next,  we  estimate $\left\| \Lambda_h l^n\right\|$ in \eqref{sta:ADIvab:innH:tol} by taking discrete inner products with $ \Lambda_h \delta_t l^{n-1/2}$ for \eqref{Priori:ADIcncfd:d}.  After multiplying the resulting equation by $2\tau$ and then summing over $n$ from $1$ to $m$, we have
\begin{equation}\label{sta:ADIvb:innL}
	\begin{aligned}
	&d_2 \| \Lambda_h l^n\|^2 
	+ \frac{d_2a_{22}\tau^3}{4} \sum_{n=1}^{m}\| \Lambda_h \delta_t l^{n-1/2}\|^2    \\
	&\quad = d_2 \| \Lambda_h l^{0}\|^2
	+2\tau\sum_{n=1}^{m}\left( \mathcal{H}_h \bar{w}^{n-1/2}, \Lambda_h \delta_t l^{n-1/2}\right) 
	+ \frac{d_2 \tau^3}{2}\sum_{n=1}^{m}\left( \delta_x^2 \delta_y^2 \delta_t l^{n-1/2}, \Lambda_h \delta_t l^{n-1/2}\right) \\ 
	&\qquad+ \frac{a_{22}^2\tau^3}{8}\sum_{n=1}^{m}\left( \mathcal{H}_h \delta_t l^{n-1/2}, \Lambda_h \delta_t l^{n-1/2}\right)  
	+2a_{22}\tau\sum_{n=1}^{m}\left( \mathcal{H}_h \bar{l}^{n-1/2}, \Lambda_h \delta_t l^{n-1/2}\right) \\ 
	&\qquad- 2\tau\sum_{n=1}^{m}\left( \lambda^{n-1/2}, \Lambda_h \delta_t l^{n-1/2}\right) 
	- 2\tau\sum_{n=1}^{m}\left( \mathcal{H}_h [\bar{g}(s)]^{n-1/2}, \Lambda_h \delta_t l^{n-1/2}\right) .  
	\end{aligned}
\end{equation}

Similarly, compared to the estimate of $\left\| \Lambda_h s^n\right\|$ in Step I, we only have to pay special attention on the last term of  \eqref{sta:ADIvb:innL}. Following summation by parts in temporal direction, we have 
\begin{equation}\label{sta:ADIvb:innH:gt}
	\begin{aligned}
		-2\tau\sum_{n=1}^{m}\left( \mathcal{H}_h [\bar{g}(s)]^{n-1/2}, \Lambda_h \delta_t l^{n-1/2}\right) 
		&=   \sum_{n=1}^{m-1}\tau\left( \mathcal{H}_h \delta_t {g(s)}^{n+1/2} + \mathcal{H}_h \delta_t {g(s)}^{n-1/2},  \Lambda_h  l^{n}\right) \\
		&\quad +2\left( \mathcal{H}_h [\bar{g}(s)]^{1/2}, \Lambda_h l^0\right) 
		- 2\left( \mathcal{H}_h [\bar{g}(s)]^{m-1/2}, \Lambda_h l^{m}\right)\\
		&=:\sum_{i=1}^3 II_i.  
	\end{aligned}
\end{equation}
In fact, the terms $II_1$--$II_3$ could  be analyzed in a very similar way as $I_1$--$I_3$ (see  \eqref{sta:ADIuab:innL:1:rhs1}--\eqref{sta:ADIuab:innL:1:rhs3}) by 
 recalling \eqref{sta:v:Lmid:gtl2}. We directly list the result as below
 \begin{equation}\label{sta:ADIvb:innH:gt:rhs}
 	\begin{aligned}
 			 | II_1  | +| II_2  | +| II_3 |
 & \leq  \frac{d_2}{2} \|\Lambda_h l^{m} \|^2 
 + \frac{9L}{4}\left( \|\mathcal{H}_h s^{0} \|^2 + \|\Lambda_h l^{0} \|^2\right) 
 + \frac{81L^2}{32d_2}\|\mathcal{H}_h s^{m} \|^2
 + C\tau\sum_{n=1}^{m}\|\mathcal{H}_h \xi^{n-1/2} \|^2		\\
 & \quad + C\tau\sum_{n=0}^{m-1}\|\mathcal{H}_h l^{n} \|^2	
 + C\tau \sum_{n=0}^{m}\|\Lambda_h l^{n} \|^2	
 + \frac{\tau}{2}\sum_{n=1}^{m}\|\mathcal{H}_h \bar{e}^{n-1/2} \|^2.
 	\end{aligned}
\end{equation}

Therefore,  similar as the estimates for $\left\| \Lambda_h s^{m}\right\|$ in Step I, by inserting \eqref{sta:ADIvb:innH:gt}--\eqref{sta:ADIvb:innH:gt:rhs} as well as the estimate $\left\|\mathcal{H}_h s^{m} \right\|$ in   \eqref{sta:ADIuab:innH:tol} into \eqref{sta:ADIvb:innL},  we  obtain a similar result as \eqref{sta:ADIuab:innL:toltol}  for $ \| \Lambda_h l^m\|$:
\begin{equation}\label{sta:ADIvab:innL:toltol}
	\begin{aligned}
		\frac{d_2}{2}\| \Lambda_h l^m\|^2  
		&\leq C\left(  \| \Lambda_h l^{0}\|^2 
		+  \|  \mathcal{H}_h s^0\|^2\right) 
		+ C \tau\sum_{n=0}^{m} \left\| \Lambda_h s^{n}\right\| ^2\\
		&\qquad+ C\tau\sum_{n=1}^{m} \left( \| \mathcal{H}_h \xi^{n-1/2}\|^2
		       + \| \Lambda_h \zeta^{n-1/2}\|^2 
		       + \| \eta^{n-1/2}\|^2\right) \\
		&\qquad  + C\left( \|\lambda^{1/2} \|^2 + \|\lambda^{m-1/2} \|^2 
		+  \tau \sum_{n=1}^{m} \| \lambda^{n-1/2}\|^2   +   \tau \sum_{n=1}^{m-1} \| \triangle_t \lambda^{n}\| ^2\right)\\
		&\qquad
		       + C\tau\sum_{n=0}^{m-1}  \| \mathcal{H}_h l^{n} \| ^2
		       + C \tau\sum_{n=0}^{m}\left\| \Lambda_h l^{n}\right\| ^2 
		       +  \tau \sum_{n=1}^{m}\|  \mathcal{H}_h \bar{w}^{n-1/2}\|^2
		       + \frac{\tau}{2}\sum_{n=1}^{m}\|  \mathcal{H}_h \bar{e}^{n-1/2}\|^2.
	\end{aligned}
\end{equation} 
%and $C_8 = d_2+a_{22}+1+9L/4$, $C_9 = 9L/4 + 81L^2/32d_2$.
Similarly, by combining \eqref{sta:ADIvab:innH:tol} and \eqref{sta:ADIvab:innL:toltol} together, we obtain
\begin{equation}\label{sta:ADIvab:innHL}
	\begin{aligned}
		&	a_{22} \| \mathcal{H}_h l^m \|^2 + \frac{d_2}{2}\left\| \Lambda_h l^m\right\|^2 \\
%		+ \frac{\tau}{2}\sum_{n=1}^{m}\|  \mathcal{H}_h \bar{w}^{n-1/2}\|^2
%		+ \frac{\tau}{8}\sum_{n=1}^{m}\|  \mathcal{H}_h \bar{e}^{n-1/2}\|^2\\
		& \quad \leq a_{22} \| \mathcal{H}_h l^{0}\|^2 + C\left( \| \Lambda_h l^0\|^2
		+ \|  \mathcal{H}_h s^0\|^2 \right) 
		+ C\tau\sum_{n=0}^{m}\left( \left\|  \mathcal{H}_h s^n\right\|^2  +\left\| \Lambda_h s^n\right\|^2\right) \\
		&\qquad+ C\tau\sum_{n=1}^{m}\left( \| \mathcal{H}_h \zeta^{n-1/2}\|^2 
		+ \| \Lambda_h \zeta^{n-1/2}\|^2 
		+ \| \mathcal{H}_h \xi^{n-1/2}\|^2
		+ \| \eta^{n-1/2}\|^2\right)\\
		&\qquad + C\left( \|\lambda^{1/2} \|^2 + \|\lambda^{m-1/2} \|^2 
		+  \tau \sum_{n=1}^{m} \| \lambda^{n-1/2}\|^2   +   \tau \sum_{n=1}^{m-1} \| \triangle_t \lambda^{n}\| ^2\right)\\
		&\qquad
		+ C\tau\sum_{n=0}^{m}\left(\left\| \mathcal{H}_h l^{n}\right\|^2 + \left\| \Lambda_h l^{n}\right\|^2\right)
		+ \frac{\tau}{2}\sum_{n=1}^{m}\|  \mathcal{H}_h \bar{e}^{n-1/2}\|^2.
	\end{aligned}
\end{equation}
%which is similar to \eqref{sta:ADIuab:innHL}, where $C_{28} = \max\left\lbrace C_{21} + C_{27},C_{23} + 11/16\right\rbrace $ and $C_{29} = C_{26}+4$.

\paragraph{\bf Step III. Estimates for $\left\| s^m\right\|_\infty $ and $\left\| l^m\right\|_\infty $}

 Merging \eqref{sta:ADIuab:innHL} and \eqref{sta:ADIvab:innHL}, 
 %we have the following result
%\begin{equation*}
%	\begin{aligned}
%		&	\| \mathcal{H}_h s^m \|^2 + \| \Lambda_h s^m\|^2
%		+ \| \mathcal{H}_h l^m \|^2 +  \| \Lambda_h l^m\|^2 \\
%		&\quad\leq C\Gamma^0 + C \tau\sum_{n=1}^{m} \Psi^{n-1/2}
%		+ C\left( \|\eta^{1/2} \|^2 + \|\eta^{m-1/2} \|^2 
%		+   \tau \sum_{n=1}^{m-1} \| \triangle_t \eta^{n}\| ^2\right)\\
%		&\qquad + C\left( \|\lambda^{1/2} \|^2 + \|\lambda^{m-1/2} \|^2 
%		+   \tau \sum_{n=1}^{m-1} \| \triangle_t \lambda^{n}\| ^2\right)\\
%		&\qquad+ C\tau\sum_{n=0}^{m}\left(\| \mathcal{H}_h s^{n}\|^2 + \| \Lambda_h s^{n}\|^2
%		+  \|  \mathcal{H}_h l^n \|^2  + \| \Lambda_h l^n\|^2\right),
%	\end{aligned}
%\end{equation*} 
% and $C_{10} = \max \left\lbrace a_{11} + C_{25}, d_1 + a_{11} + 1, a_{22}, C_8\right\rbrace $.
and for  sufficiently small $\tau$, an application of discrete Gr\"{o}nwall's inequality implies that 
\begin{equation*}
	\begin{aligned}
		&	 \| \mathcal{H}_h s^m \|^2 + \| \Lambda_h s^m\|^2 
		+  \| \mathcal{H}_h l^m \|^2 + \| \Lambda_h l^m\|^2  \\   
		&\quad\leq C \left( \Gamma^0
		    + \tau\sum_{m=1}^{n} \Psi^{n-1/2}
		    +  \sum_{\chi=\{\eta, \lambda\}}\left(\|\chi^{1/2} \|^2 + \|\chi^{m-1/2} \|^2  + \tau\sum_{n=1}^{m-1} \| \triangle_t \chi^{n}\|^2\right) \right),    
	\end{aligned}
\end{equation*}
where $\Gamma^0$ and $\Psi^{n-1/2}$ are defined by \eqref{aprioro:adi:Gamma}--\eqref{aprioro:adi:Psi}.
Then, from Lemma \ref{lem:infHL} we  complete the proof of  the theorem. \qed

With the help of Theorem \ref{a priori:ADI}, we have the following stability conclusion.
\begin{thm}\label{stability:ADI}
	The linearized CN-ADI-CFD scheme \eqref{ADI:cfd}--\eqref{ADI:cncfd:ic} is unconditionally stable with respect to the initial values in the sense that
	\begin{align*}
		\|  U^m \| _\infty + \|  V^m \| _\infty  
		\leq C\left( \|\mathcal{H}_h U^0\| + \| \Lambda_h U^0\|
		+ \|\mathcal{H}_h V^0\| + \| \Lambda_h V^0\|\right),\quad 1\leq m\leq N,
	\end{align*}
 where the constant $C$ is related to the coefficients $d_1$, $a_{11}$, $a_{12}$, $d_2$, $a_{22}$, and the Lipschitz constant $L$.
\end{thm}
{\bf Proof.}
The  equivalent scheme \eqref{bactre:ADIcncfd} can be viewed as taking $\{s^n, e^n, l^n, w^n\}=\{U^n, P^n, V^n, Q^n\}$ in \eqref{Priori:ADIcncfd:a}--\eqref{Priori:ADIcncfd:d} with $\xi^{n-1/2} = \zeta^{n-1/2} = \eta^{n-1/2}= \lambda^{n-1/2} = 0$.
Then, Theorem \ref{a priori:ADI} immediately implies the result of Theorem \ref{stability:ADI}.\qed

Theorem \ref{a priori:ADI} also implies the following convergence estimate, which is presented below without rigorous proof.
\begin{thm}\label{convergence:ADI}
	Let $u(x, y, t), v(x, y, t) \in C^{3}( J;W^{2,\infty}(\Omega)) \cap C^{2}( J;W^{4,\infty}(\Omega)) \cap C^{1}( J;H^6(\Omega))$ be the exact solutions of the bacterial system \eqref{bact:eq} and $\left\{U^n, V^n \mid 0 \leq n \leq N\right\}$ be  the numerical solutions  of the linearized CN-ADI-CFD scheme \eqref{ADI:cfd}--\eqref{ADI:cncfd:ic}. Then, the maximum-norm error estimates hold unconditionally  in the sense that
	\begin{align*}
		\| u^m- U^m \| _\infty  +\|  v^m-V^m \| _\infty  	\leq \tilde{C}\left(\tau^2+h_x^4+h_y^4\right),\quad 1\leq m\leq N,
	\end{align*}
	where the constant $\tilde{C}$ is related to the coefficients $d_1$, $a_{11}$, $a_{12}$, $d_2$, $a_{22}$, the Lipschitz constant $L$ and the bounds of the exact solutions $u$ and $v$.
\end{thm}
{\bf Proof.} Again let $\theta^n:= u^n - U^n$, $\sigma^n := p^n - P^n$, $\varphi^n := v^n - V^n$, and $\mu^n := q^n - Q^n$ for $0 \leq n \leq N$.
We have the following error equations 
\begin{subequations}
	\begin{numcases}{}
		\bar{\sigma}_{i,j}^{n-1/2} = \delta_t \theta_{i,j}^{n-1/2} - a_{12} \varphi_{i,j}^{*,n} + R_{i,j}^{t,n-1/2}, &  $(i, j) \in \omega, $ \nonumber \\
		- d_1 \Lambda_h \bar{s}_{i,j}^{n-1/2} + a_{11} \mathcal{H}_h\bar{\theta}_{i,j}^{n-1/2} 
		 + \frac{d_1 \tau^2}{4} \delta_x^2 \delta_y^2 \delta_t \theta_{i, j}^{n-1/2}
		 - \frac{d_1 a_{11}\tau^2}{8} \Lambda_h \delta_t \theta_{i, j}^{n-1/2} \nonumber \\
		\qquad+ \frac{a_{11}^2\tau^2}{16} \mathcal{H}_h \delta_t \theta_{i, j}^{n-1/2}
		= - \mathcal{H}_h \bar{\sigma}_{i,j}^{n-1/2} + \bar{R}_{i, j}^{s,n-1/2} + R_{i, j}^{P,n-1/2},  & $(i, j) \in \omega, $ \nonumber \\
		\bar{\mu}_{i,j}^{n-1/2} = \delta_t \varphi_{i,j}^{n-1/2} +Z_{i, j}^{t,n-1/2}, & $(i, j) \in \omega, $ \nonumber \\
		- d_2 \Lambda_h \bar{\varphi}_{i,j}^{n-1/2} + a_{22} \mathcal{H}_h \bar{\varphi}_{i,j}^{n-1/2} - \mathcal{H}_h\left[ \bar{g}(u) - \bar{g}(U)\right]_{i,j}^{n-1/2}  
		 + \frac{d_2 \tau^2}{4} \delta_x^2 \delta_y^2 \delta_t \varphi_{i, j}^{n-1/2} \nonumber \\
		\qquad- \frac{d_2 a_{22}\tau^2}{8} \Lambda_h \delta_t \varphi_{i, j}^{n-1/2}
		 + \frac{a_{22}^2\tau^2}{16} \mathcal{H}_h \delta_t \varphi_{i, j}^{n-1/2}
		= -\mathcal{H}_h \bar{\mu}_{i,j}^{n-1/2} + \bar{Z}_{i, j}^{s,n-1/2} + Z_{i, j}^{P,n-1/2},  &  $(i, j) \in \omega, $ \nonumber
	\end{numcases}
\end{subequations}
for $1\leq n\leq N$. It is easy to check that $\theta_{i,j}^0= \varphi_{i,j}^0= 0$, $\sigma_{i,j}^0=R_{i,j}^{s,0}$, and $\mu_{i,j}^0=Z_{i,j}^{s,0}$ for $(i,j)\in \bar{\omega}$, and $R^{s,0}$ and $Z^{s,0}$ are defined by \eqref{ext:pxy} and \eqref{ext:qxy} for $n = 0$.

Taking $\{s^n, e^n, l^n, w^n\}=\{\theta^n, \sigma^n, \varphi^n, \mu^n\}$ in \eqref{Priori:ADIcncfd:a}--\eqref{Priori:ADIcncfd:d} with input $\{ \xi^{n-1/2}, \eta^{n-1/2}, \zeta^{n-1/2}, \lambda^{n-1/2}\} = \{R^{t,n-1/2}, \bar{R}^{s,n-1/2} + R^{P,n-1/2}, Z^{t,n-1/2}, \bar{Z}^{s,n-1/2} + Z^{P,n-1/2}\}$.  Following the analysis of Theorem \ref{a priori:ADI},  and notice that only the estimates of \eqref{sta:ADIvb:innH:g} and \eqref{sta:ADIvb:innH:gt} in Step II are different from the \textit{a priori} estimate and the convergence analysis, and they can be obtained in a similar fashion to those of Theorem \ref{convergence}. Thus, we finally get the following error estimate
\begin{equation}\label{conv:ADI}
	\begin{aligned}
		&\| \mathcal{H}_h \theta^m \|^2 + \| \Lambda_h \theta^m\|^2
		+ \| \mathcal{H}_h \varphi^m \|^2 +  \| \Lambda_h \varphi^m\|^2\\
		&\quad\leq C \left[\tilde\Gamma^0 
		+  \tau\sum_{n=1}^{m} \tilde\Psi^{n-1/2} 
		+  \sum_{\chi=\{ R^{P}, \ Z^{P}\}}\Bigg(\|\chi^{1/2} \|^2 + \|\chi^{m-1/2} \|^2  + \tau\sum_{n=1}^{m-1} \| \triangle_t \chi^{n}\|^2\Bigg) \right. \\
		&\qquad \qquad \left.  +\sum_{ \chi=\{R^{s}, \ Z^{s} \}}\Bigg(\|\bar{\chi}^{1/2} \|^2 + \|\bar{\chi}^{m-1/2} \|^2  + \tau\sum_{n=1}^{m-1} \| \triangle_t \bar{\chi}^{n}\|^2\Bigg)\right], 
	\end{aligned}
\end{equation}
where
\[
\tilde\Gamma^0
=  \| \mathcal{H}_h \theta^0 \|^2 + \|\Lambda_h \theta^0\|^2 + \|\mathcal{H}_h \varphi^0\|^2 + \|\Lambda_h \varphi^0\|^2,
\]  
\[
\tilde\Psi^{n-1/2} 
= \sum_{\chi=\{R^t, Z^t\}}\left( \| \mathcal{H}_h \chi^{n-1/2}\|^2 
+ \| \Lambda_h \chi^{n-1/2}\|^2\right)   
+ \sum_{\chi=\{ R^{P}, Z^{P}\} } \|  \chi^{n-1/2}\|^2
+ \sum_{\chi=\{R^{s} , Z^{s} \}} \|  \bar{\chi}^{n-1/2}\|^2. 
\]
Noting that
\begin{equation*}
	\| \triangle_t R^{P,n}\| 
	= 	\frac{1}{\tau}\|  R^{P,n+1/2}  - R^{P,n-1/2} \|
	\leq C\tau^2 ,
\end{equation*}
\begin{equation*}
	\| \triangle_t \bar{R}^{s,n} \| 
	= 	\frac{1}{\tau}\|  \bar{R}^{s,n+1/2}  - \bar{R}^{s,n-1/2} \| 
	\leq \|  \delta_t R^{s,n}\| + \|  \delta_t R^{s,n-1}\| 
	\leq C\left( \tau^2 + h_x^4 + h_y^4\right),
\end{equation*}
and similarly, $\| \triangle_t Z^{P,n} \| \leq C\tau^2 $ and $\| \triangle_t \bar{Z}^{s,n} \| \leq C\left( \tau^2 + h_x^4 + h_y^4\right)$. Combined  the truncation errors in \eqref{ext:ut}, \eqref{ext:t1}, \eqref{trun:Rs}, \eqref{ext:vt} and \eqref{trun:Zs} with \eqref{conv:ADI}, 
we complete the proof. \qed

%%%%%%%%%%%%%%%%%%%%%%%%%%%%%%%%%%%%%%%%%%%%%%%%%%%%%%%%%%%%%%%
\section{Numerical experiments}	
In this section, we conduct several numerical experiments using the linearized, decoupled CN-CFD method   \eqref{bactre:cncfd}--\eqref{bactre:cncfd:ic} via \eqref{bact:cncfd}--\eqref{bact:cncfd:t1}, \eqref{bactre:cncfd:a}  and \eqref{bactre:cncfd:c},
as well as the linearized, decoupled CN-ADI-CFD method  \eqref{ADI:cfd}--\eqref{ADI:cncfd:ic} via \eqref{ADI:cfd:uxd}--\eqref{ADI:cfd:vyd}, \eqref{bactre:ADIcncfd:a} and \eqref{bactre:ADIcncfd:c}.  
Accuracy and efficiency of the two algorithms are demonstrated in subsection \ref{test:ac}. Meanwhile,
 two types of asymptotic equilibrium are  simulated for different input initial values  in subsection \ref{test:asy}, in which the much more efficient CN-ADI-CFD algorithm is employed. 
For all the tests below, we set $M_x = M_y=M $ and $\tau = h^2$.
%---------------------------------------------------------------------
\subsection{Accuracy and efficiency test}\label{test:ac}	
In this subsection, we set $\Omega=(0,1)^2$, $J=(0,1]$, $g(u)=u^2/(1+u^2)$ and the other parameters $d_1$, $d_2$, $a_{11}$, $a_{12}$, $a_{21}$, $a_{22}$ are set to be 1 .
\begin{example}[\cite{SL22}]\label{exp:actest} 
	In this example, we consider the following general type of  model \eqref{bact:eq} with two additional fluxes:
	\begin{equation*}
		\begin{cases}
			u_{t}= \Delta u- u+ v + f_1(x,y,t), &\qquad  (\bm x, t) \in \Omega \times J, \\
			v_{t}= \Delta v- v+g(u) + f_2(x,y,t), &\qquad  (\bm x, t) \in \Omega \times J,
		\end{cases}
	\end{equation*}
where $f_1=0$, $f_2 = 4\pi^4e^{-t} \sin (\pi x) \sin (\pi y) -g(u(x,y,t))$, and the exact solutions $u(x, y, t)=e^{-t} \sin (\pi x) \sin (\pi y)$, $v(x, y, t)=2\pi^2e^{-t} \sin (\pi x) \sin (\pi y)$. 
\end{example}

\begin{table}[!htbp]\small
	\centering \caption{Errors and convergence orders  for the CN-CFD scheme \eqref{bactre:cncfd}--\eqref{bactre:cncfd:ic} at $T=1$. }\label{tb:CN_CFD}
	\setlength{\tabcolsep}{1.5mm}
	\begin{tabular}{c |c c c c c c c c c}
		\toprule
		$M$   &$\|U-u\|$     &Order     &$\|U-u\|_\infty$    &Order    &$\|V-v\|$   &Order     &$\|V-v\|_\infty$   &Order &CPU time\\
		\midrule	
		10	  &5.60e-06	     &--	    &1.12e-05	         &--	   &1.47e-04	     &--	 &2.93e-04	 &-- &0.34{\ \rm s}\\	 
		20	  &3.49e-07	     &4.00	    &6.98e-07	         &4.00	   &9.13e-06	     &4.00	 &1.83e-05	 &4.00 &1{\ \rm s}\\	 
		40	  &2.18e-08	     &4.00      &4.36e-08	         &4.00	   &5.70e-07	     &4.00	 &1.14e-06	 &4.00 &77{\ \rm s}\\	 
		80	  &1.36e-09	     &4.00	    &2.73e-09	         &4.00	   &3.56e-08	     &4.00	 &7.12e-08	 &4.00 &2{\ \rm h} 8{\ \rm m}\\	 
		160	  &8.53e-11	     &4.00	    &1.71e-10	         &4.00	   &2.23e-09	     &4.00	 &4.46e-09	 &4.00  & 7{\ \rm d} 3{\rm h} \\	
		\bottomrule
	\end{tabular}
	\end {table}	
\begin{table}[!htbp]\small
	\centering \caption{Errors and convergence orders for the CN-ADI-CFD scheme \eqref{ADI:cfd}--\eqref{ADI:cncfd:ic} at $T=1$. }\label{tb:ADI_CN_CFD}
	\setlength{\tabcolsep}{1.5mm}
	\begin{tabular}{c |c c c c c c c c c}
		\toprule
		$M$   &$\|U-u\|$     &Order     &$\|U-u\|_\infty$    &Order    &$\|V-v\|$   &Order     &$\|V-v\|_\infty$   &Order &CPU time\\
		\midrule	
%		10	  &4.73e-05	     &--	    &9.46e-05	         &--	   &4.09e-04	     &--	 &8.18e-04	 &--   &0.02{\ \rm s} \\	 
		20	  &2.96e-06	     &--	    &5.91e-06	         &--	   &2.56e-05	     &--	 &5.11e-05	 &--  &0.06{\ \rm s} \\	 
		40	  &1.85e-07	     &4.00      &3.69e-07	         &4.00	   &1.60e-06	     &4.00	 &3.19e-06	 &4.00  &0.6{\ \rm s} \\	 
		80	  &1.15e-08	     &4.00	    &2.31e-08	         &4.00	   &9.98e-08	     &4.00	 &2.00e-07	 &4.00  &10{\ \rm s} \\
		160	  &7.22e-10	     &4.00	    &1.44e-09	         &4.00	   &6.24e-09	     &4.00	 &1.25e-08	 &4.00  &3 {\ \rm m}\\
		280	  &7.69e-11	     &4.00	    &1.54e-10	         &4.00	   &6.65e-10	     &4.00	 &1.33e-09	 &4.00  &36 {\ \rm m}\\	 
		\bottomrule
	\end{tabular}
	\end {table}		

The purpose of this example is to test the accuracy and efficiency of the linearized CN-CFD and CN-ADI-CFD methods. Table \ref{tb:CN_CFD} shows the errors and convergence orders of the discrete $L^2$ norm and maximum-norm for the average concentrations of the bacteria $u$ and the infective people $v$ using the CN-CFD method \eqref{bactre:cncfd}--\eqref{bactre:cncfd:ic}. Besides, the CPU time consumed by the CN-CFD method  is also presented in Table \ref{tb:CN_CFD}.
Meanwhile, Table \ref{tb:ADI_CN_CFD} displays the corresponding numerical errors and convergence orders of $u$ and $v$, as well as the CPU running time for the CN-ADI-CFD method \eqref{ADI:cfd}--\eqref{ADI:cncfd:ic}. In summary, the following conclusions can be observed:
\begin{enumerate}[(i)]
	\item %Both methods produce the same order of magnitude errors for $u$ and $v$ measured in the discrete $L^2$ norm and maximum norm. Furthermore, 
	Both methods can generate numerical results with second-order convergence in time and fourth-order convergence in space, which are consistent with the conclusions in Theorems \ref{convergence} and \ref{convergence:ADI} very well.
	\item With the same order of magnitude errors for $u$ and $v$ measured in the discrete norms, the CN-ADI-CFD method is significantly more efficient than the CN-CFD method. 	For instance, the CN-CFD scheme consumes more than $7$ days  to attain the order of magnitude $1.0e-10$ for the maximum-norm error  of $u$. However,  the CN-ADI-CFD method  costs only  half an hour! This indicates a substantial reduction in CPU time consumption for the efficient  CN-ADI-CFD method.
\end{enumerate}	

%---------------------------------------------------------------------
\subsection{Asymptotic equilibrium}	\label{test:asy}
Real data indicates that a bistable scenario is most likely to occur when the initial conditions are relevant. Large outbreaks tend to result in a nontrivial endemic state that invades the entire habitat. However, small outbreaks tend to lead to extinction and may remain strictly localized in space (see \cite{CVM,CW}). This could explain why, despite our exposure to numerous infections, only some diseases have evolved into an endemic state. Here, we use the CN-ADI-CFD scheme \eqref{ADI:cfd}--\eqref{ADI:cncfd:ic}  to simulate the two different states caused by distinct initial values.
\begin{table}[!htbp]\small
	\centering \caption{Errors and convergence orders for the CN-ADI-CFD scheme \eqref{ADI:cfd}--\eqref{ADI:cncfd:ic}  at $T=1$. }\label{ext:ADI_CN_CFD}
	\setlength{\tabcolsep}{1.5mm}
	\begin{tabular}{c |c c c c c c c c}
		\toprule
		$M$   &$\|e_U\|$     &Order     &$\|e_U\|_\infty$    &Order    &$\|e_V\|$   &Order     &$\|e_V\|_\infty$   &Order \\
		\midrule	
		20	  &3.51e-03	     &--	    &1.69e-02	         &--	   &1.90e-03	     &--	 &1.01e-02	 &--   \\	 
		40	  &1.30e-04	     &4.76      &7.41e-04	         &4.51	   &1.13e-04	     &4.06	 &5.17e-04	 &4.29   \\	 
		80	  &7.35e-06	     &4.14	    &4.40e-05	         &4.07	   &5.88e-06	     &4.27	 &3.14e-05	 &4.04   \\
		160	  &4.51e-07	     &4.03	    &2.71e-06	         &4.02	   &3.55e-07	     &4.05	 &1.94e-06	 &4.02  \\	 
		\bottomrule
	\end{tabular}
	\end {table}
	
	\begin{example}\label{exam:noise}
 In this example, we choose $g(u)=u^2/(1+u^2)$, $d_1=0.001$, $a_{11}=1.0$, $a_{12}=2.5$, $d_2=0.0001$, $a_{22}=1.0$, and initial data
	\begin{align*}
		&u_0(x, y)=0.3 \sin \pi x \sin \pi y+n(x, y), \quad
		v_0(x, y)=g\left(u_0 \right),
	\end{align*}
	where $n(x, y)>0$ is the "noisy" function defined by
	\[n(x, y):=\frac{1}{4} \sin ^2 5 \pi x \sin ^2 4 \pi y+\frac{1}{4} \sin ^2 3 \pi x \sin ^2 7 \pi y+\frac{1}{10} \sin ^2 9 \pi x \sin ^2 11 \pi y.\]
\end{example}
	
\begin{figure}[!htbp]\small
	\centering
	\includegraphics[width=0.4\textwidth]{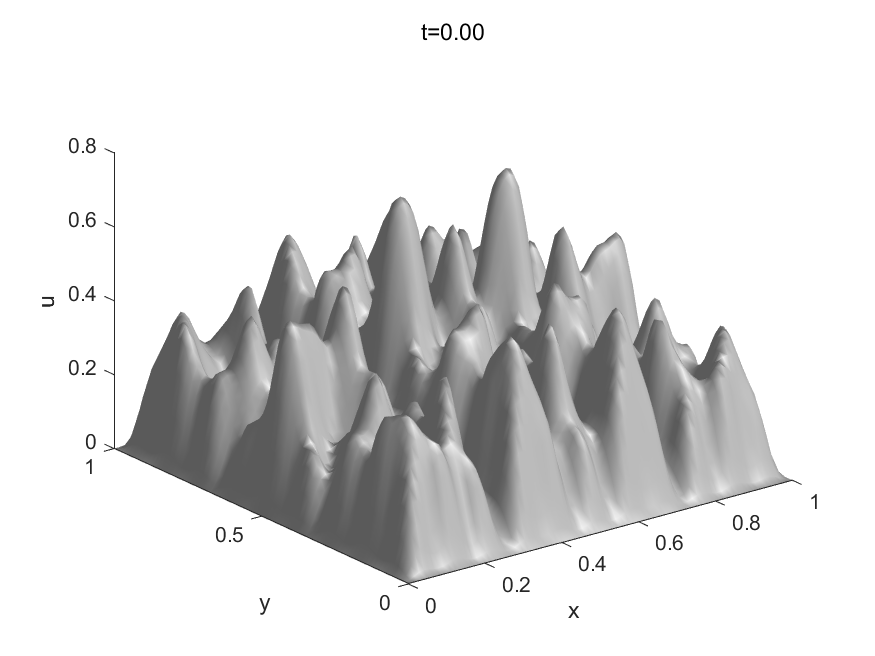}
	\includegraphics[width=0.4\textwidth]{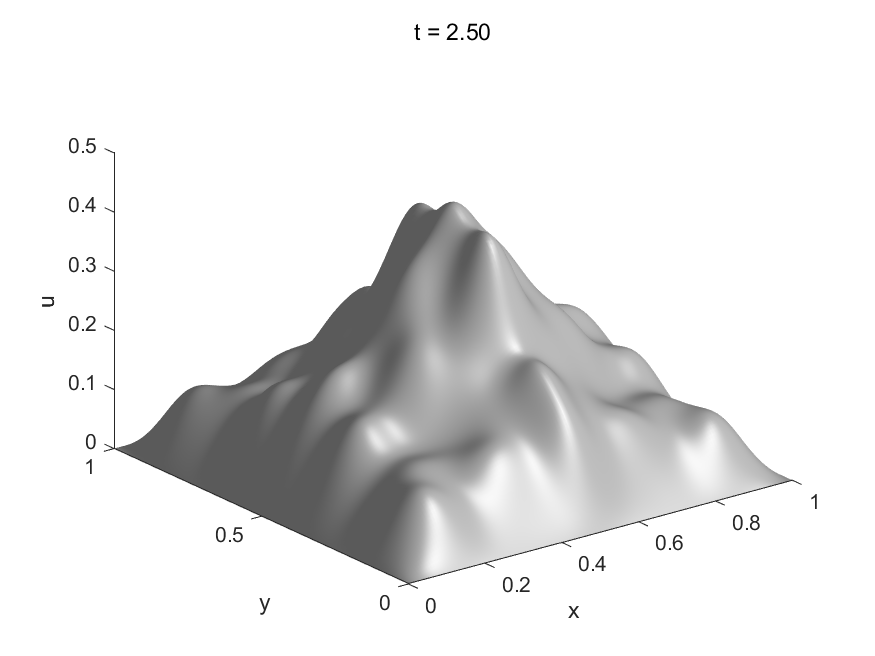}
	\includegraphics[width=0.4\textwidth]{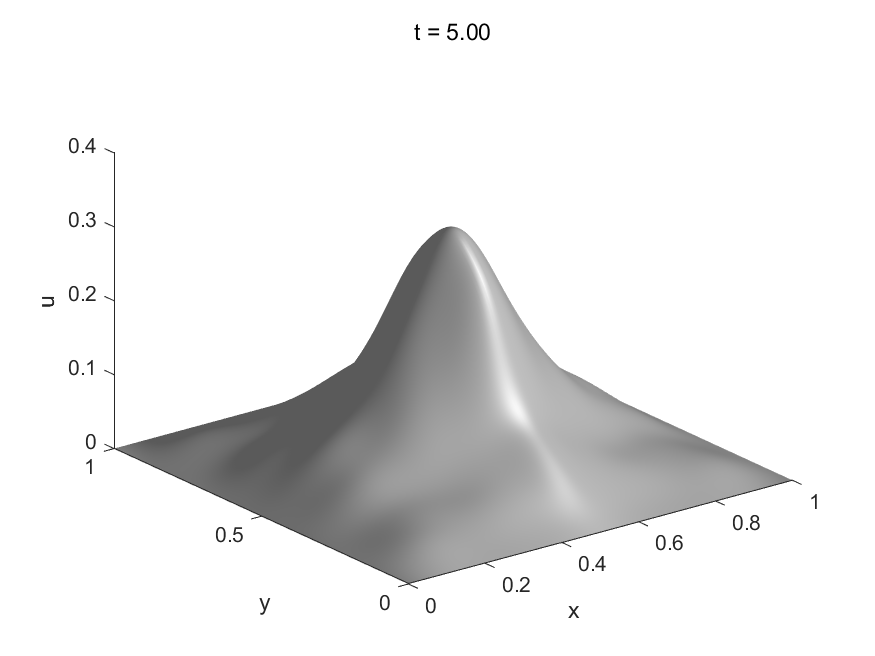}
	\includegraphics[width=0.4\textwidth]{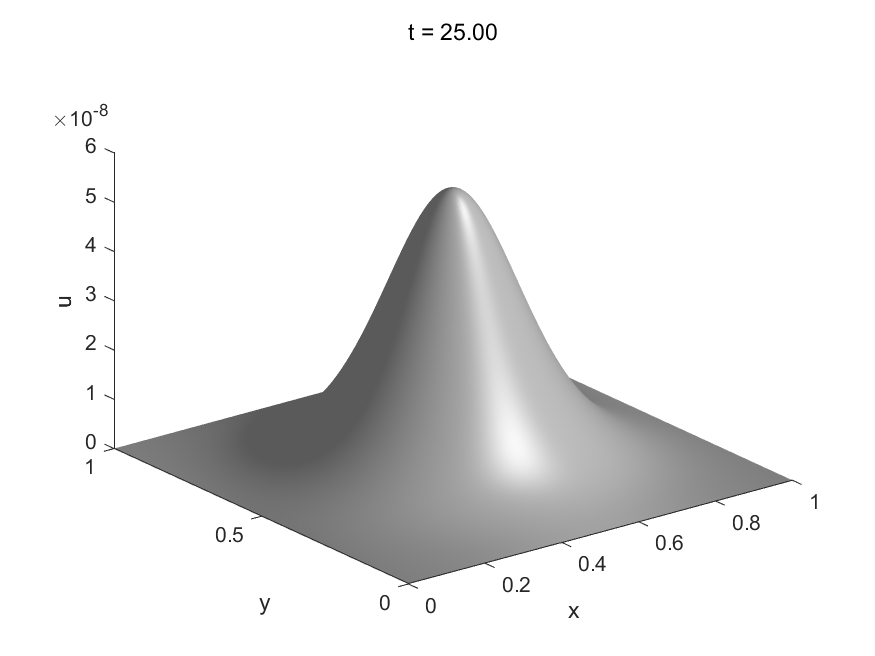}
	\caption{Extinction phenomenon of the average concentration of bacteria $u$.}\label{fig:u-}
\end{figure}

\begin{figure}[htbp]\small
	\centering
	\includegraphics[width=0.4\textwidth]{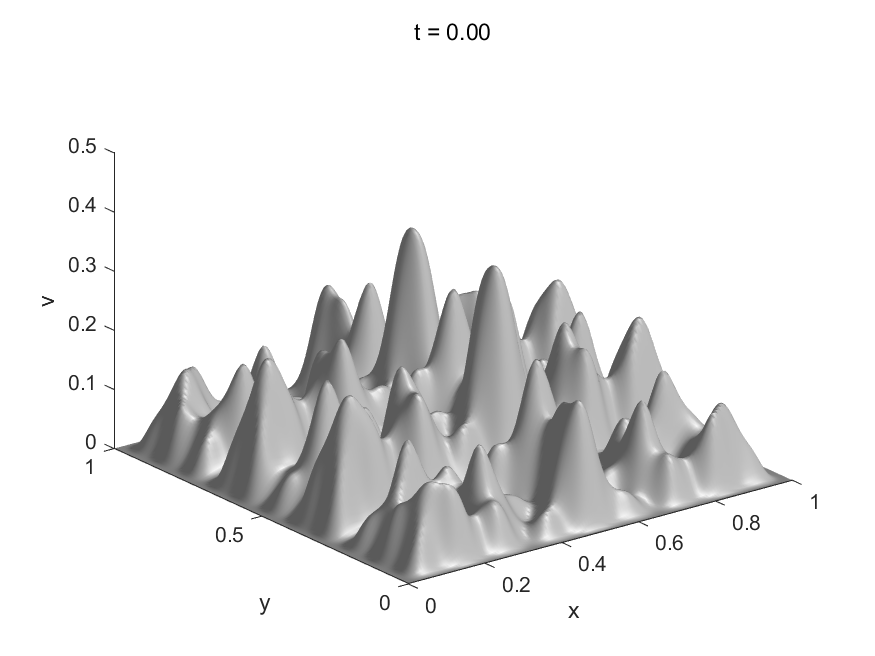}
	\includegraphics[width=0.4\textwidth]{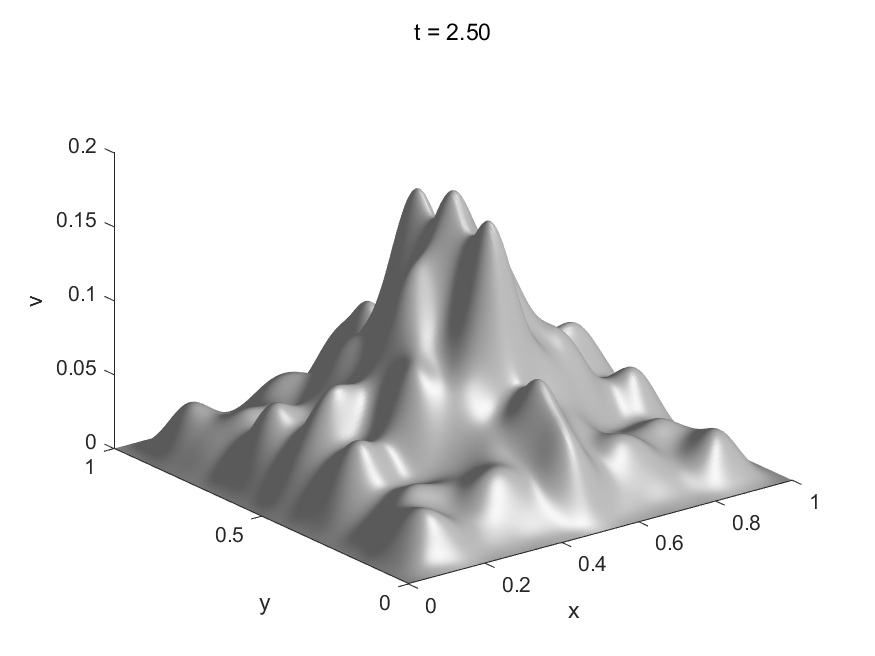}
	\includegraphics[width=0.4\textwidth]{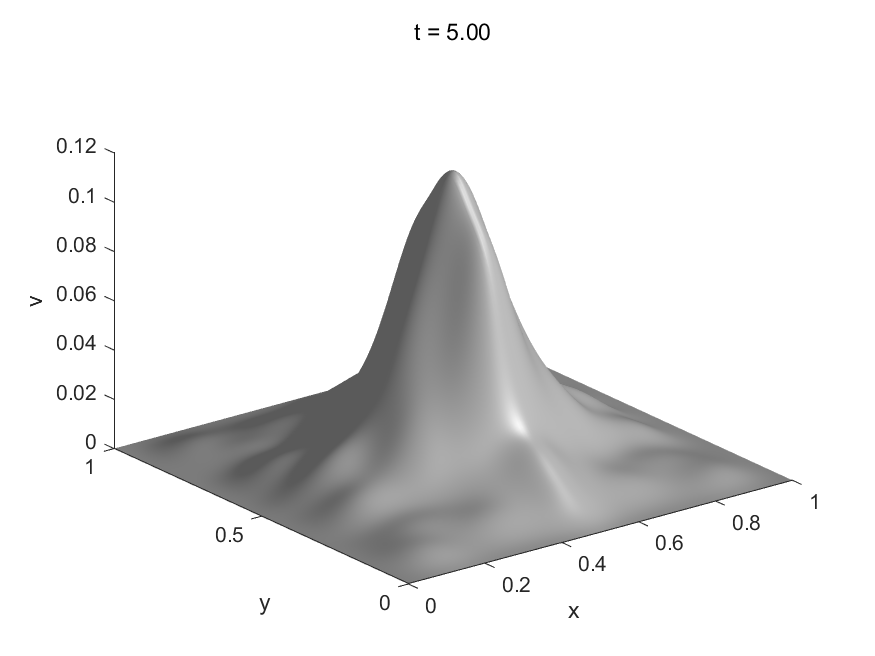}
	\includegraphics[width=0.4\textwidth]{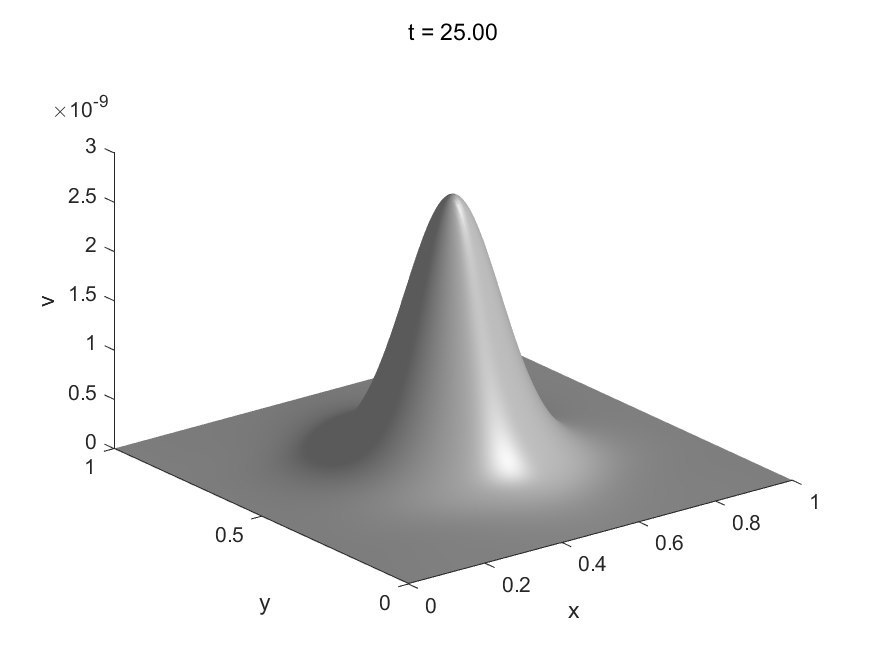}
	\caption{Extinction phenomenon of the average concentration of infective people $v$.}\label{fig:v-}
\end{figure}

 Since no exact solutions are available for this example, the errors are  measured by the so-called Cauchy error \cite{DFW}, which  is calculated between solutions obtained on two different grid pairs $\left( \tau, h\right) $ and $\left( \tau/4,h/2\right) $  by $\left\|e_U\right\|=\|U_h^\tau-U_{h / 2}^{\tau/4}\|$, where $U_h^\tau$ denotes the numerical solution under temporal stepsize $\tau$ and spatial stepsize $h$. 
The results of the discrete $L^2$ and $L^\infty$ Cauchy errors are presented in Table \ref{ext:ADI_CN_CFD}, and from which we can also confirm the second-order temporal and fourth-order spatial convergences for both the variables $u$ and $v$. 
Besides, the evolutions of the bacterial concentration $u$ and the concentration of the infective humans $v$ are respectively  illustrated in Figures \ref{fig:u-}--\ref{fig:v-} with $\tau = 0.1$ and $M=160$. 
For this chosen initial data, the bacterial concentration $u$ is attracted to zero.
Ecologically, in this situation, the bacteria will die out eventually and the humans are safe.
We can also see that the localization of the epidemic becomes evident as it is extinguished.

\begin{example}\label{exam:endemic}
	In this example, most of the parameters are exactly  the same as Example \ref{exam:noise}, except that the initial data
	\begin{align*}
		u_0(x, y)=0.5 \sin \pi x \sin \pi y+n(x, y), \quad v_0(x, y)=g\left(u_0 \right),
	\end{align*}
is chosen a little larger.
\end{example}

\begin{figure}[!htbp]\small
	\centering
	\includegraphics[width=0.32\textwidth]{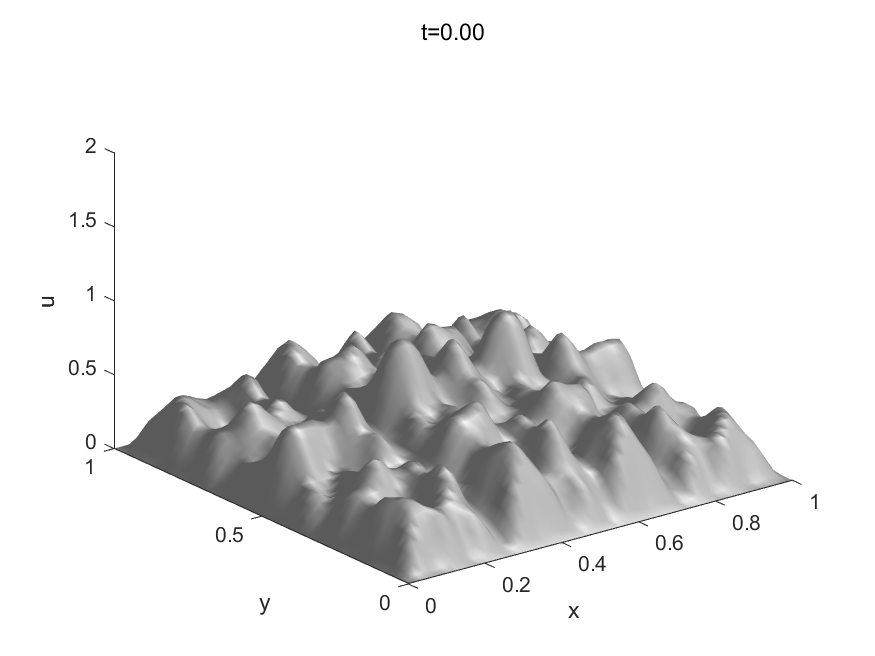}
	\includegraphics[width=0.32\textwidth]{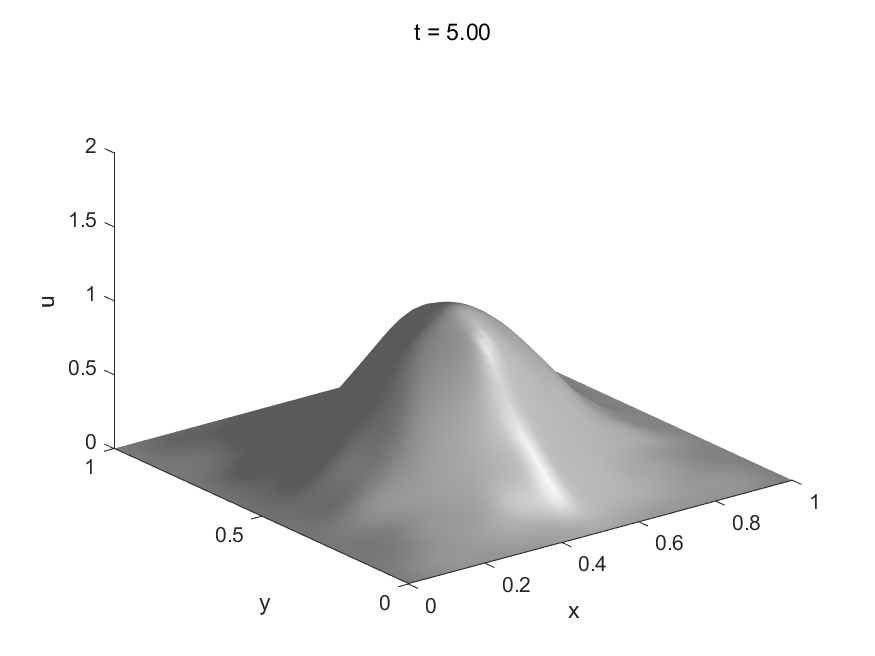}
	\includegraphics[width=0.32\textwidth]{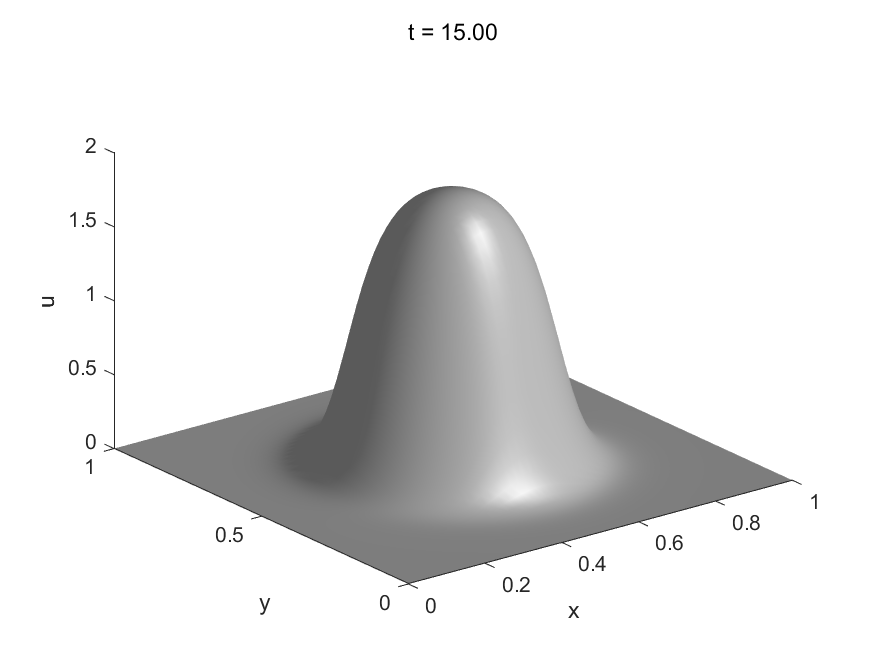}
	\includegraphics[width=0.32\textwidth]{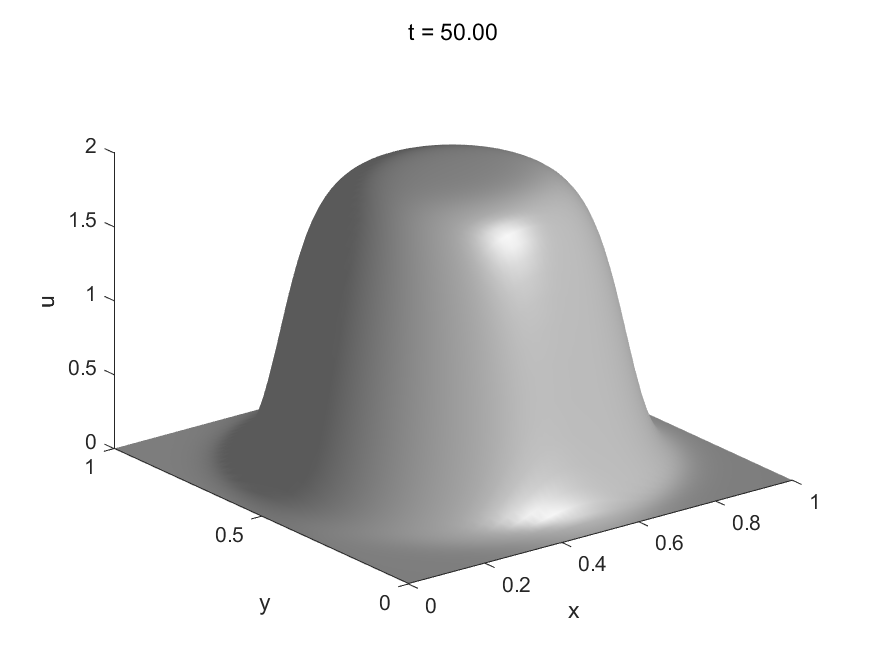}
	\includegraphics[width=0.32\textwidth]{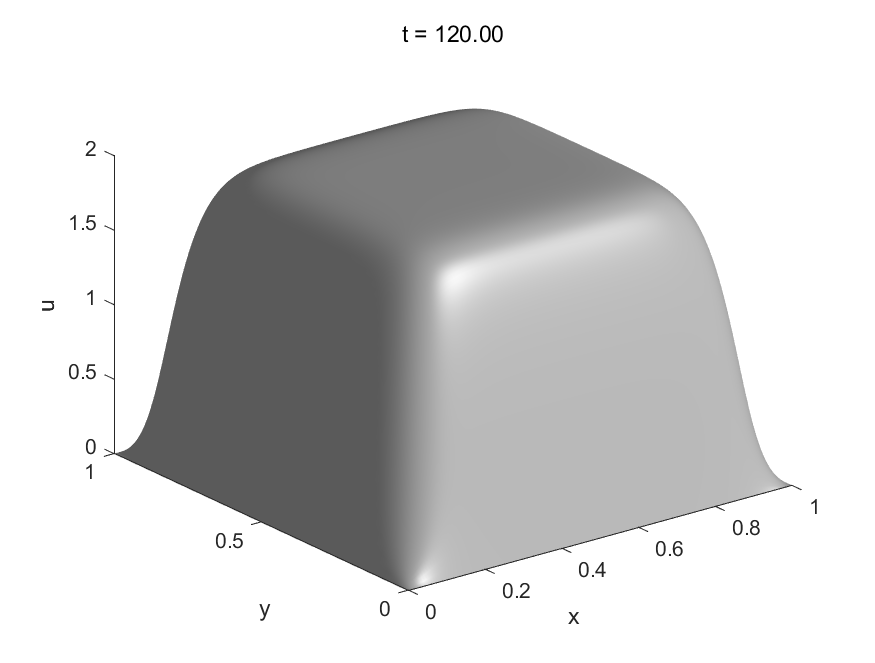}
	\includegraphics[width=0.32\textwidth]{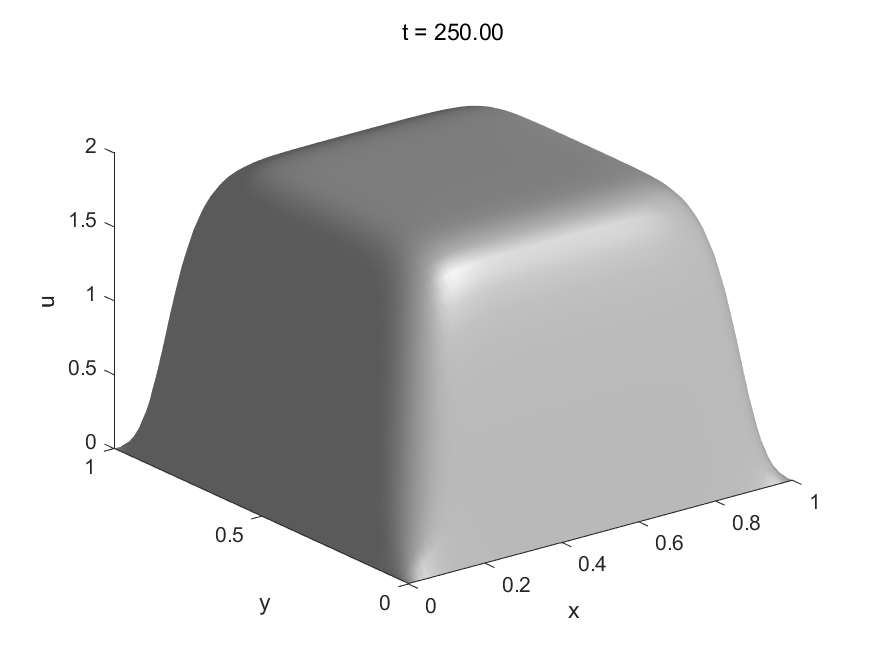}
	\caption{Endemic phenomenon of the average concentration of bacteria $u$.}\label{fig:u+}
\end{figure}

\begin{figure}[htbp]\small
	\centering
	\includegraphics[width=0.32\textwidth]{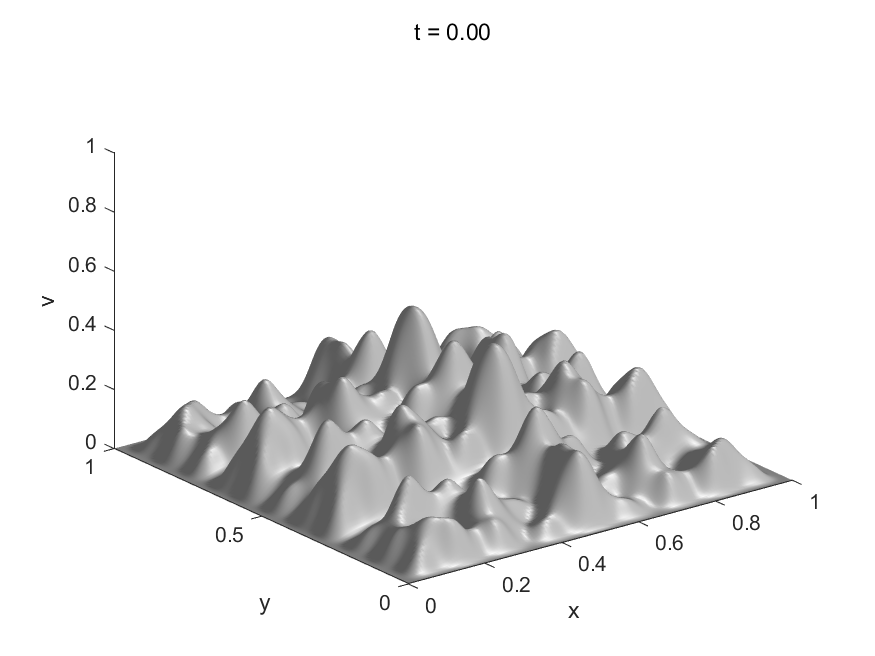}
	\includegraphics[width=0.32\textwidth]{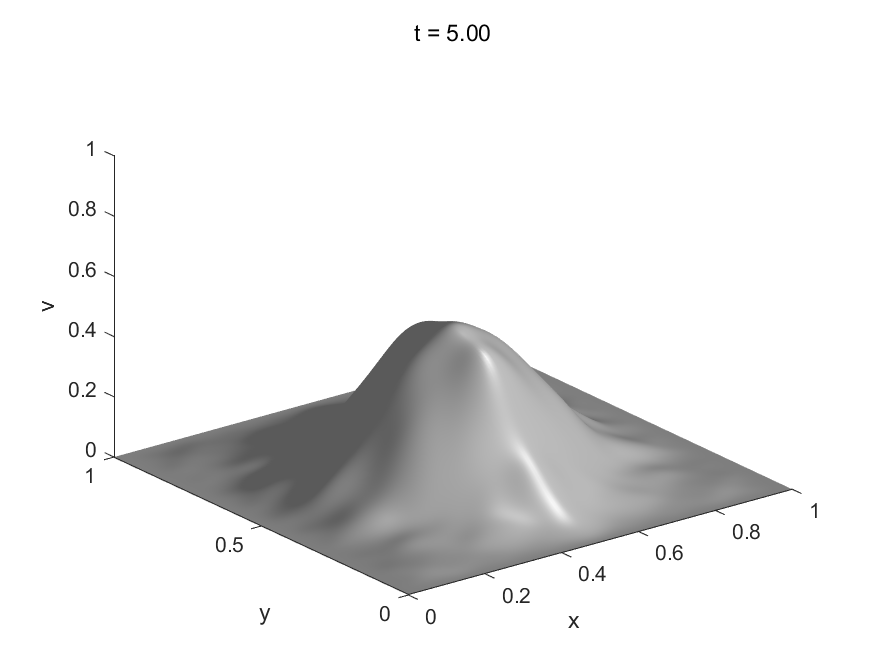}
	\includegraphics[width=0.32\textwidth]{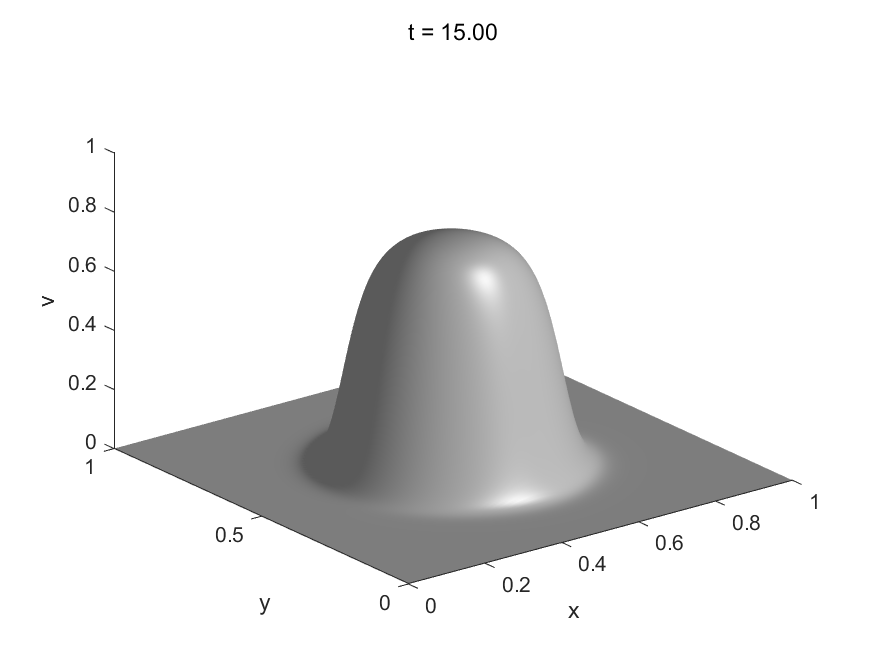}
	\includegraphics[width=0.32\textwidth]{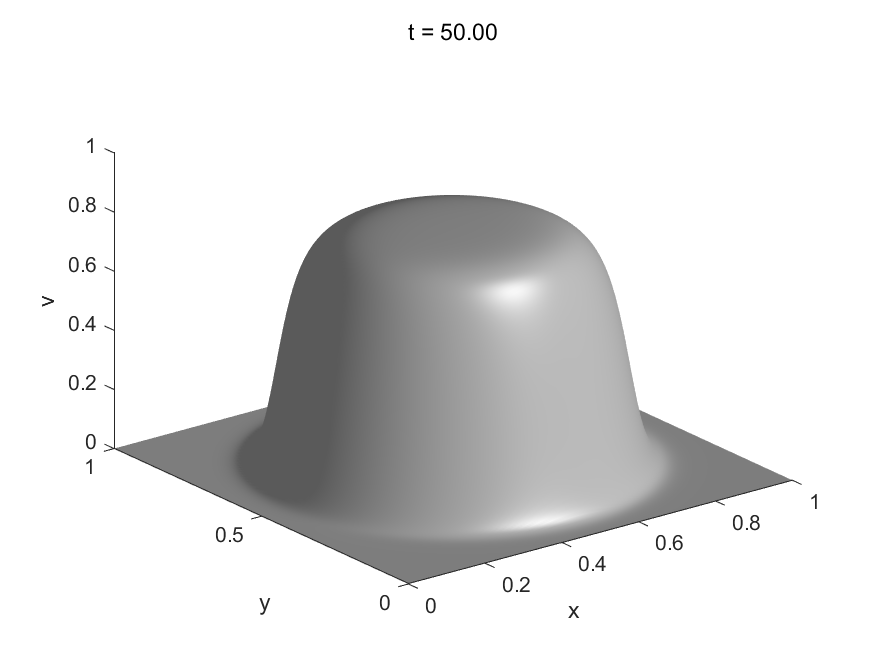}
	\includegraphics[width=0.32\textwidth]{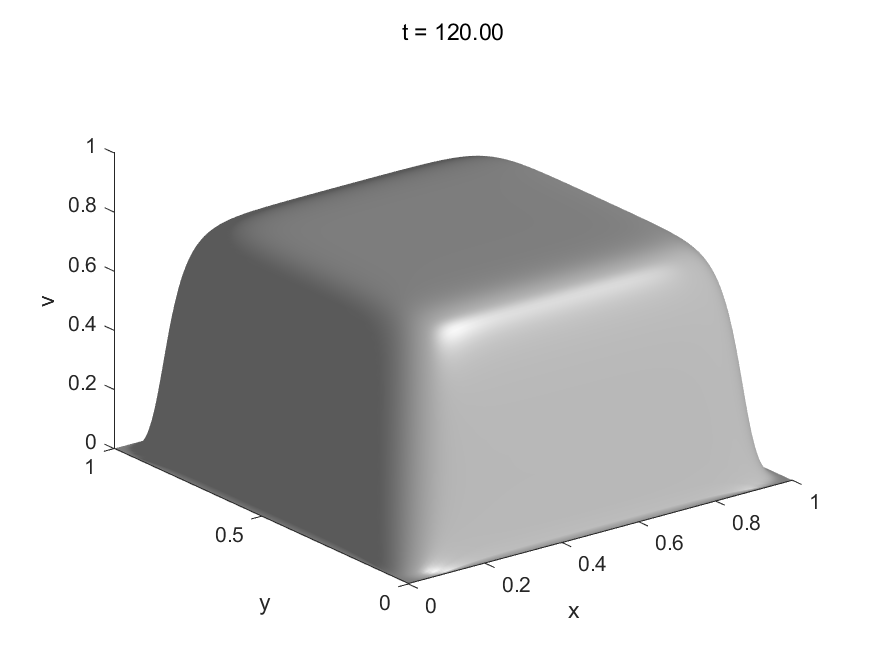}
	\includegraphics[width=0.32\textwidth]{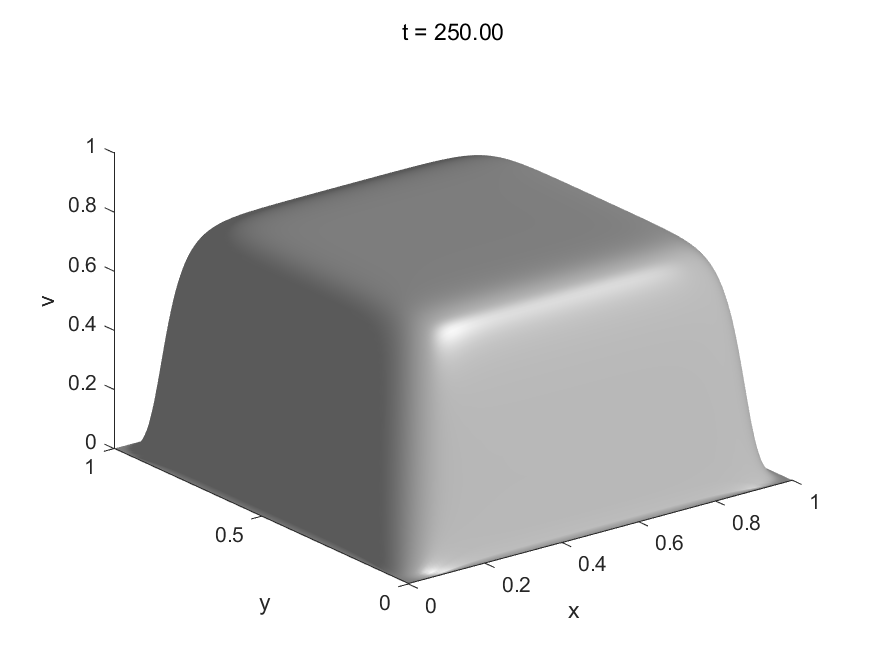}
	\caption{Endemic phenomenon of the average concentration of infective people $v$.}\label{fig:v+}
\end{figure}

In this test, we also take $\tau = 0.1$ and $M=160$. The evolutions of the bacterial concentration $u$ and the  concentration of infective people $v$ are respectively depicted in Figures \ref{fig:u+}--\ref{fig:v+}. 
For this chosen initial data, the bacterial concentration $u$ is attracted to a non-zero equilibrium state.
Ecologically, in this situation, it will ultimately lead to endemic diseases.
From Figures \ref{fig:u+}--\ref{fig:v+}, the solutions converge rapidly to a symmetric bell-shaped curve and enlarge gradually over time. 
After a much longer time period do the solutions approach a stable steady state (see $t=120$ and $t=250$ in Figures \ref{fig:u+}--\ref{fig:v+}). 
This phenomenon shows that a relatively  large initial value can indeed induce the spreading of bacteria, and as Figure \ref{fig:u+} shows the process of bacterial spreading finally reach to a stable steady state.

\begin{example}\label{exam:compact}
	In this example, we choose different function $g(u)$ and parameters $d_1$, $a_{11}$, $a_{12}$, $d_2$, $a_{22}$ as follows:
	\begin{enumerate}[(a)]
		\item $g(u)=u/(1+u)$, $d_1=0.001$, $a_{11}=1.0$, $a_{12}=2.5$, $d_2=0.0001$, $a_{22}=1.0$;
		\item $g(u)=u^2/(1+u^2)$, $d_1=0.01$, $a_{11}=1.0$, $a_{12}=2.5$, $d_2=0.0001$, $a_{22}=1.0$;
		\item $g(u)=u^2/(1+u^2)$, $d_1=0.001$, $a_{11}=1.0$, $a_{12}=2.5$, $d_2=0.001$, $a_{22}=1.0$;
		\item $g(u)=u^2/(1+u^2)$, $d_1=0.001$, $a_{11}=1.0$, $a_{12}=1.0$, $d_2=0.0001$, $a_{22}=1.0$.
	\end{enumerate}
The initial condition pairs in Example \ref{exam:noise}--\ref{exam:endemic}  are respectively denoted as $(u_0^-,v_0^-)$ and $(u_0^+,v_0^+)$.
\end{example}

\begin{figure}
	\centering
	\includegraphics[width=0.32\linewidth]{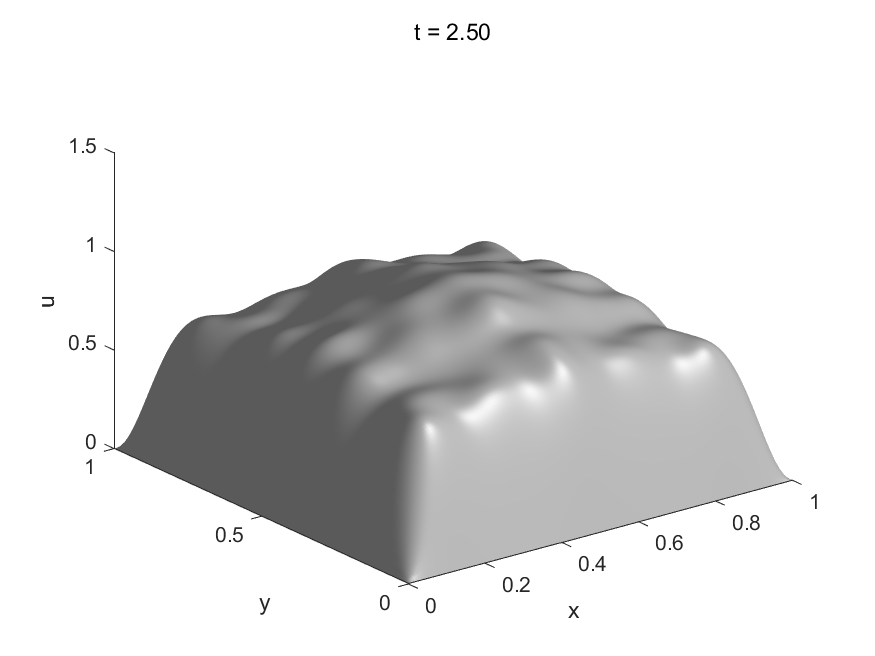}
	\includegraphics[width=0.32\linewidth]{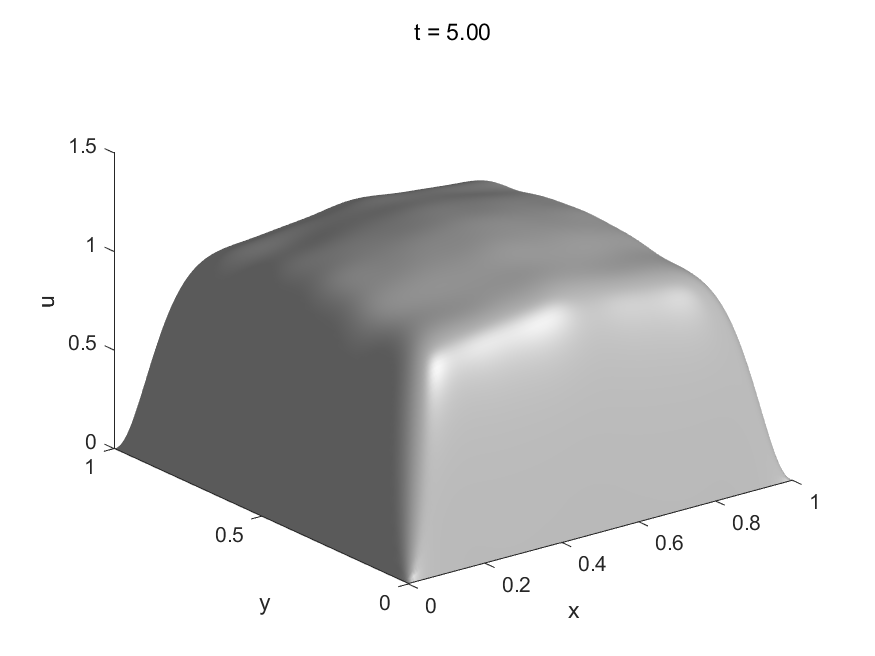}
	\includegraphics[width=0.32\linewidth]{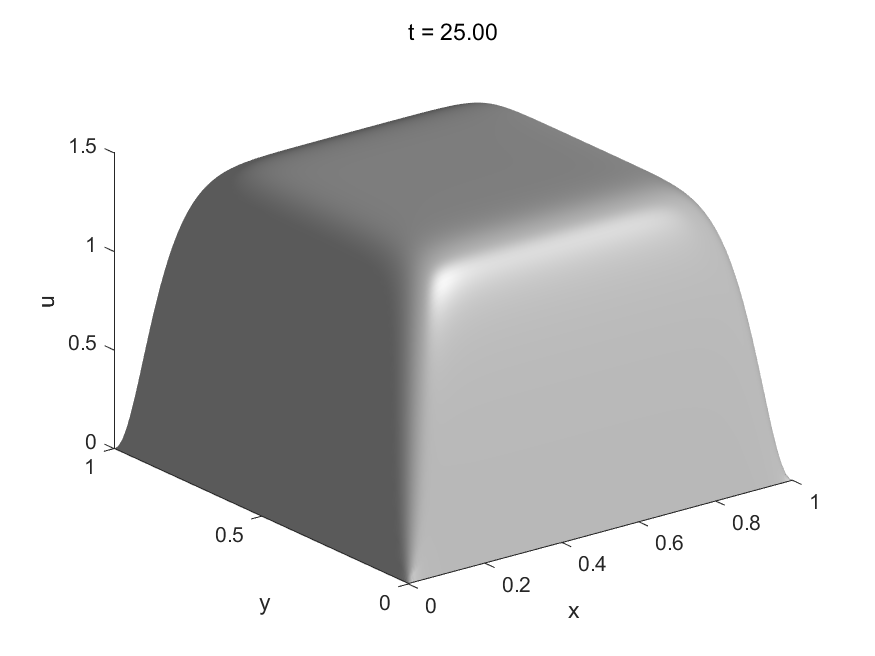}
	\includegraphics[width=0.32\linewidth]{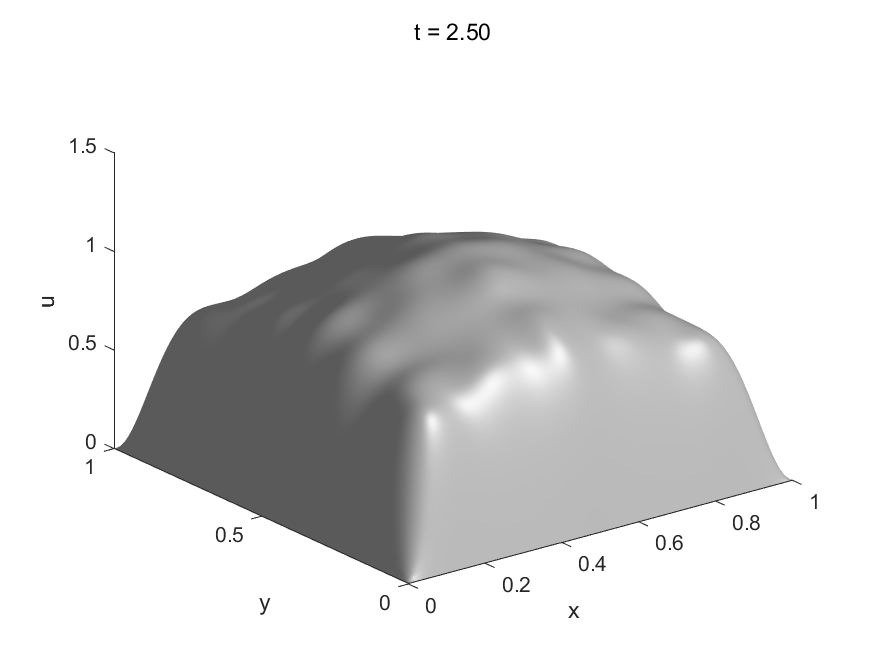}
	\includegraphics[width=0.32\linewidth]{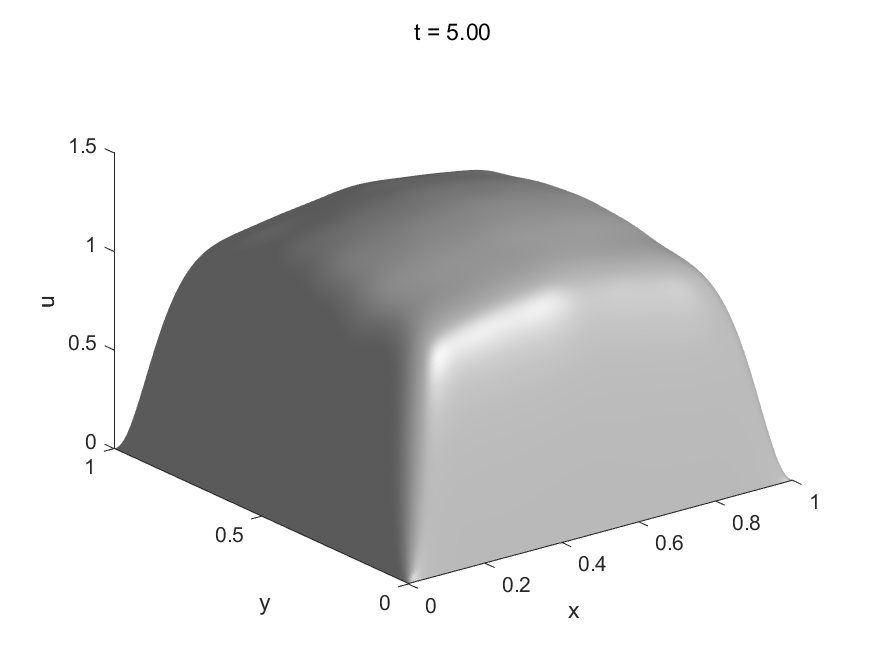}
	\includegraphics[width=0.32\linewidth]{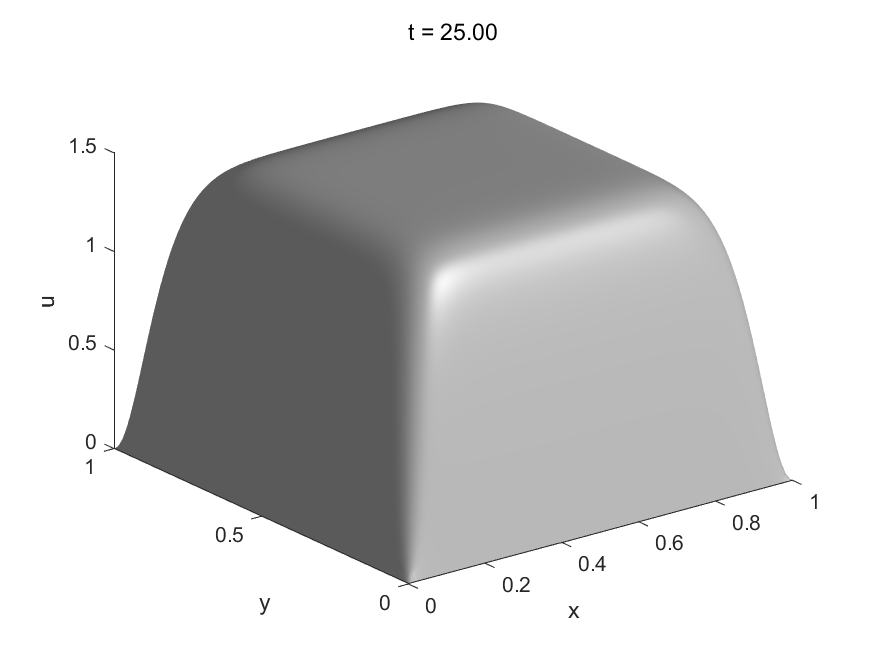}
	\caption{Evolution of $u$ with parameters (a) in Example \ref{exam:compact} and initial conditions $(u_0^-,v_0^-)$ (up) and $(u_0^+,v_0^+)$ (down).}
	\label{fig:dg}
\end{figure}

\begin{figure}
	\centering
	\includegraphics[width=0.32\linewidth]{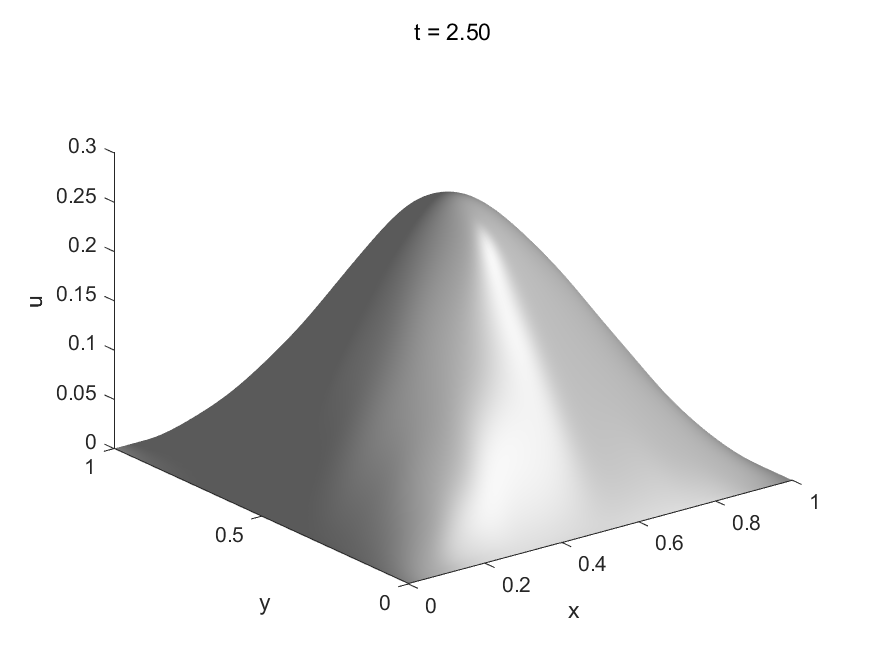}
	\includegraphics[width=0.32\linewidth]{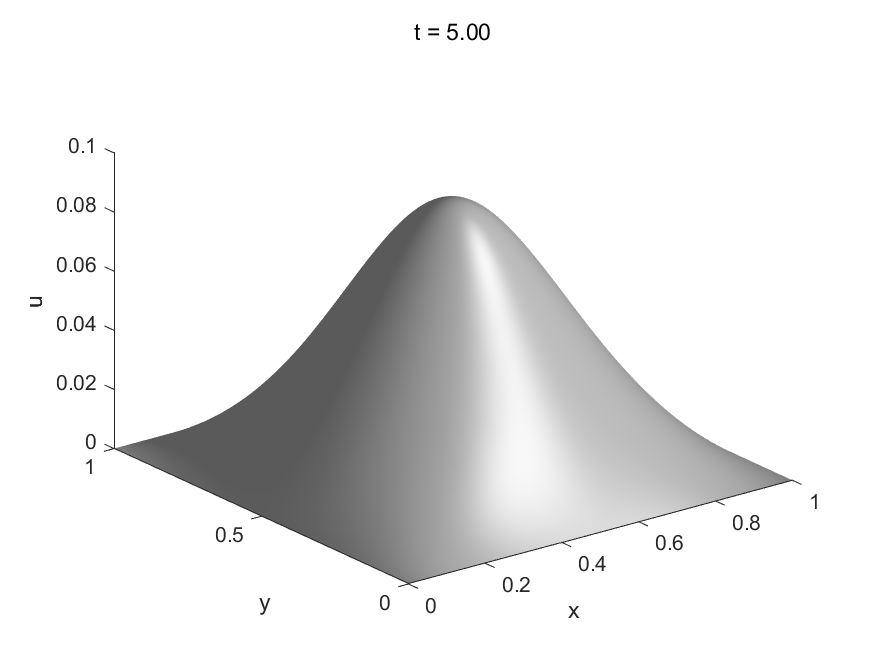}
	\includegraphics[width=0.32\linewidth]{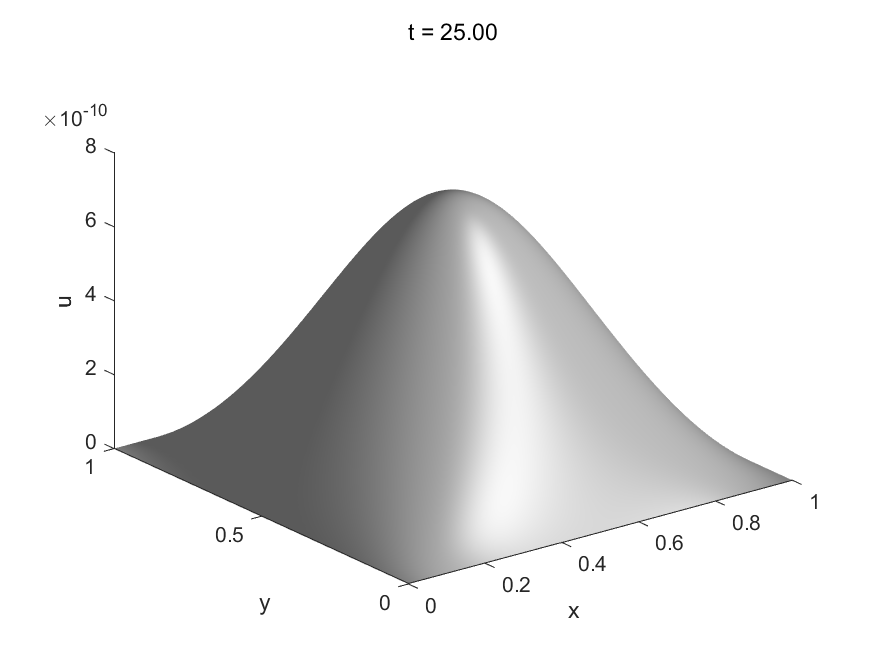}
	\includegraphics[width=0.32\linewidth]{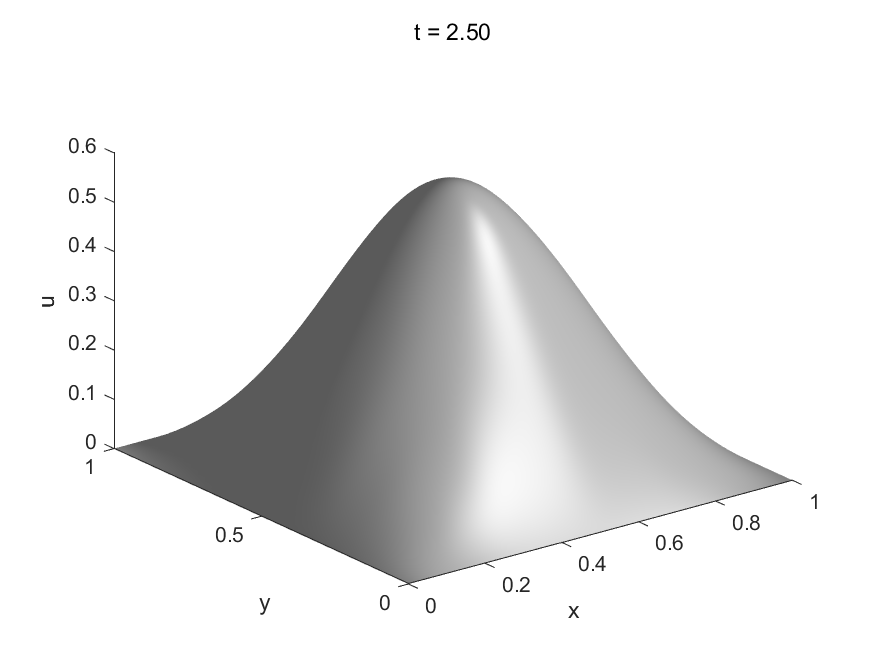}
	\includegraphics[width=0.32\linewidth]{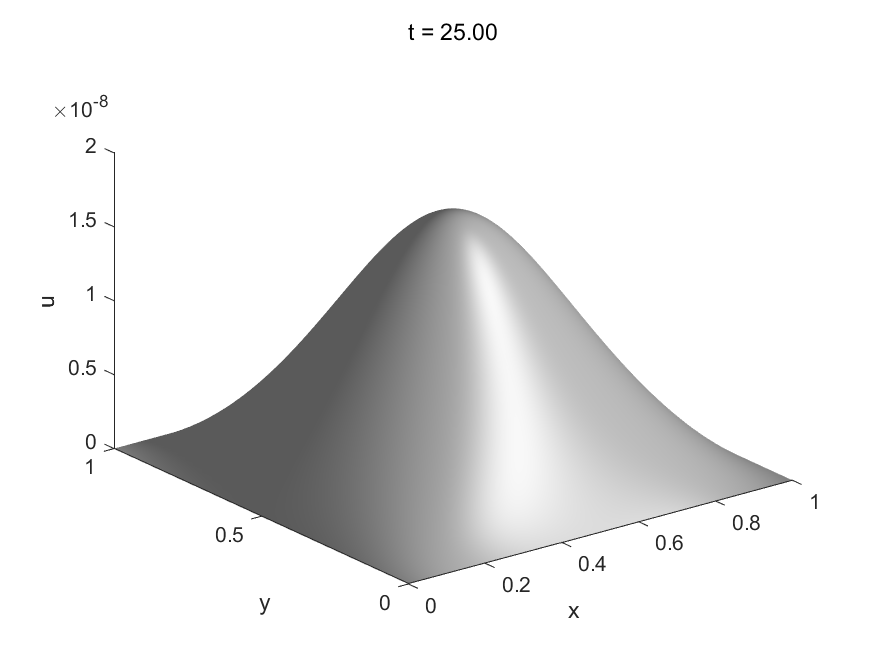}
	\includegraphics[width=0.32\linewidth]{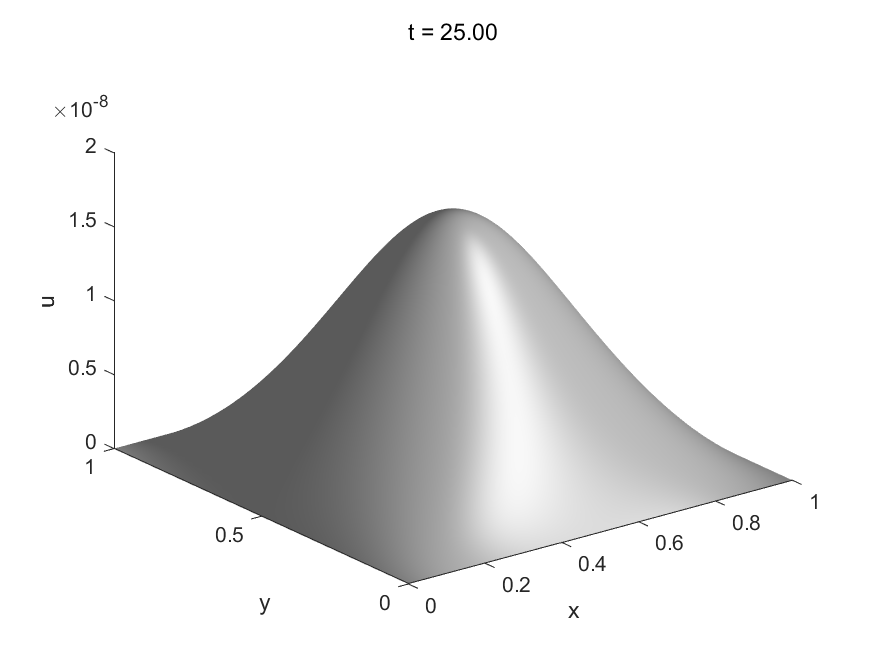}
	\caption{Evolution of $u$ with parameters (b) in Example \ref{exam:compact} and initial conditions $(u_0^-,v_0^-)$ (up) and $(u_0^+,v_0^+)$ (down).}
	\label{fig:dd1}
\end{figure}

\begin{figure}
\centering
\includegraphics[width=0.32\linewidth]{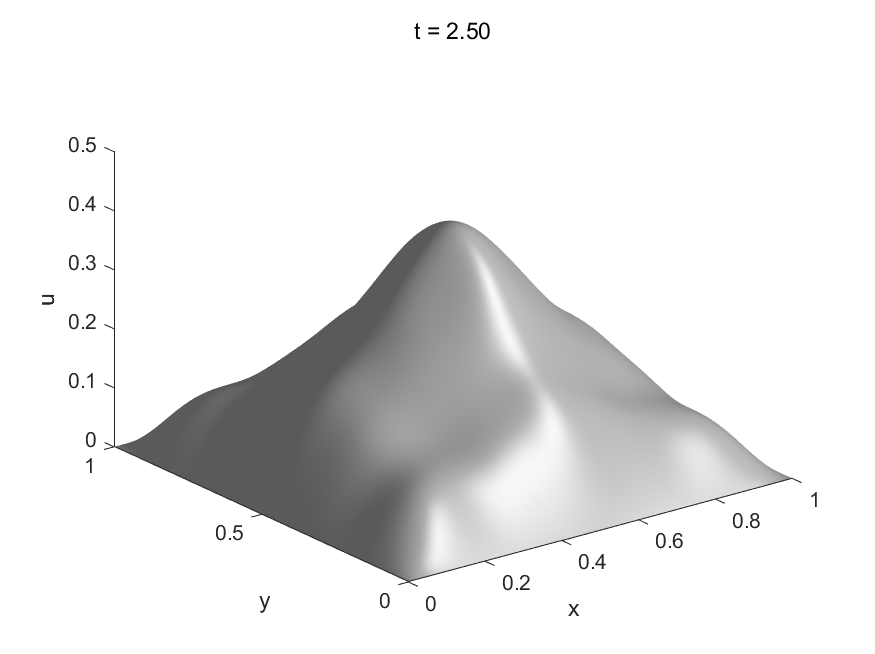}
\includegraphics[width=0.32\linewidth]{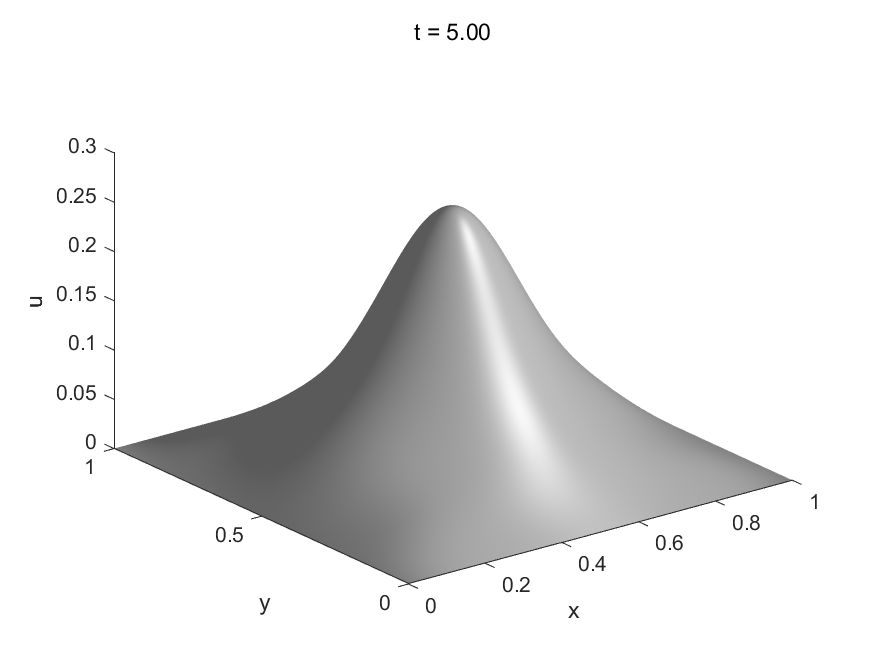}
\includegraphics[width=0.32\linewidth]{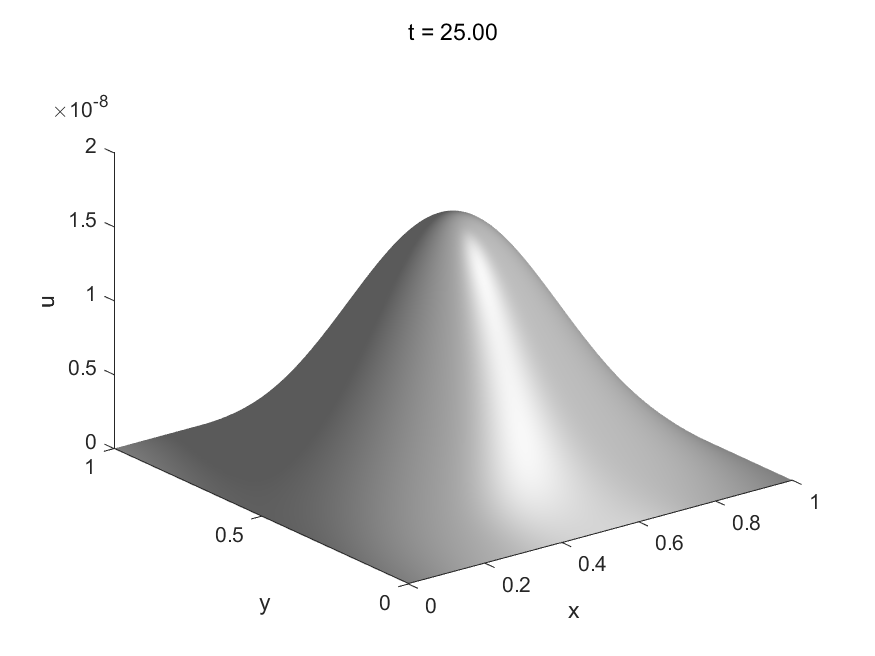}
\includegraphics[width=0.32\linewidth]{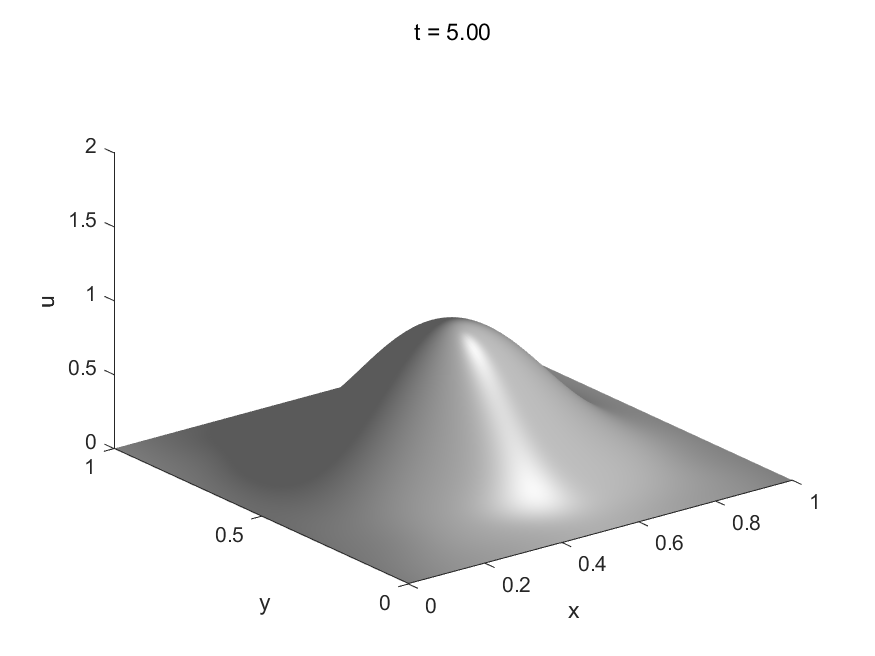}
\includegraphics[width=0.32\linewidth]{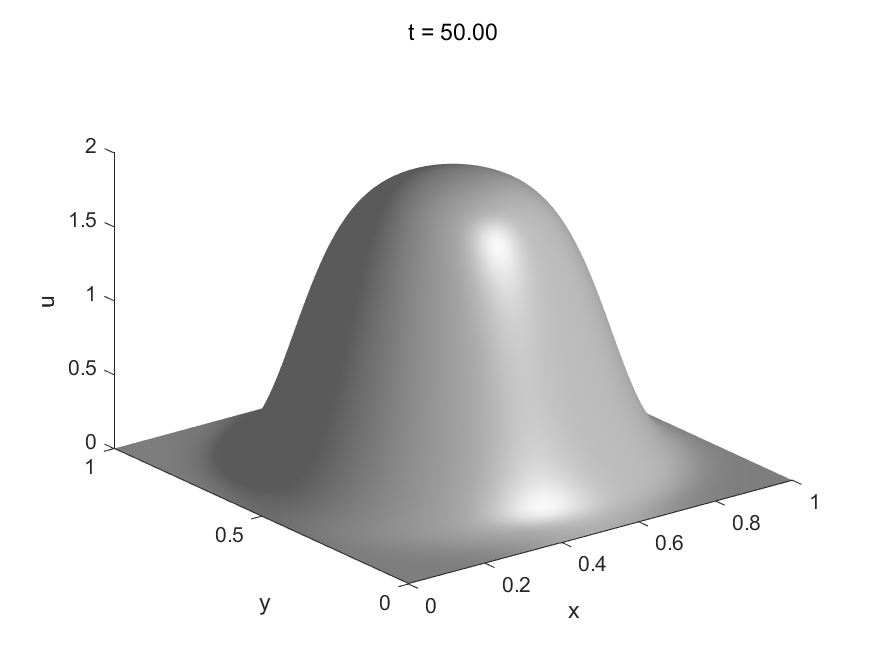}
\includegraphics[width=0.32\linewidth]{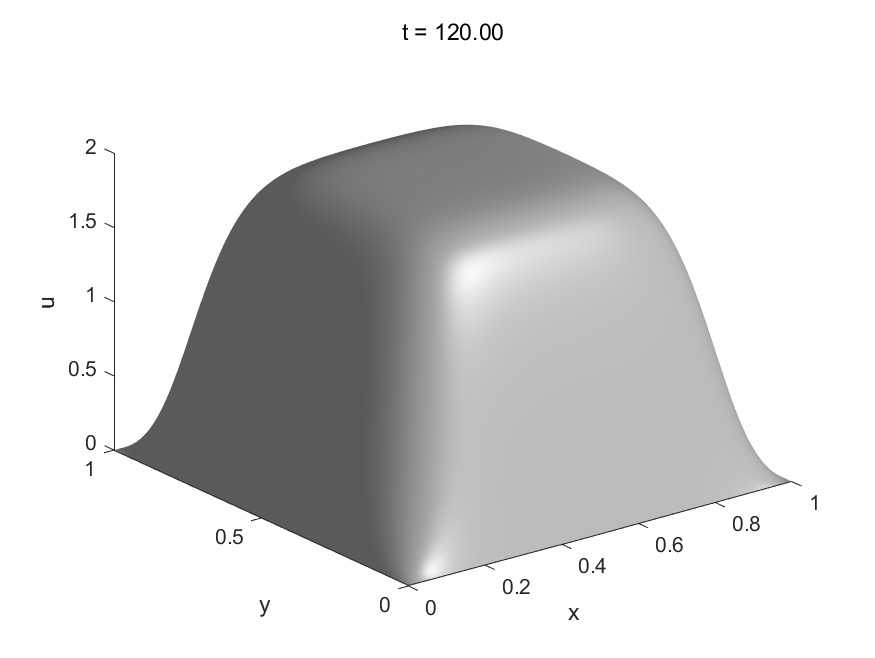}
\caption{Evolution of $u$ with parameters (c) in Example \ref{exam:compact} and initial conditions $(u_0^-,v_0^-)$ (up) and $(u_0^+,v_0^+)$ (down).}
\label{fig:dd2}
\end{figure}

\begin{figure}
	\centering
	\includegraphics[width=0.32\linewidth]{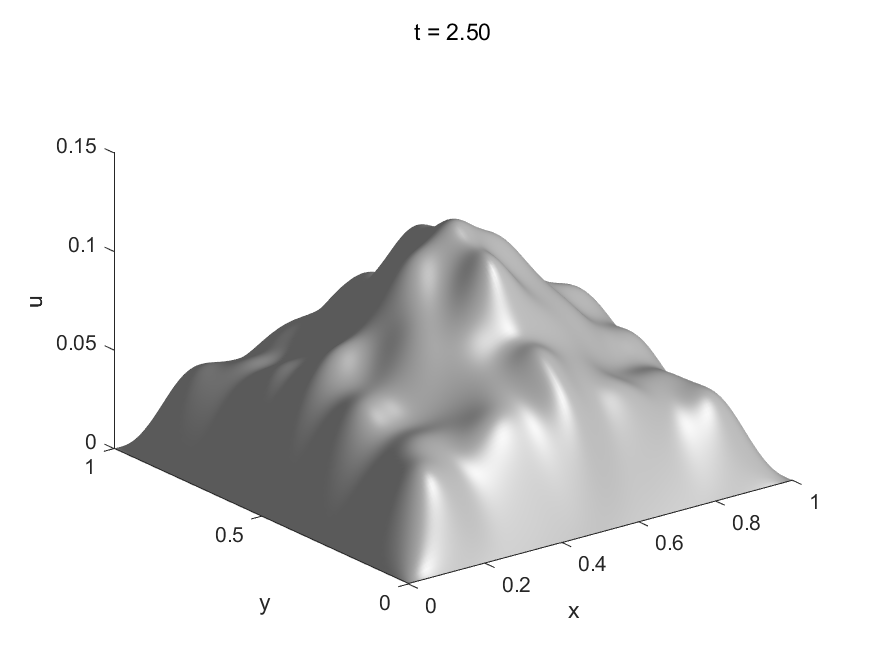}
	\includegraphics[width=0.32\linewidth]{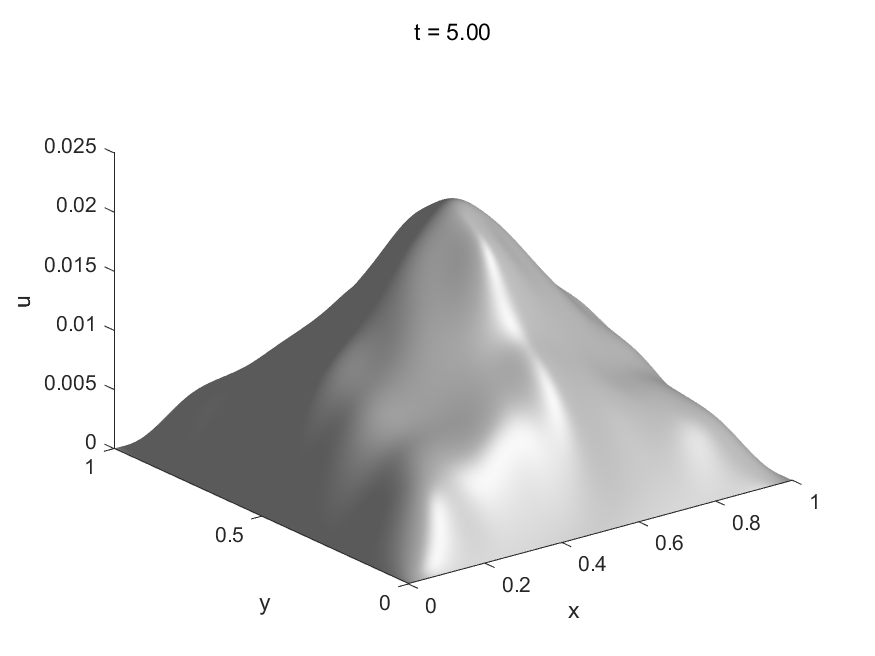}
	\includegraphics[width=0.32\linewidth]{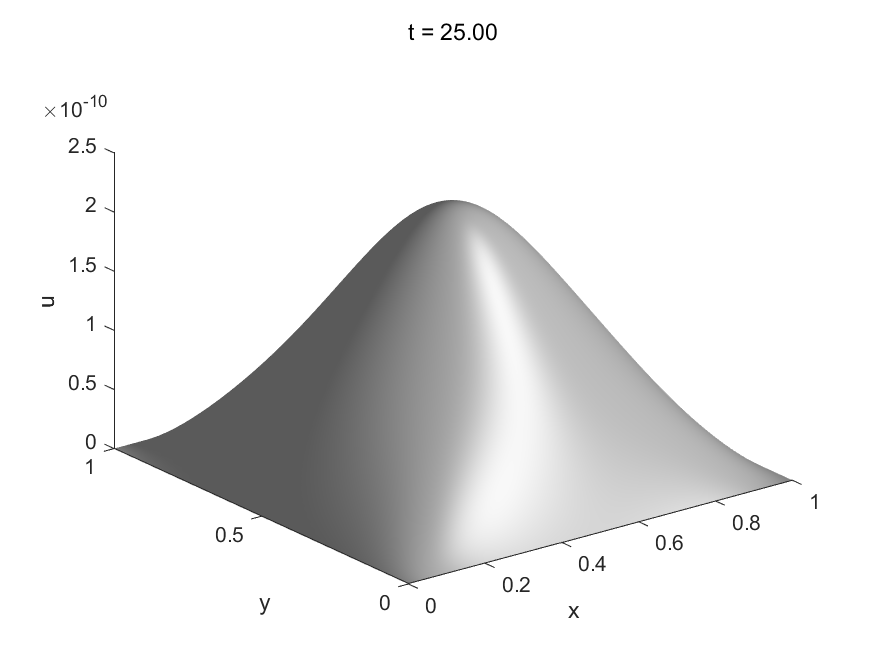}
	\includegraphics[width=0.32\linewidth]{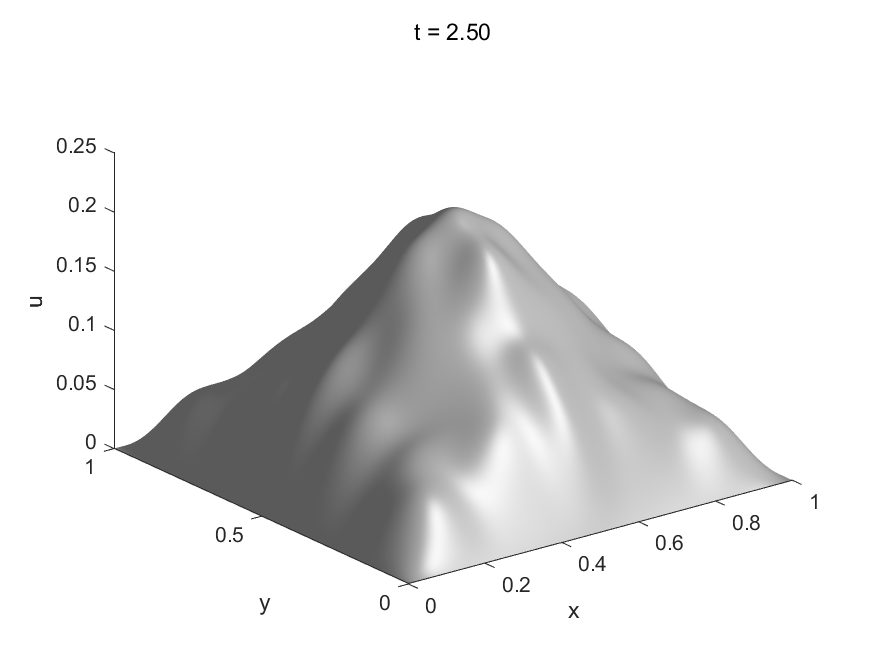}
	\includegraphics[width=0.32\linewidth]{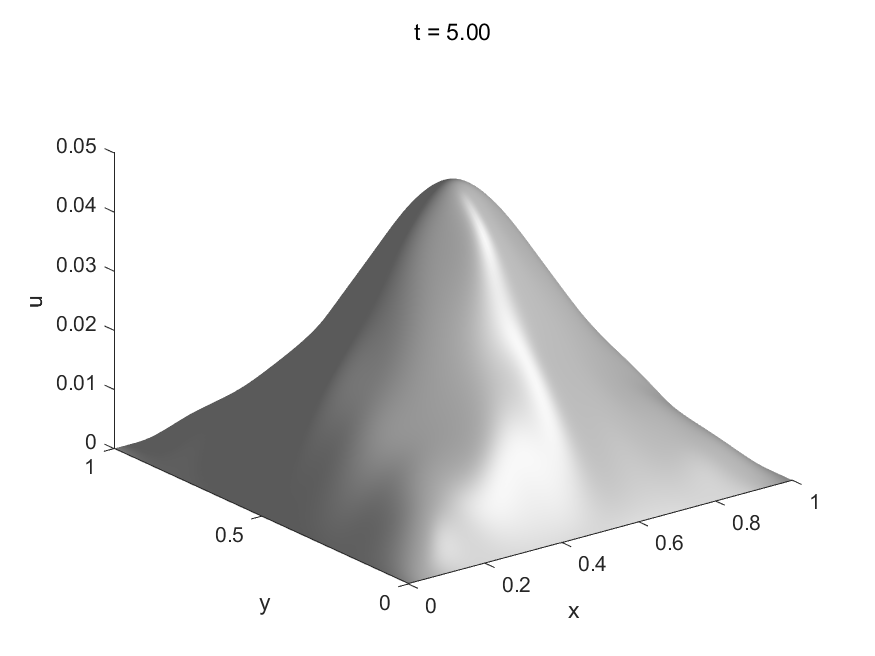}
	\includegraphics[width=0.32\linewidth]{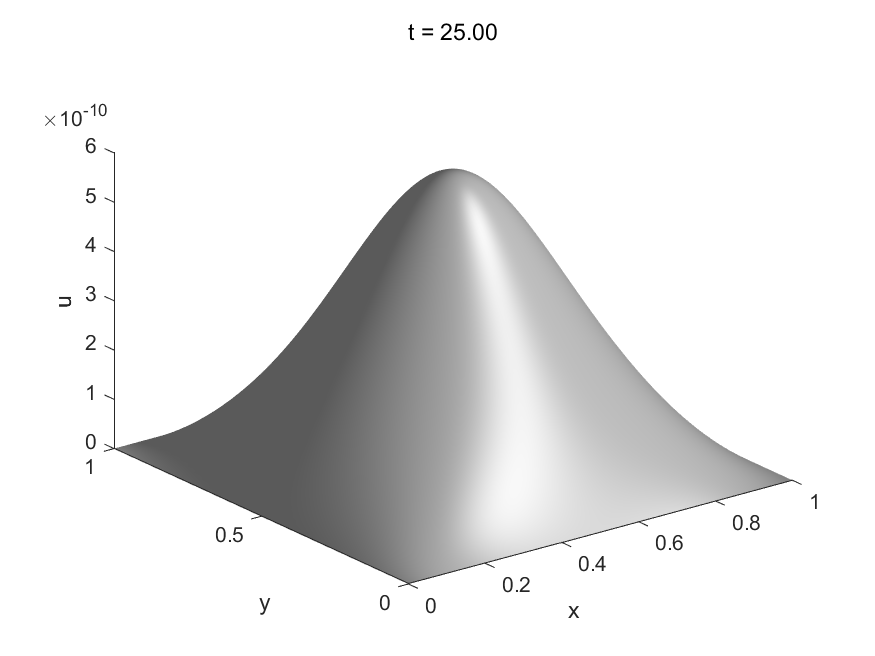}
	\caption{Evolution of $u$ with parameters (d) in Example \ref{exam:compact} and initial conditions $(u_0^-,v_0^-)$ (up) and $(u_0^+,v_0^+)$ (down).}
	\label{fig:da}
\end{figure}

As the average concentration of infective people $v$ has a similar evolution as the average concentration of bacteria $u$, we only show the evolutions of  $u$ with different parameters and initial conditions in Figures \ref{fig:dg}--\ref{fig:da}, where we still take $\tau = 0.1$ and $M=160$.  The following phenomena can be observed:
\begin{enumerate}[(i)]
	\item Compared with those parameters in Examples \ref{exam:noise}--\ref{exam:endemic}, a different infection function $g(u)$ is chosen in case (a).
	We observe from Figure \ref{fig:dg}  that the solutions converge to the same non-zero equilibrium state with different initial conditions $(u_0^-,v_0^-)$ and $(u_0^+,v_0^+)$. However, the bacteria with initial conditions $(u_0^-,v_0^-)$ in Example \ref{exam:noise} eventually becomes extinct. Moreover,
 the non-zero equilibrium states in Figures \ref{fig:u+} and \ref{fig:dg} are also different with different function $g(u)$.
	\item Compared with the diffusion coefficient $d_1$ in Examples \ref{exam:noise}--\ref{exam:endemic}, a larger parameter $d_1$ is chosen in case (b). Theoretically, with a larger value of $d_1$, the dissipation process occurs more rapidly. In this numerical simulation,  the solution $u$ is seen less than $0.1$ at $t = 5$ as shown in Figure \ref{fig:dd1}. Conversely,  Figure \ref{fig:u-} shows that, with the same initial conditions $(u_0^-,v_0^-)$, the solution is greater than $0.3$ at the same time, which is in accordance with the theoretical expectations. Furthermore, Figure \ref{fig:dd1} illustrates that in this case, the bacterial population extinction occurs even with larger initial conditions $(u_0^+,v_0^+)$.	Similarly, for the given parameters in case (c), Figure \ref{fig:dd2} demonstrates that a larger diffusion coefficient $d_2$ also results in a faster dissipation, but in this moment, the solution  will converge to a non-zero equilibrium state for larger initial conditions $(u_0^+,v_0^+)$.
	\item For the given parameter in case (d), we choose a smaller coefficient $a_{12}$ compared to Examples \ref{exam:noise}--\ref{exam:endemic}. Since  $a_{12}v$ represents the natural transmission rate. Thus, a smaller value of $a_{12}$ indicates a lower natural transmission rate, which in turn inhibits the spread of bacteria. From Figure \ref{fig:da}, it is evident that the solution $u$ with  initial conditions $(u_0^-,v_0^-)$ remains below $0.0025$ at $t = 5$. This observation aligns with the aforementioned point. Similarly, under this scenario, the solution $u$ with larger initial conditions $(u_0^+,v_0^+)$ will tend to converge towards zero, and this is totally different from the result in Example \ref{exam:endemic}, see Figure \ref{fig:u+} for details.
	
\end{enumerate} 


%%%%%%%%%%%%%%%%%%%%%%%%%%%%%%%%%%%%%%%%%%%%%%%%%%%%%%%%%%%%%%%%%%%%%%%%%%%%%%%%%
\section{Conclusions}
In this paper, two linearized, decoupled compact finite difference methods named CN-CFD  (see \eqref{bactre:cncfd}--\eqref{bactre:cncfd:ic}) and CN-ADI-CFD  (see \eqref{ADI:cfd}--\eqref{ADI:cncfd:ic}) for the bacterial system \eqref{bact:eq} are presented. By  using the temporal-spatial error splitting technique and the discrete energy method, 
we demonstrate that both methods are unconditionally stable,  and exhibit second-order temporal and fourth-order spatial accuracy with respect to the discrete maximum-norm, which are also verified by ample numerical experiments. High efficiency of the CN-ADI-CFD scheme \eqref{ADI:cfd}--\eqref{ADI:cncfd:ic} is also shown in the numerical tests.
Furthermore, we apply the efficient CN-ADI-CFD method to simulate the asymptotic equilibrium state of extinction or formation of endemic phenomenon reflected by  the average concentration of the bacteria $u$ and the  infective people $v$. 
%%%%%%%%%%%%%%%%%%%%%%%%%%%%%%%%%%%%%%%%%%%%%%%%%%%%%%%%%%%%%%%
\section*{Declarations}
\begin{itemize}
	\item \textbf{Funding}~   This work is supported in part by funds from the National Natural Science Foundation of China Nos. 11971482 and 12131014, the Fundamental Research Funds for the Central Universities Nos. 202264006 and 202261099, and the OUC Scientific Research Program for Young Talented Professionals.
	\item  \textbf{Data availability}~  Enquiries about data availability should be directed to the authors.
	\item \textbf{Conflict of interest}~  The authors declare that they have no competing interests.
	%	\item \textbf{Author contributions} Hongfei Fu: Conceptualization, Funding acquisition, Methodology, Supervision, Writing -- review editing. Bingyin Zhang: Formal analysis, Investigation, Methodology, Visualization, Writing -- original draft. Xiangcheng Zheng: Formal
	%	analysis, Funding acquisition, Investigation, Methodology, Writing -- review editing.  
\end{itemize}

%\section*{CRediT authorship contribution statement}
%\textbf{Jie Xu}: Methodology, Formal analysis, Software, Writing- Original draft.
%\textbf{Shusen Xie}: Methodology, Supervision, Writing- Reviewing and Editing.
%\textbf{Hongfei Fu}: Conceptualization, Supervision, Writing- Reviewing and Editing, Funding acquisition.
%
%\section*{Declaration of Competing Interest}
%The authors declare that they have no competing interests.
%
%\section*{Availability of supporting data}  Enquiries about data availability should be directed to the authors.
%\section*{Acknowledgements}
%%The authors would like to express their most sincere thanks to the referees for their very helpful comments
%%and suggestions, which greatly improved the quality of this paper.
%This work was supported in part by the National Natural Science Foundation of China (Nos. 11971482,12131014), 
%by the Fundamental Research Funds for the Central Universities (Nos. 202264006, 202261099) and by the  OUC Scientific Research Program for Young Talented Professionals.

\end{document}